\documentclass[11pt]{article}
\usepackage{amssymb}
\usepackage{geometry}
\usepackage{graphicx}
\usepackage{amsmath}
\usepackage{amssymb}
\usepackage{amsthm}
\usepackage{amsfonts}
\usepackage{enumerate}

\usepackage{bbm}
\usepackage{epsfig}
\usepackage{color}
\usepackage{setspace}
\usepackage{rotating}
\usepackage{array}
\usepackage{calc}

\usepackage{tcolorbox}
\usepackage{multirow}
\usepackage{multimedia}
\usepackage{epstopdf}
\AppendGraphicsExtensions{.tif}

\usepackage{float}
\usepackage{indentfirst}
\usepackage{multicol}
\usepackage{mathrsfs}
\usepackage[normalem]{ulem}
\usepackage{caption}

\newlength{\mylength}

\numberwithin{equation}{section}

\newtheorem{theorem}{Theorem}[section]

\newtheorem{proposition}[theorem]{Proposition}

\newtheorem{remark}{Remark}

\newcommand{\mA}{\mathcal{A}}

\newcommand{\mG}{\mathcal{G}}

\newcommand{\vu}{\vec{u}}
\newcommand{\vv}{\vec{v}}
\newcommand{\vx}{\vec{x}}
\newcommand{\vw}{\vec{w}}
\newcommand{\dt}{\Delta t}
 
\newcommand{\gb}{\kappa_{\rm ub}}

\renewcommand{\vec}[1]{\mbox{\boldmath$#1$}}

\definecolor{orange}{rgb}{1,0.5,0}
\definecolor{dgreen}{rgb}{0,0.5,0}

\newcommand{\du}{\, \mathrm{d}}

\newcommand{\myremarkend}{$\clubsuit$}
\newcommand{\mat}[1]{\mbox{\boldmath$#1$}}

\newcommand{\hvar}{h_0}
\newcommand{\etavar}{\eta_0}
\newcommand{\kapub}{\kappa_{\rm ub}}

\DeclareMathOperator{\sech}{sech}

\def\XXint#1#2#3{{\setbox0=\hbox{$#1{#2#3}{\int}$}
     \vcenter{\hbox{$#2#3$}}\kern-.5\wd0}}

\begin{document}

\title{Preconditioning and Linearly Implicit Time Integration for the Serre-Green-Naghdi Equations} 
\author{
Linwan Feng\thanks{Department of Mathematical Sciences, NJIT, Newark, NJ 07102, lf46@njit.edu.}, 
David Shirokoff\thanks{Department of Mathematical Sciences, NJIT, Newark, NJ 07102, david.g.shirokoff@njit.edu.}, 
Wooyoung Choi\thanks{ 
Department of Mathematical Sciences, NJIT, Newark, NJ 07102, wychoi@njit.edu}
}

\maketitle

\begin{abstract}
The treatment of the differential PDE constraint poses a key challenge in computing the numerical solution of the Serre-Green-Naghdi (SGN) equations.  In this work, we introduce a constant coefficient preconditioner for the SGN constraint operator and prove rigorous bounds on the preconditioned conditioning number.  The conditioning bounds incorporate the effects of bathymetry in two dimensions, are quasi-optimal within a class of constant coefficient operators, highlight fundamental scalings for a loss of conditioning, and ensure mesh independent performance for iterative Krylov methods.  

Utilizing the conditioning bounds, we devise and test two time integration strategies for solving the full SGN equations.  The first class combines classical explicit time integration schemes (4th order Runge-Kutta and 2nd--4th order Adams-Bashforth) with the new preconditioner. The second is a linearly implicit scheme where the differential constraint is split into a constant coefficient implicit part and remaining (stiff) explicit part.  The linearly implicit methods require a single linear solve of a constant coefficient operator at each time step.  We provide a host of computational experiments that validate the robustness of the preconditioners, as well as full solutions of the SGN equations including solitary waves traveling over an underwater shelf (in 1d) and a circular bump (in 2d).  

\end{abstract}
\medskip\noindent

    

\section{Introduction}\label{Sec:intro}
The shallow water equations (SWEs) provide an effective model for simulating coastal waves and large-scale fluid circulation in ocean basins. This includes important physical phenomena such as tides, storm surges, and the effect of extreme weather, such as hurricanes, on wave generation.  

From a modeling standpoint, the SWEs arise from an asymptotic expansion of the (depth-averaged) Euler equations truncated to the leading $O(1)$ term in $\beta \equiv \bar h / \lambda \ll 1$. Here $\beta$ is a small parameter measuring the ratio of the mean water depth $\bar h$ to the characteristic wavelength $\lambda$.  Carrying out the asymptotic expansion (in $\beta$) to the next order results in dispersive corrections to the SWEs. Physically, the incorporation of dispersion to the SWEs plays an important role in the propagation of long waves at larger (but still small) $\beta$. In particular, solitary waves may form through the balancing of dispersion and (steepening effects of) nonlinearity, while dispersion may also regularize the formation of shocks. 

The simplest models that incorporate dispersion arise under a weakly nonlinear assumption of small amplitude waves, i.e., $\alpha \equiv a/\bar{h} \ll 1$ with $a$ the amplitude of the wave. These include models for long waves such as the Koreteweg de-Vries (KdV) equation \cite{kdv95} and the Boussinesq equations \cite{Boussinesq1872} for uni- and bi-directional waves.  Several important phenomena such as extreme events observed in coastal regions, however, do not satisfy the weakly non-linear assumption, thus motivating the need for extending the models to include high-order nonlinear effects. To this end, Serre \cite{Serre1953}, Su and Gardner \cite{sg69} and Green and Naghdi \cite{gn76} derived a new class of long wave models with the leading-order dispersive terms that appear at $O(\beta^2)$ without any assumption on $\alpha$. Asymptotically, the Serre/SG/GN equations (often referred to as the SGN equations) modify the shallow water equations by incorporating dispersive terms truncated at an order of $\beta^2$. Therefore, the resulting SGN equations can also be considered as the \emph{dispersive shallow water equations} (DSWEs).

There are numerous asymptotic models closely related to the SGN equations. For instance, the SGN equations for the depth-averaged velocity can be formally extended to any arbitrary order in $\beta^2$, e.g., see the work by Wu and others \cite{wu99,mbl02,ms98,mat15,c22}. However, it has been shown that the models truncated to even-orders in $\beta^2$ (i.e., the even-order models) are unstable when the wavenumber $k$ is greater than some critical values, at which the growth rate becomes infinity \cite{mat16,c22}. Therefore, the even-order models must be avoided while the odd-order models (including the SGN equations) are useful for numerical studies. In fact, the dispersion relation for the depth-averaged DSWEs shows that the wave frequency is bounded as $k\to \infty$ so that there is a less severe constraint on the step size for any conditionally-stable time integration scheme.

Related to the SGN equations, Nwogu \cite{nwo93} introduced an alternative weakly nonlinear long wave model aimed at more closely approximating the linearized Euler equations dispersion relation (cf. the model solved in \cite{BonnetonChazelLannesMarcheTissier2011}). Following Nwogu, Wei {\it et al.} \cite{wei95} proposed the first-order DSWEs written in terms of the horizontal velocity at a vertical level (as opposed to a depth-averaged velocity). 
In particular, it was shown that the higher order DSWEs for the horizontal velocity at the bottom are also well-posed at any order of approximation \cite{c22}. 

There are, however, numerical advantages to solving the depth averaged SGN equations over the DSWE equations in terms of the horizontal bottom velocity. For the depth-averaged SGN equations, the mass conservation equation contains a first order spatial derivative (see the detailed discussion in \S\ref{Sec:DSWE}). On the other-hand, for other DSWEs, high-order spatial derivatives need to be evaluated accurately to avoid numerical instability \cite{c22}.

For work on the existence of solutions to the three-dimensional water wave equations we refer to  \cite{SamaniegoLannes2008}; a general summary of the derivation for different equivalent SGN models can be found in \cite{LannesBonneton2009,c22}. The physics of water waves and their derivation of their models have been extensively studied, e.g., \cite{LiHymanChoi2004, DutykhClamondMilewskiMitsotakis2013, MadsenMei1969, Boussinesq1872, c22}.

\subsection{Prior Numerical Work on the SGN Equations}\label{Subsec:PriorWork}
The dispersive terms in the SGN equation contain mixed space-time derivatives resulting in a differential constraint. This differential constraint poses a key computational challenge that can dominate the computational cost of an SGN solver.  In one dimension, the constraint can be handled in a straightforward fashion without significant cost, either via classical iterative methods, 
direct solvers, or integral equation approaches, for example, in \cite{DutykhClamondMilewskiMitsotakis2013,LiHymanChoi2004,JangSungPark2024}. 
In two dimensions the differential constraint results in a couple system for the velocity components with off diagonal blocks introduced by a varying topography. 

One approach to time-discretize the SGN equations is to split the equation into a hyperbolic step (which can be treated with hyperbolic solvers) and a step for the elliptic (differential) constraint.  These approaches include Strang splitting \cite{BonnetonChazelLannesMarcheTissier2011, LannesMarche2015, SamiiDawson2018}, and predictor-corrector type approaches \cite{MetayerGavrilyukHank2010, Parisot2019, GavrilyukShyue2023, NoelleParisotTscherpel2022}. A correct treatment of the boundary conditions \cite{NoelleParisotTscherpel2022} is then important when splitting methods are used in non-periodic geometries. Most of these aforementioned works focus on practical details of spatial discretizations (e.g., with finite element methods) and structure preservation (e.g., well-balanced and positive preserving schemes).  

Approaches that address the SGN elliptic constraint include the development of quadtree multigrid methods \cite{Popinet2015}.  An alternative approach, \cite{LannesMarche2015, DuranMarche2017} instead modifies the SGN equation---to an asymptotically consistent one in $\beta$---by replacing the elliptic constraint with a block diagonal constant coefficient operator. The constraint can then be solved via two time-independent decoupled scalar elliptic systems.  

In contrast to splitting methods or predictor-corrector methods, an alternative approach is to introduce a hyperbolic relaxation of the SGN \cite{FavrieGavrilyuk2017} where the differential constraint is replaced by a hyperbolic system.  The relaxation approach can then be solved using explicit structure preserving methods \cite{GuermondKeesPopovTovar2019} that leverage the extensive body of existing hyperbolic PDE solution techniques. The original work of Favrie and Gavrilyuk has since been extended to incorporate the effects of topography \cite{BustoDumbserEscalanteFavrieGavrilyuk2021, GuermondKeesPopovTovar2022b}; numerical solutions include high order Galerkin spatial discretizations \cite{BustoDumbserEscalanteFavrieGavrilyuk2021} and up to second order in space and time structure preserving methods (e.g., well-balanced, positivity preserving) properties \cite{GuermondKeesPopovTovar2022b}. Practical details regarding convex limiters \cite{GuermondKeesPopovTovar2022a} and schemes that preserve mass and energy \cite{RanochaRicchiuto2024} have also recently been devised.  Structure preservation can improve accuracy over long times, and provide a pathway for ensuring numerical stability. 

\subsection{Contribution of this Work and Organization of Paper}\label{Subsec:ContThisWork}
The differential constraint in the SGN requires solving linear systems involving an operator $\mG_{\eta,h}$ with variable coefficients in both space and time.  The focus of this paper is on introducing formulas for constant coefficients $(\sigma, \alpha)$ that define an operator $\mA = \sigma \mat{I} - \alpha \nabla (\nabla \cdot)$, together with quasi-optimal bounds on the spectrum of $\mA^{-1} \mG_{\eta,h}$.  The bounds then enable time-stepping discretizations for the SGN that avoid fully constructing and inverting $\mG_{\eta,h}$, namely:
\begin{enumerate}
    \item For time integration methods that require solving systems involving $\mG_{\eta,h}$, the operator $\mA$ may be used as a preconditioner in a conjugate gradient method.  Linear systems involving $\mG_{\eta,h}$ arise in (standard explicit) time-stepping discretizations of the SGN equations, and must be solved in large numbers, making fast solvers desirable.
    \item The operator $\mG_{\eta,h}$ may be treated in a linearly implicit multi-step method where only $\mA$ is treated implicitly. Such linearly implicit multistep methods are efficient as they require only one solves of $\mA$ per time-step.
\end{enumerate}

This paper is organized as follows. We first introduced the SGN equations (\S\ref{Sec:DSWE}) and background for the paper (\S\ref{sec:linoper_background}). We then present (\S\ref{Sec:BoundonG}) new formulas for the constant coefficient operator $\mA$ along with rigorous bounds on the conditioning number and generalized eigenvalues. Utilizing the bounds, \S\ref{Sec:TimeStepping} devises zero-stable linearly implicit multistep time integration schemes.  Numerical tests are then performed in (\S\ref{Sec:numerical_tests}) showcasing the effectiveness of the preconditioner $\mA$ in a suite of progressively challenging test cases. The full SGN equations are then solved in dimensions one and two, comparing the performance and accuracy of different time integration approaches. The Appendices provide proofs of the main theoretical results.

\section{The Serre-Green-Naghdi Equations}\label{Sec:DSWE}
In \S\ref{Subsec:Background} we review the derivation of the SGN equations and introduce in \S\ref{Subsec:DAEform} an alternative \emph{constraint} form of the equations (equivalent for smooth solutions).  Spatial discretizations of the constraint form will naturally lead to a differential algebraic equations (DAE) which will play a role in the development and stability analysis of the numerical schemes. 

\begin{figure}
	\centering 
	\includegraphics[width=.40\textwidth]{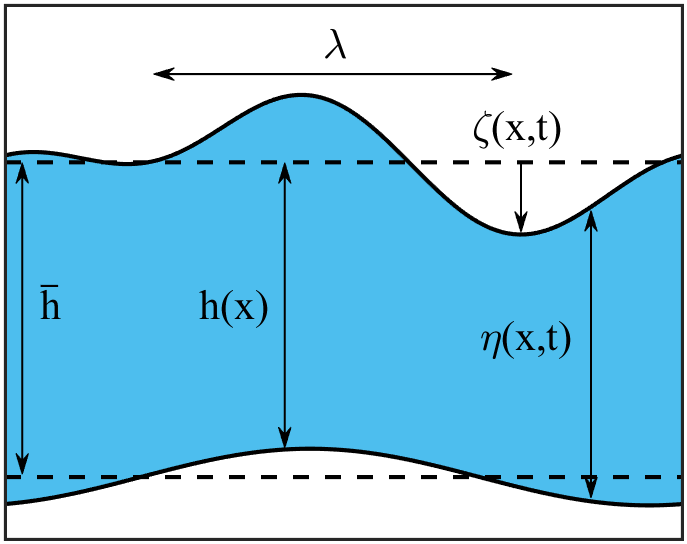}
	\caption{Domain (not to scale) with length scales and variables.}\label{Fig:Water} 
\end{figure}

\subsection{SGN Equations as Asymptotic Models of the Euler Equations}\label{Subsec:Background}
As with the shallow water equations (SWE), the SGN equations arise as asymptotic expansions of the depth averaged incompressible Euler equations.  In a domain $(\vec{x}, z) \in \mathbb{R}^{D + 1}$, where $D = 1$ or $D = 2$, the incompressible Euler equations for a homogeneous fluid layer of constant density $\rho$ take the form:
\begin{subequations}\label{Eq:Euler}
\begin{alignat}{2}\label{Eq:Euler1}
	\vw_t + \vw\cdot \nabla_3 \vw &= -\frac{1}{\rho} \nabla_3 p - \vec{g} \;, \\ \label{Eq:Euler2} 
	\nabla_3 \cdot \vw &= 0 \; .  
\end{alignat}
\end{subequations}
For notational purposes, we separate out the horizontal spatial coordinates $\vec{x} \in \Omega \subseteq \mathbb{R}^{D}$ from the vertical coordinate $z \in \mathbb{R}$.  In \eqref{Eq:Euler}, $\vw = (\vu, w) \in \mathbb{R}^{D+1}$ is the fluid velocity with $\vu \in \mathbb{R}^D$ being the horizontal component, $w$ the vertical component; $p(\vec{x},z,t)$ is the pressure; $\vec{g}$ is the gravitational force. Meanwhile $\nabla_3 = (\nabla, \partial_z)$ is the (full) gradient in $\mathbb{R}^{D+1}$, where $\nabla$ is the gradient in the variables $\vec{x}$. 

As shown in Su and Gardner \cite{sg69}, the SGN (or DSWEs) arise from: (I) integrating the Euler equations over the vertical column of fluid (i.e., $z$ coordinate); and (II) imposing boundary conditions at the free surface $\zeta(\vec{x},t)$ and the bottom bathymetry $h(\vec{x})$. The derivation results in evolution equations for the new dependent variables which are the surface elevation $\zeta(\vec{x},t)$, and the depth-averaged horizontal velocity $\bar{\vec u}(\vec{x},t)$ defined by 
\begin{align}\label{Eq:DepthAvV}
	\bar{\vec u}(\vec{x},t) = \frac{1}{\zeta +h}\int_{-h}^{\zeta}\vec{u}(\vec{x},z,t) \; \du z \; .
\end{align}
It is worth noting that the process of depth-averaging the Euler equations results in expressions of the form $\overline {\vec{u}\otimes\vec{u}}$ arising from the (depth-averaged)  nonlinear terms.  To obtain a closed system of PDEs for $\zeta, \bar{\vec{u}}$, the nonlinear term $\overline {\vec{u}\otimes\vec{u}}$ is approximated in terms of $\zeta$ and $\bar{\vec u}$ under the following assumptions
\begin{align*}
    \beta \equiv\frac{\bar{h}}{\lambda}\ll 1 \; , \quad\quad \zeta/\bar h = O( a/\bar h)=O(1) \;, \quad\quad \vert \bar \vu\vert / (g \bar h)^{1/2}=O(1) \; , \quad\quad \nabla h=O(\bar h/\lambda)=O(\beta) \;.
\end{align*}
Here $\bar{h}$ is the average water depth, while $a$ and $\lambda$ denote the characteristic amplitude and wavelength of the wave.  Typical values for $a$ and $\lambda$ may be estimated in the initial data $\vw_0$ (see Figure~\ref{Fig:Water}). The assumption $\beta \ll 1$ is referred to as the shallow water, or long wave assumption. Note that the assumptions allow for a ``large amplitude'' $a$ to be comparable in magnitude to the mean water depth $\bar h$. Likewise, the associate fluid motion $\vert \bar \vu\vert$ relative to $g\bar{h}$ does not need to be asymptotically small in $\beta$ and may be $O(1)$. 

With the asymptotic assumptions, one can then obtain a closed system for $\zeta$ and $\bar{\vec u}$ correct to $O(\beta^2)$. This yields the (dimensional form) of the Serre-Green-Naghdi (SGN) equations
\begin{subequations}\label{Eq:dswe}
    \begin{align}\label{Eq:dsweA}
        &\eta(\vec{x},t) =h(\vec{x})+\zeta(\vec{x},t)\;, \phantom{\left(\frac{1}{3}\right)}  \\ \label{Eq:dsweB}
        &\eta_t + \nabla\cdot\left(\eta\bar{\vec{u}}\right) =0\;, \\ \label{Eq:dsweC}
        &\bar{\vec{u}}_t+\bar{\vec{u}}\cdot\nabla \bar{\vec{u}}+g\nabla\zeta =
	\frac{1}{\eta}\nabla\left( \tfrac{1}{3}\eta^3 G +\tfrac{1}{2}\eta^2F\right) 
	-\left(\tfrac{1}{2}\eta G+F\right )\nabla h \;,
    \end{align}
\end{subequations}
where $G$ and $F$ are defined as
\begin{align}
	G=\nabla\cdot \bar \vu_t+
	\bar \vu\cdot\nabla\left(\nabla\cdot\bar \vu\right)-\left(\nabla\cdot \bar \vu\right)^2\;,\qquad\qquad
	F=\bar \vu_t\cdot\nabla h
	+\left(\bar \vu\cdot \nabla\right)^2h\;.
	\label{Eq:defGF}
\end{align}
Note that when the right-hand side of the SGN velocity equation is neglected, one recovers the (non-dispersive) shallow water equations.  We therefore, also refer to \eqref{Eq:dswe} as the \emph{dispersive shallow water equations}.  If an additional small-amplitude assumption of $\alpha=O(a/\bar h)\ll 1$ is adopted, the system \eqref{Eq:dswe} can be reduced to the weakly nonlinear Boussinesq equations \cite{Boussinesq1872}.

For the remainder of the paper, we non-dimensionalize all variables with respect to $g$ and $\bar h$, or, equivalently, set $g= 1$, $\bar h = O(1)$, and drop the bar on $\bar{\vu}$ and use $\vu$ for the depth-averaged velocity.  

\subsection{Constraint Form of the SGN Equations}\label{Subsec:DAEform}
We now recast \eqref{Eq:dswe} into two evolution equations plus a PDE constraint. For classical solutions with sufficient regularity, the resulting constraint form is equivalent to \eqref{Eq:dswe}.  

First, multiply the velocity equation in \eqref{Eq:dswe} by $\eta$ and group terms involving $\vu_t$ on the left to obtain
\begin{alignat}{1}\nonumber
    \mG_{\eta, h}\vu_t + \eta \vu \cdot\nabla \vu+ \eta \nabla\zeta &= 
    \nabla\left\{\tfrac{1}{3}\eta^3\left(
	\vu\cdot\nabla\left(\nabla\cdot\vu\right) - \left(\nabla\cdot \vu\right)^2\right)
	+\tfrac{1}{2}\eta^2 \left(\vu\cdot \nabla\right)^2h \right\} \\  \label{Eq:MomEq}
	&\phantom{=}
	-\eta\left\{\tfrac{1}{2}\eta\left(\vu\cdot\nabla\left(\nabla\cdot\vu\right)-\left(\nabla\cdot \vu\right)^2\right)+\left(\vu\cdot \nabla\right)^2h\right\}\nabla h\,,     
\end{alignat}
where
\begin{align}\label{Def:G}
    \mG_{\eta, h}\vu &\equiv \eta \vu - \nabla \left( \tfrac{1}{3}\eta^3 \;\nabla\cdot \vu \right) -\nabla\left(\tfrac{1}{2}\eta^2\nabla h \cdot \vu \right) + \tfrac{1}{2}\eta^2\nabla h \; (\nabla \cdot \vu ) +\eta\nabla h \; \left(\nabla h \cdot \vu \right) \, .
\end{align}
Note that the last three terms in $\mG_{\eta,h}$ encapsulate the bathymetry effect and vanish when $\nabla h = 0$.  

Next, commute the time derivative $\mG_{\eta, h} (\vu_t) = (\mG_{\eta,h} \vu)_t - \mG_{\eta_t, h} \vu$ in \eqref{Eq:MomEq}, and use the $\eta$ equation in \eqref{Eq:dswe} to express $(\mG_{\eta, h})_t$ in terms of $\eta, \vu, h$.  After introducing the new variable $\vec{U} = \mG_{\eta,h} \vu$, one obtains the SGN equations (or DSWEs) in \emph{constraint form} 
\begin{subequations}\label{Eq:DSWEcf}
\begin{align} \label{Eq:DSWEcf1}
        &\eta(\vec{x},t) =h(\vec{x})+\zeta(\vec{x},t) \, , \phantom{\left(\frac{1}{3}\right)} \\ \label{Eq:DSWEcf2}
        &\eta_t + \nabla\cdot\left(\eta\vec{u}\right) = 0 \, , \\ \label{Eq:DSWEcf3}
        &\vec{U}_{t} + \eta\nabla\zeta+\nabla(\eta\vec{u}\otimes\vec{u}) = \vec{F}(\eta, \vec{u}, h) \, , \phantom{\left(\frac{1}{3}\right)}\\ \label{Eq:DSWEcf4}
        &\vec{U} - \mG_{\eta, h}\vu = 0 \, ,   
\end{align}
\end{subequations}
where the forcing $\vec{F}(\eta, \vec{u}, h)$ is
\begin{align}\nonumber
 	\vec{F}(\eta, \vec{u}, h) &\equiv 
 	\nabla\left\{\tfrac{1}{3}\eta^3\left(\vu\cdot\nabla\left(\nabla\cdot\vu\right)
 	-\left(\nabla\cdot \vu\right)^2\right)\right. 
        +\eta^2\nabla\cdot\left(\eta\vu\right)\nabla\cdot\vu 
	\\ \nonumber
	&\hspace{3cm}\left. +\tfrac{1}{2}\eta^2\left(\vu\cdot \nabla\right)^2h 	
 	+\eta\nabla\cdot\left(\eta\vu\right)\nabla h\cdot\vu\right\}
	\\ \nonumber
	&-\nabla h\left\{\tfrac{1}{2}\eta^2\left(\vu\cdot\nabla\left(\nabla\cdot\vu\right)
	-\left(\nabla\cdot \vu\right)^2\right)\right. 
        +\eta\nabla\cdot\left(\eta\vu\right)\left(\nabla\cdot\vu\right)
	\\ \label{Def:bigF}
	&\hspace{3cm}\left.+\eta\left(\vu\cdot \nabla\right)^2h	
	+\nabla\cdot\left(\eta\vu\right)\left(\nabla h\cdot\vu\right)\vphantom{\frac{1}{2}}\right\} \, .
\end{align}
In equation \eqref{Eq:DSWEcf}, the term involving $\vu \otimes \vu$ has components 
\begin{align*}
    \nabla(\eta \vu \otimes \vu) = \begin{pmatrix} (\eta u u)_x+(\eta u v)_y \\ (\eta u v)_x+(\eta v v)_y \end{pmatrix} \qquad 
\textrm{where} \qquad \vec{u}=(u,v) \, .
\end{align*}
The presence of the constraint in \eqref{Eq:DSWEcf} highlights that spatial discretizations (e.g., via method of lines) will lead to an index-1 differential algebraic equation (DAE).

For the numerical computations (\S~\ref{Sec:numerical_tests}) of the SGN equations, we adopt periodic ($d = 1$) or doubly periodic ($d = 2$) boundary conditions for $\eta$ and $\vu$. However, the main results for the operator $\mG_{\eta,h}$ in \S~\ref{Sec:BoundonG} will be formulated for both periodic and homogeneous Neumann (i.e., $\vec{n} \cdot \vu = 0$ on $\partial \Omega$) boundary conditions. 

\begin{remark}\label{rmk:dswe_as_constraint} (Conservation equation when $\nabla h = 0$)
    When the bathymetry is flat, i.e., $\nabla h = 0$, equations \eqref{Eq:DSWEcf} yield a conservative system.  In this case, $\vec{F} = \nabla W$ reduces to a gradient of a scalar $W(\eta, \vu)$, and $\nabla\eta = \nabla \zeta$. The time-evolution equations become
    \begin{subequations}\label{Eq:ConservativeForm}
    \begin{align}\label{EQ:ConservativeForm1}
        \eta_t + \nabla\cdot\left(\eta\vec{u}\right) &= 0  \,, \\ \label{EQ:ConservativeForm2}
        \vec{U}_{t} + \nabla \left( \tfrac{1}{2}\eta^2 + \eta\vec{u}\vec{u} - W\right) &= 0  \,, \\ \label{Eq:ConservativeForm3}
        \vec{U} - \mG_{\eta} \vu &= 0 \,.
    \end{align}  
    \end{subequations}
    The system \eqref{Eq:ConservativeForm} conserves the mass $m$, horizontal momentum $\vec{M}$, and total energy $E$ defined by
    \begin{align*}
    	m=\int_{\Omega}\zeta\du \vec{x}\,,\qquad \vec{M}=\int_{\Omega}\vec{U} \du \vec{x}\, , 
        \qquad     
        E=\tfrac{1}{2}\int_{\Omega} \left [ \zeta^2+ \eta(\vu\cdot\vu)
	   +\tfrac{1}{3}\eta^3(\nabla\cdot \vu)^2 \right ]\du \vec{x}\,. 	   
    \end{align*}\myremarkend
\end{remark}

\section{Mathematical Background for the Constraint Equation}\label{sec:linoper_background}
This section collects the mathematical notation and background relevant for the linear operators in the paper. Specifically, we introduce (\S~\ref{Subsec:linop}) the bilinear (weak) form of the operator $\mG_{\eta, h}$, and the operator $\mA$. Subsection~\ref{Sub:PCG} provides background on conditioning numbers relevant for the use in preconditioned conjugate gradient (PCG) methods.  

\subsection{Linear Operators and Bilinear Forms}\label{Subsec:linop} 
Throughout, $\Omega \subseteq \mathbb{R}^2$ is either an open connected domain with Lipschitz boundary $\partial \Omega$, or $\Omega = \mathbb{T}^2$ is the (doubly-periodic) torus with side lengths $L_x, L_y$.  Let 
\begin{align*}
    \mA = \sigma I - \alpha \nabla \cdot \nabla \,
\end{align*}
where $\sigma, \alpha > 0$ are real constants (\S~\ref{Sec:BoundonG} will provide formulas for $(\sigma, \alpha)$).  

The bilinear forms associated to the differential operators $\mG_{\eta,h}$ and $\mA$ will be defined on the Hilbert space $H_{\rm div}(\Omega)$ given by
\begin{align*}
    H_{\rm div}(\Omega) &\equiv \left\{ \vu \in [L^2(\Omega)]^2 \; | \; \nabla \cdot \vu \in L^2(\Omega), \; 
	\vec{n}\cdot \vec{u} = 0 \textrm{ on } \partial \Omega \right\} \,,
\end{align*}
which is endowed with the norm 
\begin{align*}
	\|\vu\|_{\rm div}^2 \equiv \int_{\Omega} \left[ |\vu|^2 + (\nabla \cdot \vu)^2 \right] \du \vec{x} \,.
\end{align*}
When $\partial \Omega$ is Lipschitz, the outward boundary unit normal, $\vec{n}(\vec{x})$, is defined for almost every $\vec{x} \in \partial \Omega$, and the normal trace of functions $\vec{u} \in H_{\rm div}(\Omega)$ is understood to be $\vec{n}\cdot \vec{u} = 0$ in $H^{-1/2}(\partial \Omega)$, e.g., see \cite[Chapter 3.5.2]{Monk2003}.  When $\Omega = \mathbb{T}^2$, the boundary conditions in $H_{\rm div}(\Omega)$ (i.e., $\vec{n}\cdot \vec{u} = 0 \textrm{ on } \partial \Omega$) are replaced by periodic boundary conditions on $\vu$.  

Although the SGN equation solution $\eta = \eta(\vec{x}, t)$ is a function of $t$, for the purposes of preconditioning, $t$ enters only as a parameter in $\mG_{\eta,h}$. When discussing the weak forms we therefore consider $\eta(\vx)$ as a function of $\vx$ only. Furthermore, we assume that $\eta(\vx), \nabla h(\vx) \in L^{\infty}(\Omega)$ and that $\eta(\vx) \geq \eta_{\rm min} > 0$ is bounded strictly away from zero.  This assumption excludes some physically relevant scenarios such as the modeling of ``dry'' states where $\eta(\vec{x},t) = 0$ in regions of the domain.  Throughout we keep $\sigma > 0, \alpha > 0$ as arbitrary constant values and use $\sigma^\star, \alpha^\star$ to refer to our concrete formulas used in the operator $\mA$.  

The bilinear forms $A : H_{\rm div}(\Omega) \times H_{\rm div}(\Omega) \rightarrow \mathbb{R}$, and $B : H_{\rm div}(\Omega) \times H_{\rm div}(\Omega) \rightarrow \mathbb{R}$ corresponding to the weak forms of the operators $\mA$ and $\mG_{\eta, h}$ respectively are:
\begin{align}
    A[\vu, \vv] &\equiv \int_{\Omega} \left[\sigma \, \vu \cdot \vv  + \alpha \, \big(\nabla \cdot \vv\big) \big(\nabla \cdot \vu\big) \right] \, \du \vx \quad \quad \forall \vu, \vv \in H_{\rm div}(\Omega) \, ,     
\end{align}
and
\begin{align}\nonumber
    B[\vu, \vv] &\equiv \int_{\Omega} \Big[ \eta \, \vu \cdot \vv  
    + \tfrac{1}{3} \eta^3 \, \big(\nabla\cdot \vu\big) \big(\nabla \cdot \vv\big) 
    + \tfrac{1}{2} \eta^2 \, \big(\nabla \cdot \vv \big) \big(\nabla h \cdot \vu \big)  \\
    & \quad \phantom{\int}+ \tfrac{1}{2} \eta^2 \, \big(\nabla \cdot \vu\big) \big(\nabla h \cdot \vv \big) 
    + \eta \, \big(\nabla h \cdot \vu \big) \big(\nabla h \cdot \vv \big) \Big] \, \du \vx \, \quad \quad \forall \vu, \vv \in H_{\rm div}(\Omega) \,.
\end{align}
The bilinear forms arise in the usual fashion from the strong form of the PDE through integration against a test function, e.g., $A[\vu, \vv] = (\vu, \mA \vv)$ and $B[\vu, \vv] = (\vu, \mG_{\eta,h} \vv)$, for $\vu,\vv$ (where $( \cdot, \cdot)$ is the $L^2$ inner product) being twice continuously differentiable with vanishing normal component (cf. \cite[Chapter~12]{Schechter2001} for additional details on operator extensions via bilinear forms).

It is clear from the definitions that both $A[\cdot, \cdot]$ and $B[\cdot, \cdot]$ are symmetric, i.e., $A[\vu, \vv] = A[\vv, \vu]$ and $B[\vu,\vv] = B[\vv, \vu]$. Moreover, $A[\cdot, \cdot]$ is positive definite ($A[\vu,\vu] > 0$ for $\vu \neq 0$), and it is the case $B[\cdot, \cdot]$ is positive definite as well.

We will often write $\mA^{-1} \mG_{\eta, h}$ etc. These expressions are understood via the following definitions
\begin{enumerate}
    \item (The operator $\mA^{-1} \mG_{\eta,h}$) For $\vec{f} \in H_{\rm div}(\Omega)$, the function $\vec{w} = \mA^{-1} \mG_{\eta,h} \vec{f}$ is defined by:
        Find $\vec{w} \in H_{\rm div}(\Omega)$ such that 
        \begin{align}\label{Eq:Weakform}            
            A[\vec{w}, \vv] &= B[\vec{f}, \vv] \; \quad\quad \forall \; \vv \in H_{\rm div}(\Omega) \, .
        \end{align}
    \item (The operator $\mG_{\eta,h}^{-1} \mA$) For $\vec{f} \in H_{\rm div}(\Omega)$, the function $\vec{w} = \mG_{\eta,h}^{-1}\mA \vec{f}$ is defined by: Find $\vec{w} \in H_{\rm div}(\Omega)$ such that 
        \begin{align}\label{Eq:WeakformInvG}
            B[\vec{w}, \vv] &= A[\vec{f}, \vv] \; \quad\quad \forall \; \vv \in H_{\rm div}(\Omega) \, .
        \end{align}
    \item (Generalized eigenvalues) The generalized eigenvalue problem $\mG_{\eta,h} \vu = \lambda \mA \vu$ is: Find $(\vu, \lambda) \in H_{\rm div}(\Omega) \times \mathbb{R}$ ($\vu \neq 0$) such that
        \begin{align}\label{Eq:WeakGenEV}
            \lambda A[\vec{u}, \vv] &= B[\vu, \vv ] \; \quad\quad \forall \; \vv \in H_{\rm div}(\Omega) \, .
        \end{align}
\end{enumerate}
The fact that both $\lambda \in \mathbb{R}$ and $\vu$ can be restricted to a real-valued function in $H_{\rm div}(\Omega)$ in \eqref{Eq:WeakGenEV} follows from the symmetry of the bilinear forms.  

\subsection{Conjugate Gradient and Conditioning Numbers}\label{Sub:PCG}

Conjugate gradient (CG) is an iterative Krylov method for solving linear matrix systems of the form $\mat{A} \vec{x} = \vec{f}$ where $\mat{A} \in \mathbb{R}^n$ is symmetric positive definite (in exact arithmetic, CG converges after a finite number of iterations).  A straight-forward, but not necessarily sharp \cite[Chapter 38]{TrefethenBau1997}, bound on the convergence of CG is given in terms of the conditioning number $\kappa(\mat{A}) \equiv \lambda_{\rm max}/\lambda_{\rm min}$ where $\lambda_{\rm max}$, $\lambda_{\rm min}$ are the largest and smallest eigenvalues of $\mat{A}$:
\begin{align}\label{eq:CGRate}
    \| \vec{e}_k \|_{A} \leq 2 \left( \frac{\sqrt{\kappa} - 1}{\sqrt{\kappa} +1 }\right)^{k} \| \vec{e}_0\|_{A}
    \, , \qquad \|\vec{x}\|_{A} \equiv \vec{x}^T \mat{A} \vec{x} \, .    
\end{align}
In \eqref{eq:CGRate} $\vec{e}_k = \vec{x}_k - \vec{x}^*$ is the error after $k$ iterations of CG, $\vec{x}^* = \mat{A}^{-1} \vec{f}$ is the exact solution. 

The convergence of CG  may be improved (such as in settings where $\mat{A}$ has a large conditioning numbers) by solving a preconditioned system with a symmetric positive definite $\mat{P}$.  Applying CG to the following system, which shares the same solution as $\mat{A} \vec{x} = \vec{f}$, 
\begin{align}\label{Eq:PCGLinSyst}
    \mat{P}^{-1/2} \mat{A} \mat{P}^{-1/2} \vec{y} = \mat{P}^{-1/2} \vec{f}  \quad \quad \textrm{with} \quad \quad 
    \vec{y} = \mat{P}^{1/2} \vec{x} \, ,
\end{align}
results in a preconditioned conjugate gradient (PCG) method.  It is worth emphasizing that practical implementations of PCG avoid computations of matrix square roots and only require linear system solves involving $\mat{P}$.  The conditioning number of the PCG matrix in \eqref{Eq:PCGLinSyst} is then modified from that of $\mat{A}$ to be the ratio of maximum to minimum eigenvalues of $\mat{P}^{-1} \mat{A}$. Namely, the conditioning number, and upper bound on the convergence rate may be (significantly) reduced if $\mat{P}$ is ``close'' to $\mat{A}$. 

When $\mA$ is used as a preconditioner for $\mG_{\eta, h}$, the operator analogue of $\kappa$ for symmetric operators is
\begin{align}\label{Def:PCG_cond_num}
    \kappa_{\mA}(\mG_{\eta, h}) \equiv \frac{\mu_{\rm \max}}{\mu_{\rm \min}} \, ,
\end{align}
where $\mu_{\min}$ and $\mu_{\max}$ are the (generalized) Rayleigh quotients 
\begin{align}\label{Exp:RayleighQuotients}
    \mu_{\min} = \inf_{0 \neq \vu \in H_{\rm div}(\Omega)} \frac{B[\vu, \vu]}{A[\vu, \vu]} \quad\quad \textrm{and} \quad\quad
    \mu_{\max} = \sup_{0 \neq \vu \in H_{\rm div}(\Omega)} \frac{B[\vu, \vu]}{A[\vu, \vu]} \, .
\end{align}
The quantity $\kappa_{\mA}(\mG_{\eta, h})$ measures the ratio of the largest to smallest values in the spectrum of $\mA^{-1} \mG_{\eta,h}$.  Not only does $\kappa_{\mA}(\mG_{\eta,h})$ play the operator analogue of the PCG matrix conditioning number, but it also appears as a bound in the convergence of operator/Hilbert space versions of preconditioned gradient descent (see \cite{Kirby2010}) and PCG (e.g., variable coefficient Laplacian \cite{GillesTownsend2019}). A smaller value of $\kappa_{\mA}(\mG_{\eta,h})$ will yield a tighter bound on the convergence of PCG. 

Bounds on the continuum operators also carry over to the discrete setting as well.

\begin{remark}[Conditioning number bound for spatially discrete problems]\label{rmk:meshInd} When the operator bilinear forms $A[\cdot,\cdot]$ and $B[\cdot,\cdot]$ are discretized via finite element methods, the quantity $\kappa_{\mA}(\mG_{\eta,h})$ is a mesh-independent upper bound on the discrete matrix preconditioned conditioning number (see also \cite{Kirby2010}). 

Namely, let $V \subset H_{\rm div}(\Omega)$ be any finite dimensional subspace with basis $\{v_i\}_{i=1}^n \in V$. The discrete matrices that arise from the weak form \eqref{Eq:WeakGenEV} are symmetric and positive definite
\begin{align*}
    \mat{A}_{ij} = A[v_i, v_j] \qquad \mat{G}_{ij} = B[v_i, v_j] \qquad \textrm{for}\; i,j = 1, \ldots, n \, .
\end{align*}
Then $\kappa(\mat{A}^{-1} \mat{G}) \leq \kappa_{\mA}(\mG_{\eta,h})$. In particular, the discrete conditioning number $\kappa(\mat{A}^{-1} \mat{G})$ is bounded independent of the mesh. This ensures that  the convergence rate of PCG is bounded independent of the mesh. \myremarkend
\end{remark}

\section{Preconditioner for the Operator $\mG_{\eta, h}$}\label{Sec:BoundonG}
In this section we introduce the constant coefficient preconditioner 
\begin{align}\label{Def:A}
    \mA \vu = \sigma^\star \, \vu - \alpha^\star \, \nabla (\nabla \cdot \vu ) \, ,
\end{align}
and establish bounds on the generalized eigenvalues and conditioning number of $\mA^{-1} \mG_{\eta,h}$.  These bounds will play an important role in ensuring convergence for preconditioned conjugate gradient methods or linear stability in semi-implicit time-stepping. 

The coefficients in \eqref{Def:A} for $\mA$ are defined in terms of pointwise estimates of $\eta$ and $h$
\begin{align}\label{Eq:MaxMinWaveFields}
	|\nabla h|_{\rm max} \equiv \max_{\vec{x} \in \Omega} | \nabla h(\vec{x}) | \,, \quad\quad
	\eta_{\rm \min} \equiv \min_{\vec{x} \in \Omega} \eta(\vec{x}) \,, \quad\quad
	\eta_{\rm \max} \equiv \max_{\vec{x} \in \Omega} \eta(\vec{x}) \,.
\end{align}
\begin{itemize}
    \item For flat bathymetry  ($|\nabla h|_{\rm max} = 0$) we set
        \begin{align}\label{Eq:FlatbathCoeff}
            \sigma^\star \equiv \eta_{\rm max} \, , \qquad \qquad \alpha^\star \equiv \tfrac{1}{3} \eta_{\rm max}^3 \, ,
            \qquad \qquad 
            \gb &\equiv \left(\frac{\eta_{\rm max}}{\eta_{\rm min}} \right)^{3} \, .
        \end{align}
    \item For variable bathymetry  ($|\nabla h|_{\rm max} \neq 0$) we take
        \begin{align}\label{Eq:Optcoeff}
            \sigma^\star &\equiv \, \max_{\vx \in \Omega} \, \left( \eta(\vx) + \tfrac{1}{6}(4 + \sqrt{13}) \eta(\vx)|\nabla h(\vx) |^2\right) \, , \quad  \quad 
            \alpha^\star \equiv \tfrac{1}{6}(4 + \sqrt{13}) \, \eta_{\rm max}^3 \, , \\ \label{Eq:Optcoeffb}
            \gb &\equiv \, \max\left\{ \frac{\sigma^\star}{\eta_{\rm min}}, \left(\frac{4+\sqrt{13}}{4-\sqrt{13}}\right)\left(\frac{\eta_{\rm max}}{\eta_{\rm min}} \right)^{3} \right\}  \, .
    \end{align}    
\end{itemize}
Note that $|\nabla h|_{\rm max}$ quantifies the maximum slope of the bathymetry with $|\nabla h|_{\rm max} = 0$ corresponding to a flat bottom. When solving the SGN equations, $|\nabla h|_{\rm max}$ can be determined as a one-time computation since $h(\vec{x})$ is prescribed as part of the problem data.

The significance of the coefficients is highlighted by the next theorem.

\begin{theorem}\label{Thm:CondNumBnd} Assume that $\eta, \nabla h \in L^{\infty}(\Omega)$, and $\eta_{\rm min} > 0$. Let $\mG_{\eta, h}$ be defined as above and $\mA$ defined with the coefficients \eqref{Eq:FlatbathCoeff}--\eqref{Eq:Optcoeff}. Then the condition number $\kappa_{\mA}(\mG_{\eta,h})$ of $\mA^{-1}\mG_{\eta,h}$ is bounded by 
\begin{align}\label{Eq:CondNumBound}
    \kappa_{\mA}( \mG_{\eta,h}) \leq \gb \, . 
\end{align}
    Moreover, the generalized eigenvalues $\lambda$ of $\mG_{\eta, h} \, \vu = \lambda \, \mA \vu$ are real and lie in the interval
    \begin{align}\label{Eq:EigValueBound}
        (\gb)^{-1} \leq \lambda \leq 1 \, .         
    \end{align}
\end{theorem}
The proof of Theorem~\ref{Thm:CondNumBnd} is deferred to Appendix~\ref{Appendix:Pf_MainResult}.

The inequalities \eqref{Eq:CondNumBound} and \eqref{Eq:EigValueBound} are important for several reasons. Firstly, smaller values of $\kappa_{\mA}( \mG_{\eta,h})$ are desirable as they provide stronger performance bounds for preconditioned conjugate gradient when solving systems involving $\mG_{\eta,h}$. 

Without the use of a preconditioner $\mA$, the eigenvalues of $\mG_{\eta,h}$ would be unbounded and bounds such as \eqref{Eq:CondNumBound} and \eqref{Eq:EigValueBound} would not exist. This would render methods such as conjugate gradient far less effective (or fail entirely) for solving $\mG_{\eta,h} \vec{u} = \vec{f}$.  While the bounds \eqref{Eq:CondNumBound} and \eqref{Eq:EigValueBound} are formulated in terms of PDE operators, they provide mesh-independent PCG convergence rates for the spatially discrete problem as well (see Remark~\ref{rmk:meshInd})

Secondly, when semi-implicit time-stepping approaches are used to solve the SGN equations, the clustering of eigenvalues $\lambda$ near $1$ is important to ensure linear stability. 

The formula \eqref{Eq:Optcoeffb} highlights two ways in which the upper bound on $\kappa_{\mA}(\mG_{\eta,h})$ may diverge and lead to a loss in practical effectiveness of $\mA$
\begin{itemize}
    \item (Large contrast ratio) If $\eta_{\max}/\eta_{\min} \gg 1$, and the second term in \eqref{Eq:Optcoeffb} dominates then $\gb \sim (\eta_{\max}/\eta_{\min})^3$. 

    \item (Large variations in bathymetry) If $\max_{\vec{x} \in \Omega}\eta(\vec{x}) |\nabla h(\vec{x})|^2 \gg 1$ so that the first term in \eqref{Eq:Optcoeffb} dominates, then at worst $\gb \sim \eta_{\max} |\nabla h|_{\max}^2$.         
\end{itemize}

We conclude this section with several remarks.

\begin{remark} (Robustness of the coefficients $\sigma^{\star}, \alpha^{\star}$) One may ask if alternative formulas may be used for $(\sigma, \alpha)$ and still yield an upper bound on $\kappa_{\mA}(\mG_{\eta,h})$? It turns out that  $\kappa_{\mA}(\mG_{\eta,h})$ is bounded for any positive coefficients $\sigma, \alpha > 0$ in $\mA$. Alternative choices of $\sigma, \alpha$ may however yield larger values of $\kappa_{\mA}(\mG_{\eta,h})$ and decrease the efficiency of the preconditioner $\mA$ or increase the error constant in semi-implicit time integration. In practice, this robustness allows for some flexibility in under or overestimating the formulas for $\sigma^{\star}, \alpha^{\star}$. It turns out, see Appendix~\ref{Appendix:QuasiOpt}, that the choice of $\mA$ is \emph{quasi-optimal} within a class of constant coefficient operators.
\myremarkend
\end{remark}

\begin{remark}\label{rmk:replacementcoeff} (Replacing $\eta(\vec{x}) |\nabla h(\vec{x})|^2$ with $\eta_{\rm max} |\nabla h|_{\rm max}^2$ in $\mA$) The coefficient $\sigma^\star$ relies on the maximum pointwise bound for $\eta(\vec{x}) |\nabla h(\vec{x})|^2$. In some practical settings, it may be simpler to individually estimate $|\nabla h|_{\rm max}$ and $\eta_{\rm max}$ separately.  Instead of \eqref{Eq:Optcoeff}, one can take:
$\mA = \sigma^{\star\star} I - \alpha^{\star\star} \nabla \nabla \cdot $ where
\begin{align}
    \sigma^{\star\star} = \eta_{\rm max} \left( 1 + \tfrac{1}{6}(4+\sqrt{13}) |\nabla h|_{\rm max}^2\right) \, , 
    \qquad 
    \alpha^{\star\star} = \tfrac{1}{6}(4 + \sqrt{13}) \, \eta_{\rm max}^3 \, .
\end{align}
These coefficients still ensure the eigenvalue bound \eqref{Eq:EigValueBound} and conditioning number bound \eqref{Eq:CondNumBound} where $\gb$ is replaced with a (less optimal) bound 
\begin{align}
    \gb &=  
    \max\left\{ \frac{\eta_{\rm max}}{\eta_{\rm min}} \left(1 + \tfrac{1}{6}(4+\sqrt{13}) |\nabla h|_{\rm max}^2\right),  \left(\frac{4+\sqrt{13}}{4-\sqrt{13}}\right)\left( \frac{\eta_{\rm max}}{\eta_{\rm min}} \right)^{3} \right\}  \, .
\end{align}
Unlike \eqref{Eq:Optcoeff}, these coefficients are not \emph{quasi-optimal} in the sense discussed in Appendix~\ref{Appendix:QuasiOpt}, and can result in a drop in PCG performance (see \S\ref{subsec:1dpreconditioner}). However, they are still effective coefficients for both PCG and to define the implicit $\mA$ in semi-implicit time-stepping. \myremarkend
\end{remark}

\begin{remark}[Riesz Preconditioner]\label{rmk:Rieszpreconditioner}    
    The choice of $A[\cdot, \cdot]$ can be viewed as designing a Riesz preconditioner for the operator $B[\cdot, \cdot]$; for a discussion of Riesz preconditioners in the context of bilinear forms see \cite{Kirby2010, MalekStrakos2015, GillesTownsend2019}. \myremarkend
\end{remark}

\section{Time-stepping for the SGN Equations}\label{Sec:TimeStepping}
This section provides details on the time integration schemes we use to solve \eqref{Eq:DSWEcf}.  Subsections \ref{Subsec:ERK_PCG} and \ref{Subsec:LMM} present explicit Runge-Kutta (ERK) and multistep methods that handle the constraint through preconditioned conjugate gradient (PCG) methods. Subsection \ref{Subsec:IMEX_LMM} presents linearly implicit multistep (ImEx-type) schemes that treat $\mA$ implicitly and $\mG_{\eta,h}$ explicitly. 

The linearly implicit multistep schemes admit several advantages---they require only one linear solve involving $\mA$ per time step and they do not suffer from order reduction (cf. Remark~\ref{Rmk:OrderReduction}).  However, linear implicit schemes do require a precise choice of $\mA$ to ensure zero-stability, which will be discussed in subsection \ref{subsec:zerostability}.

In this section, it is useful to view the SGN system \eqref{Eq:DSWEcf} as an evolution equation with a constraint
\begin{subequations}\label{Eq:DAE_DSWE}
\begin{align}\label{Eq:DAE1_DSWE}
        \frac{\partial}{\partial t} \begin{pmatrix}
            \eta \\
            \vec{U} 
        \end{pmatrix} &= \begin{pmatrix}
            -\nabla\cdot(\eta\vec{u}) \\
            -\eta \nabla \zeta - \nabla \left( \eta \vu \otimes \vu \right) + \vec{F} 
    \end{pmatrix} \, , \\ \label{Eq:DAE2_DSWE}
    \mG_{\eta, h} \, \vu &= \vec{U} \, , 
\end{align}
\end{subequations}
where $\vec{F}$ is defined in \eqref{Def:bigF}.  Spatial discretizations of \eqref{Eq:DAE_DSWE} via a method-of-lines give rise to a differential algebraic equation (DAE) 
\begin{subequations}\label{Eq:DAE}
\begin{align}\label{Eq:DAE1}
	\vec{w}_t &= \vec{f}(\vec{w}, \vec{p}) \, , \\ \label{Eq:DAE2}
	\vec{g}(\vec{w}, \vec{p}) &= 0 \, ,
\end{align} 
\end{subequations}
where $\vec{w} \in \mathbb{R}^n$, $\vec{p} \in \mathbb{R}^m$ and $\vec{f} : \mathbb{R}^n \times \mathbb{R}^m \rightarrow \mathbb{R}^n$ and $\vec{g} : \mathbb{R}^n \times \mathbb{R}^m \rightarrow \mathbb{R}^m$.  Specifically, the ODE variables are spatial approximations given by $\vec{w} \approx (\eta, \vec{U})^T$, $\vec{p} \approx \vec{u}$, $\vec{f} \approx (-\nabla\cdot(\eta\vec{u}), \; -\eta\nabla \zeta - \nabla(\eta \vec{u}\otimes\vec{u})+ \vec{F})^T$ and $\vec{g} \approx \mathcal{G}_{\eta,h} \vec{u} - \vec{U}$. 

In general, \eqref{Eq:DAE} will be an \emph{index}-1 DAE since the constraint equation \eqref{Eq:DAE2} can be solved to implicitly define $\vec{p}$ as a function of $\vec{w}$.  This is because $\mathcal{G}_{\eta,h}$ is an invertible operator (\S\ref{Sec:BoundonG}), and matrix approximations to $\mG_{\eta,h}$ constructed by testing $B[\cdot , \cdot]$ with Galerkin or finite element subspaces of $H_{\rm div}$ are also invertible. 

For simplicity we adopt constant time steps $\Delta t > 0$ where $t_n = n \Delta t$, however the results for Runge-Kutta methods naturally extend to variable step sizes. The fluid variables at time $t_n$ are written as $\vu_n \approx \vec{u}(t_n)$, $\eta_n \approx \eta(t_n)$, $\zeta_n \approx \zeta(t_n)$, and  $\vec{U}_n \approx \vec{U}(t_n)$.  

\begin{remark}
    While the focus here is on time discretizations of \eqref{Eq:DSWEcf}, the same discretization approaches can be applied to the non-constraint form of the equations \eqref{Eq:dswe}. \myremarkend
\end{remark}

\subsection{Explicit Runge-Kutta with Preconditioned Conjugate Gradient}\label{Subsec:ERK_PCG}
We discretize the SGN equations \eqref{Eq:DSWEcf} via the classical explicit RK4 scheme defined by the Butcher coefficients \cite[Chapter~5]{Leveque2007}
\begin{equation*}
	(a_{ij})_{i,j=1}^4 = \begin{pmatrix}
	    & & & \\
        \tfrac{1}{2} &		& 	& \\
        0	&  \tfrac{1}{2} 	& 	& \\ 
        0   &   0 	& 1 & 
	\end{pmatrix} \, , \qquad 
    (b_i)_{i=1}^4 = \begin{pmatrix}
        \tfrac{1}{6} & \tfrac{1}{3} & \tfrac{1}{3} & \tfrac{1}{6}
    \end{pmatrix}^T \, ,
    \qquad 
    (c_i)_{i=1}^{4} = \begin{pmatrix} 
		                  0 &
		                \tfrac{1}{2} &
		              \tfrac{1}{2} & 
		              1	
                      \end{pmatrix}^T \, .
\end{equation*}
This yields the following equations for the stages $(i = 1, \ldots, 4$)
\begin{subequations}\label{Eq:ERK}
\begin{align}    
    \eta_n^{(i)} &= \eta_n + \dt  \sum_{j = 1}^{i-1} a_{ij} \left( 
    -\nabla\cdot\left(\eta_n^{(j)} \vu_n^{(j)} \right)  \right) \, , \\ \label{Eq:ERK_uStage}
    \vec{U}_n^{(i)} &= \vu_n + \dt \sum_{j=1}^{i-1} a_{ij} \; \left(
    -\eta_n^{(j)} \nabla \zeta_n^{(j)} - \nabla \left( \eta_n^{(j)} \vu_n^{(j)} \otimes \vu_n^{(j)} \right) + \vec{F}_n^{(j)} 
    \right)
    \,, \\ \label{Eq:ERK_CStage}
    \mG_{\eta_n^{(i)}, h} \; \vu_n^{(i)} &=  \vec{U}_n^{(j)} \, , \phantom{\sum_{i=1}^{i-1}} 
\end{align}
\end{subequations}
where  $\zeta_n^{(i)} =\eta_n^{(i)} - h$. The final update is
\begin{align*}
    \eta_{n+1} &= \eta_n + \dt \sum_{j=1}^s b_j \left( -\nabla\cdot\left(\eta_n^{(j)} \vu_n^{(j)} \right)  \right) \,,\\
    \vu_{n+1} &= \vu_n + \dt \sum_{j=1}^s b_j \left(
    -\eta_n^{(j)} \nabla \zeta_n^{(j)} - \nabla \left( \eta_n^{(j)} \vu_n^{(j)} \otimes \vu_n^{(j)} \right) + \vec{F}_n^{(j)} 
    \right)\,. 
\end{align*}
While RK4 is an explicit method, the algebraic constraint in \eqref{Eq:ERK_CStage} still requires a linear solve for $\mG_{\eta,h}$ at each stage. To solve \eqref{Eq:ERK_CStage}, we apply PCG with the preconditioner $\mA$ defined by the coefficients \eqref{Eq:Optcoeff}. Since each solve is independent, one can change the coefficients of $\mA$ as needed for each stage or step.  

\begin{remark}[Order reduction in RK schemes]\label{Rmk:OrderReduction} Runge-Kutta methods may exhibit a reduction in order when applied to initial boundary value problems. The reduction often stems from the presence of boundary conditions or a lack of regularity in the solution. Due to the periodic domain, the test cases we employ do not exhibit order reduction. However, the RK schemes presented in this section could exhibit order reduction when applied to problems with time-dependent boundary conditions (multistep methods do not exhibit order reduction). \myremarkend
\end{remark}

\subsection{Explicit Multistep with Preconditioned Conjugate Gradient}\label{Subsec:LMM}
Similar to explicit RK methods, we also test explicit Adams-Bashforth with PCG used to solve the constraint equation at each step. A general s-stage linear multistep method discretizes \eqref{Eq:DAE} with coefficients $(a_j, b_j)_{j=1}^s$ as
\begin{subequations}\label{Eq:LMM_DAE}
\begin{align}\label{Eq:LMM_DAEa}
        \frac{1}{\Delta t}\sum_{j=0}^{s} a_j \vec{w}_{n+j} &= 
        \sum_{j = 0}^{s-1} b_j \vec{f}(\vec{w}_{n+j}, \vec{p}_{n+j}) \, , \\ \label{Eq:LMM_DAEb}
         \vec{g}(\vec{w}_{n+j}, \vec{p}_{n+j}) &= 0 \, .
\end{align}
\end{subequations}
The coefficients for Adams-Bashforth (with orders 2, 3 and 4) are included in Table~\ref{Table:LMMCoeffs}.  

Applying the discretization to the SGN yields
\begin{subequations}\label{Eq:AB}
\begin{align}\label{Eq:AB_1}	
    \zeta_{n} &= \eta_{n} - h \, , \phantom{\sum_{j=1}^s}  \\
    \frac{1}{\dt}\sum_{j=0}^{s} a_j \eta_{n+j} &= 
    -\sum_{j=0}^{s-1} b_j \nabla\cdot\left(\eta_{n+j}  \vu_{n+j} \right)   \, ,\\   \label{Eq:AB_2}
    \frac{1}{\dt}\sum_{j=0}^{s} a_j \vec{U}_{n+j} &= 
    \sum_{j=0}^{s-1} b_j\Big( 
    - \eta_{n+j} \nabla \zeta_{n+j} - \nabla(\eta_{n+j} \vu_{n+j} \otimes \vu_{n+j} ) + \vec{F}(\eta_{n+j}, \vu_{n+j}, h) 
    \Big) \, , \\   \label{Eq:AB_3}
    \mG_{\eta_{n+j}, h} \vu_{n+j} &= \vec{U}_{n+j} \, . 
\end{align}
\end{subequations}
When solving this system, we again use the preconditioner $\mA$ with PCG to solve the constraint equation at each time step.

\begin{table}[!htb]
\setlength{\mylength}{0.15\textwidth}
\bigskip
\noindent
\begin{tabular}{ |@{~}>{\centering\arraybackslash}p{0.08\textwidth} @{~} |@{~}>{\centering\arraybackslash} p{0.03\textwidth} @{~} || @{~}>{\centering\arraybackslash}p{\mylength} @{~}| @{~}>{\centering\arraybackslash}p{\mylength} @{~}|@{~} >{\centering\arraybackslash}p{\mylength} @{~} |@{~} >{\centering\arraybackslash}p{\mylength} @{~} |@{~} >{\centering\arraybackslash}p{\mylength} @{~} |}
\hline
	  &  & $j=4$ & $j = 3$ & $j = 2$ & $j = 1$ & $j = 0$  \\ \hline
    AB2 & $a_j$ & $\cdot$ & $\cdot$ & 1 & -1 & 0 \\          
      & $b_j$  &  $\cdot$ & $\cdot$ & 0 &  $\frac{3}{2}$ & $-\frac{1}{2}$ \\ \hline
    AB3 & $a_j$ & $\cdot$ & 1 & -1 & 0 & 0 \\      
      & $b_j$  &  $\cdot$ & 0 & $\frac{23}{12}$ & $-\frac{16}{12}$ & $\frac{5}{12}$ \\ \hline
    AB4 & $a_j$ & 1 & -1 & 0 & 0 & 0 \\      
      & $b_j$  &  0 & $\frac{55}{24}$ & $-\frac{59}{24}$ & $\frac{37}{24}$ & $-\frac{9}{24}$ \\ \hline
    SBDF2 & $a_j$  & $\cdot$ & $\cdot$ & $\frac{3}{2}$ & $-2$ & $\frac{1}{2}$ \\
          & $c_j$  & $\cdot$ & $\cdot$ & 1 & 0 & 0 \\
          & $b_j$  & $\cdot$ & $\cdot$ & 0 & $2$ & $-1$ \\ \hline
\end{tabular}
    \caption{Multistep coefficients for $r$th order Adams–Bashforth (ABr) with $r=2,3,4$, and $2$nd order semi-implicit backward differentiation (SBDF2).}
	\label{Table:LMMCoeffs}
\end{table}

\subsection{A Linearly Implicit Method for the SGN Equations}\label{Subsec:IMEX_LMM}
In contrast to the previous two subsections, this subsection provides a linearly implicit treatments of $\mG_{\eta,h}$ in \eqref{Eq:AB}.  The starting point is to split the algebraic equation by first adding and subtracting a linear term $\mat{A} \vec{p}$ to obtain
    \begin{align}\label{Eq:SplitDAE1}
        \mat{A}\vec{p} - \mat{A}\vec{p} + \vec{g}(\vec{w}, \vec{p}) = 0 \, .
    \end{align}
    Next, using a set of coefficients $(c_j)_{j=0}^{s}$, extrapolate the first term in \eqref{Eq:SplitDAE1} at time $t + \tau$,
    \begin{align}\label{Eq:ApproxDAE1}
        \mat{A}\,\vec{p}(t + \tau) \approx \sum_{j=0}^s c_j \mat{A} \vec{p}_{n+j} + \mathcal{O}(\Delta t^{p}) \, . 
    \end{align}
    The coefficient $c_s \neq 0$ ensures the approximation of $\mat{A}\vec{p}$ will lead to an implicit equation in $\vec{p}$. The remaining term $\mat{A}\vec{p} - \vec{g}(\vec{w}, \vec{p})$ is extrapolated to time $t + \tau$ via 
    \begin{align}\label{Eq:ApproxDAE2}
        \mat{A}\vec{p}(t+\tau) - \vec{g}\big(\vec{w}(t+\tau), \vec{p}(t+\tau)\big)  \approx 
        \sum_{j=0}^{s-1} b_j \mat{A} \vec{p}_{n+j} - b_j \vec{g}( \vec{w}_{n+j}, \vec{p}_{n+j} ) 
        + \mathcal{O}(\Delta t^{p}) \, . 
    \end{align}        
    The equation \eqref{Eq:LMM_DAEb} in \eqref{Eq:LMM_DAE} is then replaced by 
    \begin{align}\label{Eq:IMEX_DAE}
            &0 = \sum_{j=0}^{s} c_j \mat{A} \vec{p} +\sum_{j=0}^{s-1} b_j \Big( \vec{g}(\vec{w}_{n+j}, \vec{p}_{n+j}) - \mat{A} \vec{p}_{n+j} \Big) \,,  
    \end{align}    
    which is obtained from equating \eqref{Eq:ApproxDAE1} and \eqref{Eq:ApproxDAE2}.
	By construction the scheme \eqref{Eq:LMM_DAEa} with \eqref{Eq:IMEX_DAE} is linearly implicit in $\mat{A}$; no solution involving $\vec{g}$ is required.  

    We use the SBDF2 coefficients for $a_j, b_j, c_j$ presented in Table~\ref{Table:LMMCoeffs}, together with the time-discretization \eqref{Eq:IMEX_DAE} where $\mat{A} \approx \mA$ to yield the linearly implicit scheme
\begin{subequations}\label{Eq:MM}
\begin{align}\label{Eq:MM_1}	
    \zeta_{n} &= \eta_{n} - h \, , \phantom{\sum_{j=1}^s}  \\
    \frac{1}{\dt}\sum_{j=0}^{s} a_j \eta_{n+j} &= 
    -\sum_{j=0}^{s-1} b_j \nabla\cdot\left(\eta_{n+j}  \vu_{n+j} \right)   \, ,\\   \label{Eq:MM_2}
    \frac{1}{\dt}\sum_{j=0}^{s} a_j \vec{U}_{n+j} &= \phantom{-} 
    \sum_{j=0}^{s-1} b_j\Big( 
    - \eta_{n+j} \nabla \zeta_{n+j} - \nabla(\eta_{n+j} \vu_{n+j} \otimes \vu_{n+j} ) + \vec{F}(\eta_{n+j}, \vu_{n+j}, h) 
    \Big) \, , \\   \label{Eq:MM_3}
    \sum_{j=0}^s c_j \; \mA \: \vu_{n+j} &= 
    \phantom{-}\sum_{j=0}^{s-1} b_j \left( \mA \vu_{n+j} - \mG_{\eta_{n+j}, h} \vu_{n+j} + \vec{U}_{n+j} \right) \, .
\end{align}
\end{subequations}
For extensions of the scheme to orders higher than 2, linearly implicit SBDF will not in general yield a zero-stable scheme (cf. \cite{Prac} where the simultaneous choice of $\mA$ and scheme coefficients must be used to obtain orders $3$ and higher).

\begin{remark}[Alternative derivation via ImEx scheme applied to a singular ODE]\label{Rmk:SngLimit} The scheme \eqref{Eq:MM} may be formally derived by applying an implicit explicit multistep method to a singularly perturbed ODE related to the DAE as follows:
    
    First replace the DAE by the singularly perturbed equation (cf. \cite[Chapter VI.2]{wanner1996solving})
    \begin{align}\label{Eq:SingularSys}
        \vec{w}_t = \vec{f}(\vec{w}, \vec{p}) \qquad  \qquad 
        \epsilon \vec{p}_t = \mat{A} \vec{p} + \big( \vec{g}(\vec{w}, \vec{p}) - \mat{A} \vec{p} \big) \, .
    \end{align}
    The scheme \eqref{Eq:IMEX_DAE} then arises from an implicit-explicit (ImEx) time-discretization of \eqref{Eq:SingularSys} that treats (in the second equation) $\mat{A}\vec{p}$ implicitly and $\vec{g} - \mat{A}\vec{p}$ explicitly, followed by the formal limit $\epsilon \rightarrow 0$.
    As a result, the scheme coefficients $(a_j, b_j, c_j)$ in \eqref{Eq:IMEX_DAE} are guaranteed to be a $p$th order approximation provided they satisfy the order conditions \cite[Equation (7)]{AscherRuuthWetton1995} of a multistep ImEx method. \myremarkend 
\end{remark}

\subsection{Zero-stability for Linear Index-1 Algebraic Equations}\label{subsec:zerostability}
In this subsection we establish, in the spirit of absolute stability, conditions for zero-stability of \eqref{Eq:IMEX_DAE} applied to linear time-independent problems. Although we focus on constant time-independent problems, the result provides a guide for ensuring zero-stability of the equations \eqref{Eq:MM}. 
We say that the scheme is \emph{zero-stable} if the dynamics \eqref{Eq:IMEX_DAE} with $\Delta t = 0$ yield bounded iterates; zero-stability is a necessary condition for convergence.  

Herein, assume a linear constraint $\vec{g}(\vec{w}, \vec{p}) = \mat{G} \vec{p} - \vec{w}$, with invertible matrix $\mat{G} \in \mathbb{R}^{m \times m}$.  The scheme \eqref{Eq:IMEX_DAE} with $\Delta t = 0$ reads:
\begin{subequations}\label{Eq:ZeroStability}
\begin{align}\label{Eq:ZeroS1}
	\sum_{j=0}^{s} a_j \vec{w}_{n+j} &= 0 \, ,  \\ \label{Eq:ZeroS2}
        \sum_{j = 0}^{s} c_j \mat{A} \vec{p}_{n+j} &= 
        \sum_{j=1}^{s-1} b_j \big( \mat{A} \vec{p}_{n+j} - \mat{G}\vec{p}_{n+j} + \vec{w}_{n+j} \big) \,,
\end{align}
\end{subequations}
where $(\vec{p}_j, \vec{w}_j)$ for $j = 1, \ldots, s$ are given as initial data.  Next, introduce the following polynomials defined in terms of the time-stepping coefficients
\begin{align} \label{Eq:PolynomialCoeff}
	a(z) = \sum_{j = 0}^s a_j z^j \,, \qquad
	b(z) = \sum_{j = 0}^{s-1} b_j z^j \,, \qquad
	c(z) = \sum_{j = 0}^{s} c_j z^j \, .
\end{align}
Equation \eqref{Eq:ZeroS1} defines a stable linear recurrence for $\vec{w}_n$ provided $a(z)$ satisfies the \emph{root condition}: the roots of $a(z)$ lie inside the unit disk with only simple roots on the boundary. Every practical linear multistep method (including both AB2, AB3, AB4 and SBDF2) has an $a(z)$ satisfying the root condition. 

With $\vec{w}_n$ completely determined by \eqref{Eq:ZeroS1}, equation \eqref{Eq:ZeroS2} defines a forced linear recurrence for $\vec{p}_n$. The stability of $\vec{p}_n$ then depends on the homogeneous equation \eqref{Eq:ZeroS2} (i.e., $\vec{w}_n \equiv 0$) which takes the form
\begin{align}\label{Eq:ZeroS3}
    \sum_{j = 0}^{s} c_j \mat{A} \vec{p}_{n+j} &= \sum_{j=1}^{s-1} b_j \big( \mat{A} \vec{p}_{n+j} - \mat{G}\vec{p}_{n+j} \big) \,.
\end{align}
A family of eigenvector solutions to \eqref{Eq:ZeroS3} then follow as: let $(\vec{\bar{p}}, \lambda)$ satisfy
\begin{align}\label{Eq:GenEV}
    \lambda \mat{A} \vec{\bar{p}} = \mat{G} \vec{\bar{p}} \,  \qquad (\mat{A} \vec{\bar{p}} \neq 0 ) \,.
\end{align}
Substituting the ansatz $\vec{p}_n = z^n \vec{\bar{p}}$ into \eqref{Eq:ZeroS3} yields 
\begin{align}\label{Eq:CharEq}
    c(z) - b(z) + \lambda b(z) = 0 \, .
\end{align}
A sufficient condition for stability of \eqref{Eq:ZeroS3} is then ensured if each solution of \eqref{Eq:CharEq} has $|z| < 1$. Following in the spirit of absolute stability, introduce the region of zero stability as 
\begin{align}
    \mathcal{D} \equiv \big\{ \lambda \in \mathbb{C} \; \big | \; |z| < 1 \textrm{ holds for each $z$ satisfying } \eqref{Eq:CharEq} \big\} \, . 
\end{align}
Note that $\mathcal{D}$ is a property of the time-stepping coefficients only. One then has the following.

\begin{theorem}[Sufficient conditions for zero-stability of DAEs with linear constraints] 
If $a(z)$ satisfies the root condition, and $\lambda \in \mathcal{D}$ for each $\lambda$ satisfying \eqref{Eq:GenEV}, then the dynamics \eqref{Eq:ZeroStability} are stable.     
\end{theorem}

\begin{figure}[htp!]
    \includegraphics[width=.325\textwidth]{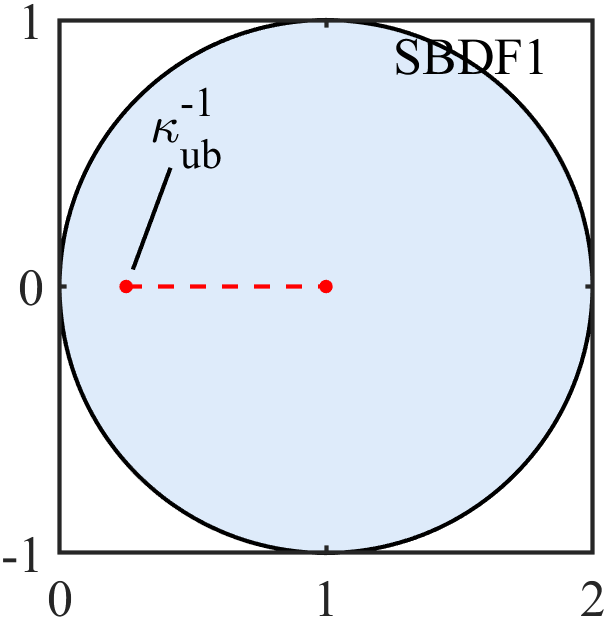}
    \includegraphics[width=.325\textwidth]{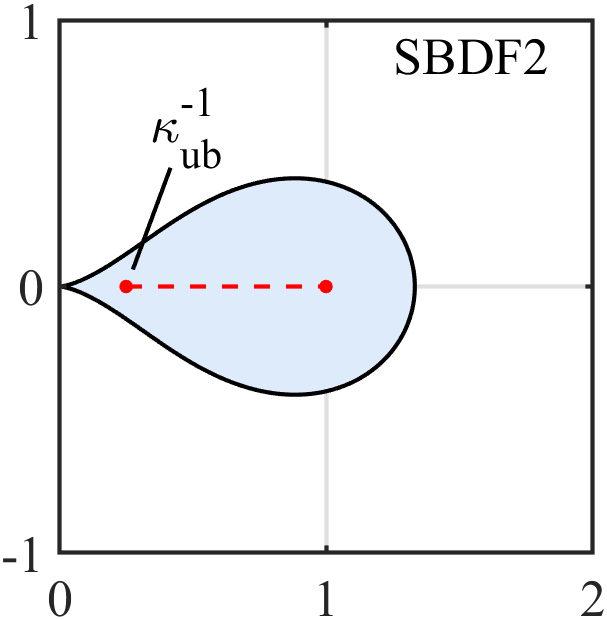}
    \includegraphics[width=.325\textwidth]{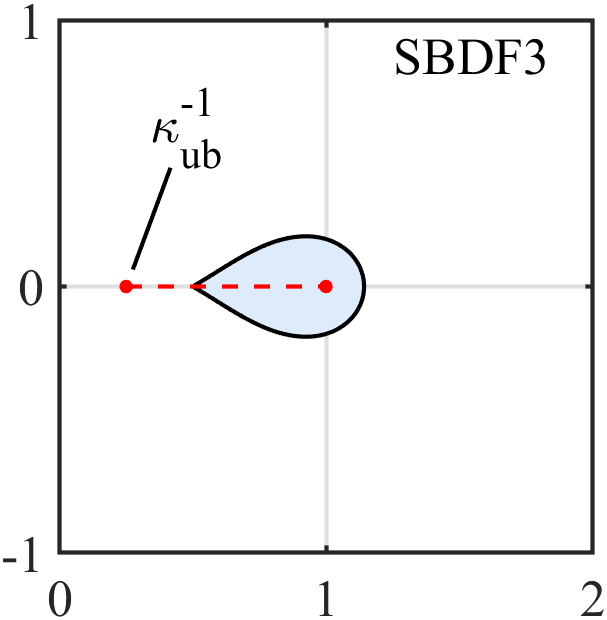}
    \caption{Figures show stability region $\mathcal{D}$ (shaded blue) in relation to the bound on eigenvalues of $\mathcal{A}^{-1}\mathcal{G}_{\eta,h}$ from Theorem~\ref{Thm:CondNumBnd} (red). Note that SBDF1 and SBDF2 are stable, while SBDF3 will not be stable if the largest eigenvalue $\lambda > 2$.}\label{Fig:StabilityRegion}
\end{figure}

Figure~\ref{Fig:StabilityRegion} presents the stability regions $\mathcal{D}$ for SBDF2. For reference, the stability regions for SBDF1 and SBDF3, which are not used in Section~\ref{Sec:numerical_tests}, are shown.  Up to a reflection in the imaginary axis and shift, the sets $\mathcal{D}$ are identical to those characterized in \cite{Theory}. The figures show that SBDF1 and SBDF2 always contains the eigenvalues of $\mA^{-1}\mG_{\eta,h}$. In contrast, depending on the variables $\eta$ and $h$, it is possible for the largest eigenvalue to lie outside the SBDF3 stability region---rendering SBDF3 unstable.  Note that by modifying the SBDF coefficients for orders $3,4$ and $5$ as in \cite{Prac}, it is possible to always devise time-stepping schemes where the eigenvalues of $\mA^{-1}\mG_{\eta,h}$ lie in the stability regions $\mathcal{D}$. We do not explore these approaches here.

\section{Numerical Tests}\label{Sec:numerical_tests}
In this section we numerically test the performance and robustness of $\mA$ as a preconditioner in dimensions $d = 1$ (\S~\ref{subsec:1dpreconditioner}) and $d=2$ (\S~\ref{subsec:2dpreconditioner}). We then demonstrate the accuracy and effectiveness of the time integration schemes to solve the SGN in $d=1$ (\S~\ref{subsec:dswe_1d}) and $d=2$ (\S~\ref{subsec:dswe_2d}).  All numerical computations are performed in \textsc{Matlab}.  For the purposes of presentation, all computations in this section are done with discretizing the constraint system \eqref{Eq:DSWEcf}. Numerical solutions to the \eqref{Eq:dswe} are similar but not shown. 

The following exact solution will be used in various test cases. In free space with a flat bathymetry, $h(\vx) = h_0$, the SGN equations admit a solitary wave solution of a single parameter $a > 0$ \cite{ray76,sg69}:

\begin{align}\nonumber
    \xi(x,y, t) &\equiv x \cos(\theta)+ y \sin(\theta)- c t + \xi_0 \, , \\ \label{ExSleta}
    \eta_{\rm s}(x,y, t) &= h_0 + a \sech^2\big(\gamma \xi(x,y) \big) \, ,
\end{align}
with velocity
\begin{align}\label{ExSlu}
    \vec{u}_{\rm s}(x,y,t) &= w(x,y,t) \, \big(\cos \theta , \sin \theta \big) \, ,\\
    w_{\rm s}(x,y,t) &= c \frac{1}{\eta_{\rm s}(x,y,t)} \big( \eta_{\rm s}(x,y,t) - h_0 \big) \, . 
\end{align}
The speed $c$ and \emph{width} $\gamma$ are given in terms of $a$ and $h_0$ by
\begin{align*}    
    c=\sqrt{h_0+a} \qquad \textrm{and} \qquad 
    \gamma=\frac{\sqrt{3a}}{2 h_0 c} \, .    
\end{align*}
The variables $\theta \in [0, 2\pi)$ and $\xi_0 \in \mathbb{R}$ are independent of $a$ and determine the orientation and center of the wave crest.  In $d=2$, the solitary wave is a one dimensional profile. A solitary wave in one dimension is obtained from the two dimensional wave by taking $\theta = 0$, $y = 0$ and a velocity $w$.

\subsection{Spatial Discretizations}\label{Subsec:Spatialdiscretization}
Throughout this section we adopt a domain with periodic (in $d = 1$) or doubly-periodic ($d=2$) boundary conditions, e.g., 
\begin{align} 
	\label{BC:periodic} 
	\eta(x, y + L_y) = \eta(x + L_x, y) = \eta(x, y)\;, 
	\qquad
	\vu(x, y+ L_y) = \vu(x + L_x, y) = \vu(x, y)\;,   
\end{align}
where  $\Omega = (0, L_x)\times (0,L_y)$ with $L_x$ and $L_y$ being the periods in the $x$- and $y$-directions, respectively.  We then discretize space with an equispaced grid of $n$ points in each direction, e.g., in dimension $d = 1$, the grid spacing $\Delta x$ and mesh points are 
\begin{align*}
    \Delta x = \frac{L}{n} \qquad x_j = j \, \Delta x \qquad 0 = x_0 < x_1 < \ldots < x_{n-1} = L - \Delta x \, .    
\end{align*}
Functions $\vu(x)$ are discretized as $\vu(x_j) \approx u_j$ with 
\begin{align*}
    \vec{u} = (u_0, u_1, \ldots, u_{n-1})^T \in \mathbb{R}^n \, .
\end{align*}
Function derivatives are obtained spectrally as 
\begin{align*}
    u_x(x_j) \approx  (  \mathbf{D} \mathbf{u} )_j \qquad \textrm{and} \qquad
    \mathbf{D} = \mathbf{F}^{-1} (\imath \mathbf{K}) \mathbf{F} \, , 
\end{align*}
where $\mathbf{F}$ is the discrete Fourier transform matrix and $\mathbf{K}$ is a diagonal matrix 
consisting of discrete wavenumbers of the form 
\begin{align*}
    k_j = \frac{2\pi}{L} j \qquad j = -\frac{n}{2}+1, \ldots, \frac{n}{2}-1 \, .
\end{align*}
We adopt the convention that the $n/2$ wavenumber is $0$.  The operators $\mG_{\eta, h}$ and $\mA$ are then discretized with spectrally accurate approximations (we suppress the parameter dependence of $\eta$ and $h$ on $\mathbf{G}$)
\begin{align}\label{def:discrete_operators}
    (\mG_{\eta, h} \vu) (x_j) \approx (\mathbf{G} \mathbf{u})_j \qquad
    (\mA \vu) (x_j) \approx (\mathbf{A} \mathbf{u})_j \qquad (j = 0, \ldots, n-1)\,.
\end{align}
For example, the matrix $\mathbf{A}$ has the form
\begin{align*}
    \mathbf{A} = \sigma^\star \mathbf{I} + \alpha^\star \mathbf{D}^{\dagger} \mathbf{D} \, ,
\end{align*}
where $\mathbf{D}^{\dagger} = \overline{\mathbf{D}}^T$ is the conjugate transpose of the derivative matrix $\mathbf{D}$. By construction, both $\mathbf{A}$ and $\mathbf{G}$ are real symmetric and positive definite. 

In presenting numerical results, we make use of both the (discrete) $2-$ and $\infty-$norms
\begin{align*}
    \| \mathbf{u} \|_{2, \Delta x}^2 = \Delta x \, \sum_{j = 1}^n u_i^2 
    \qquad \qquad
    \| \mathbf{u} \|_{\infty} = \max_{0 \leq i \leq n-1} |u_i| \, .
\end{align*}
We also use the $2-$norm induced by a matrix $\mathbf{M} \succ 0$ (which will generally be set to $\mathbf{G}$)
\begin{align*}
    \| \mathbf{u} \|_{2,  \Delta x, \mathbf{M}}^2 \equiv \left(\mathbf{u}^{T} \mathbf{M} \mathbf{u}\right) \, \Delta x  \qquad (\mathbf{u} \in \mathbb{R}^n) \, .
\end{align*}
With these notations, note that  $\| \mathbf{u} \|_{2, \Delta x}^2 =  \| \mathbf{u} \|_{2, \Delta x, \mathbf{I}}^2$. For the equation
\begin{align}\label{eq:matrixequation}
    \mathbf{G} \mathbf{u} = \mathbf{b} \, ,
\end{align}
let $\mathbf{u}^* \equiv \mathbf{G}^{-1} \mathbf{b}$ be the exact solution. In convergence studies where preconditioned conjugate gradient is used to solve \eqref{eq:matrixequation}, we track the error 
\begin{align}\label{def:discrete_Gerr}
    \varepsilon(\mathbf{u}) \equiv \frac{\| \mathbf{u} - \mathbf{u}^* \|_{2, \Delta x, \mathbf{G}}}{\| \mathbf{b} \|_{2, \Delta x}} = 
    \frac{\| \mathbf{G}\mathbf{u} - \mathbf{b} \|_{2, \Delta x, \mathbf{G}^{-1}}}{\| \mathbf{b} \|_{2, \Delta x}} \,.
\end{align}
Here $\varepsilon$ measures the discrete error of $\mathbf{u}$ (in the $\mathbf{G}$ norm), or equivalently the discrete residual $\mathbf{G}\mathbf{u} - \mathbf{b}$ (in the $\mathbf{G}^{-1}$ norm).   Note that we use $\varepsilon$ as a primary measure of the PCG error, as opposed to the $2$-, or $\infty$-norm, since PCG minimizes $\varepsilon$ over progressively larger Krylov subspaces. This ensures that $\varepsilon$ is non-increasing with each iteration of PCG.  

The generalized eigenvalues of $\mathbf{G}$ and $\mathbf{A}$ are then 
\begin{align}\label{Eq:DiscreteGenEigvalue}
    \mathbf{G} \mathbf{u}_k = \lambda_k \mathbf{A} \mathbf{u}_k \qquad (\mathbf{u}_k \neq 0) \, 
\end{align}
and, without loss of generality, can be ordered as
\begin{align*}
    0 < \lambda_1 \leq \lambda_2 \leq \ldots \leq \lambda_n \, .
\end{align*}
For dimension $d = 2$, we extend all the definitions in the section in the natural way via tensor products. 

\subsection{Preconditioner Tests in 1D}\label{subsec:1dpreconditioner}
In this section we test the performance of the preconditioner $\mA$ with coefficients \eqref{Eq:FlatbathCoeff}--\eqref{Eq:Optcoeff} using PCG to solve for $u$ on the test problem ($d = 1$, $\Omega = [0, 1]$)
\begin{align}\label{Eq:Testproblem}
    \mG_{\eta,h} u = b \qquad \textrm{where} \qquad b(x) = \cos(4 \pi x) \, . 
\end{align}
Here the variable coefficient functions are taken to be
\begin{align}\label{Eq:Testproblem1_var}
    h(x) &= 1 + \hvar \exp\left( - (x - 1/2)^2/\sigma^2 \right) \, , \qquad \sigma = 1/20 \\
    \eta(x) &= 1 + \etavar \cos(4\pi x)^2 \, . 
\end{align}
With these functions, we then vary the parameters $\etavar, \hvar \, \in (-1, \infty)$ which control the water depth contrast ratio, and bathymetry gradients respectively ($\hvar = 0$ yielding a flat bathymetry).  The preconditioner performance bound then depends on
\begin{align*}
    \frac{\eta_{\rm max}}{\eta_{\rm min}} = \etavar + 1 \qquad \textrm{and} \qquad
    |h_x|_{\rm max} = \frac{\sqrt{2}}{e^{1/2} \sigma} |\hvar| \approx 17.2 |\hvar| \, .
\end{align*}
When $\eta_0$ is \emph{large} relative to $|h_x|_{\rm max}$, $\gb \sim (\eta_0 + 1)^3$. Alternatively, when $|h_x|_{\rm max}$ is large relative to $\eta_0$, $\gb \sim (17.2)^2 \eta_0 h_0^2$.

The discrete form of \eqref{Eq:Testproblem} then is
\begin{align}\label{Eq:DiscreteGEq}
    \mathbf{G} \mathbf{u} = \mathbf{b} \qquad \textrm{where} \qquad \mathbf{b}_j = b(x_j) \,
\end{align}
with $\mathbf{G}$ defined in \eqref{def:discrete_operators}.  

\begin{figure}[htp!]
    \includegraphics[width=.43\textwidth]{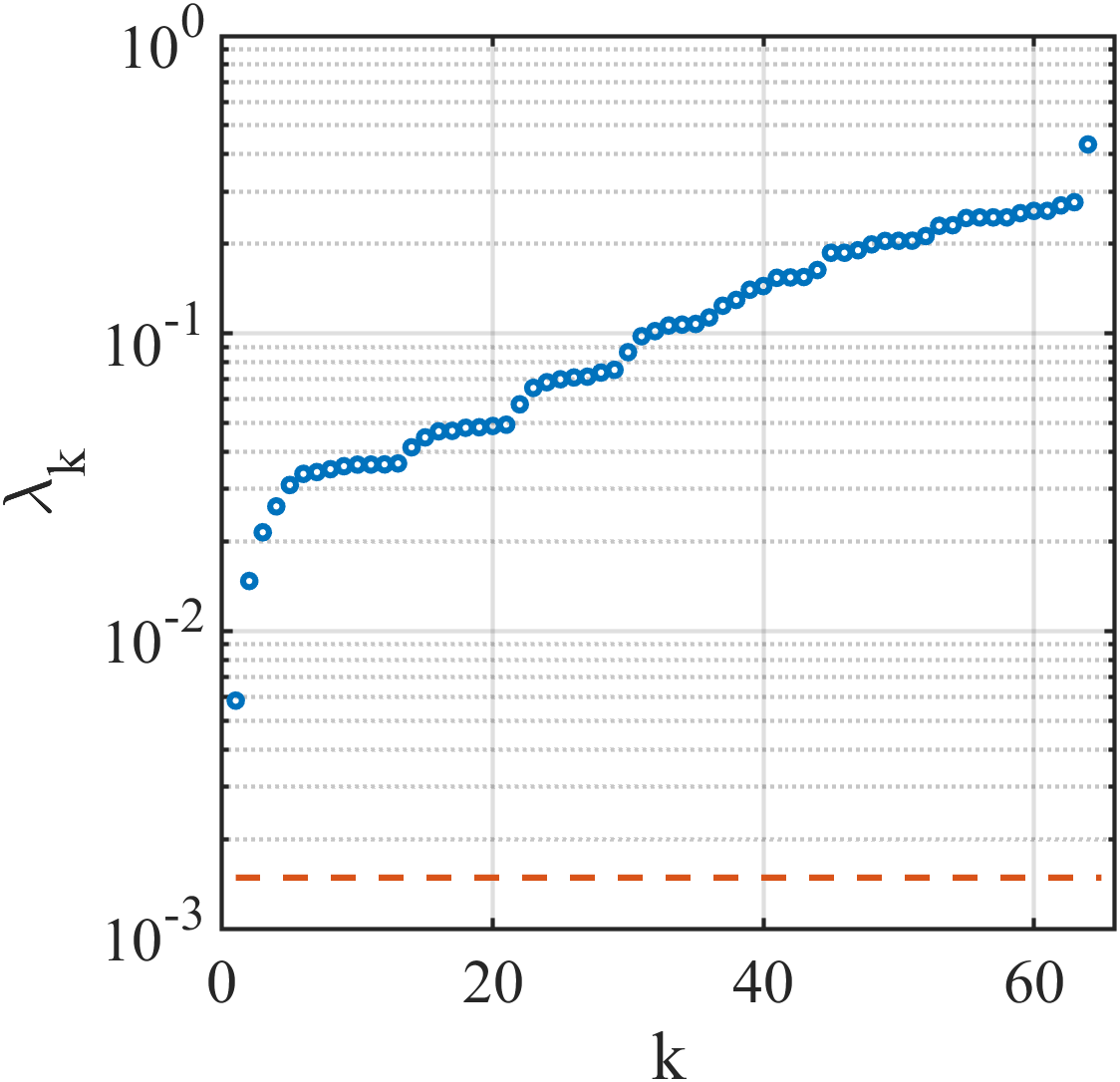}
    \includegraphics[width=.43\textwidth]{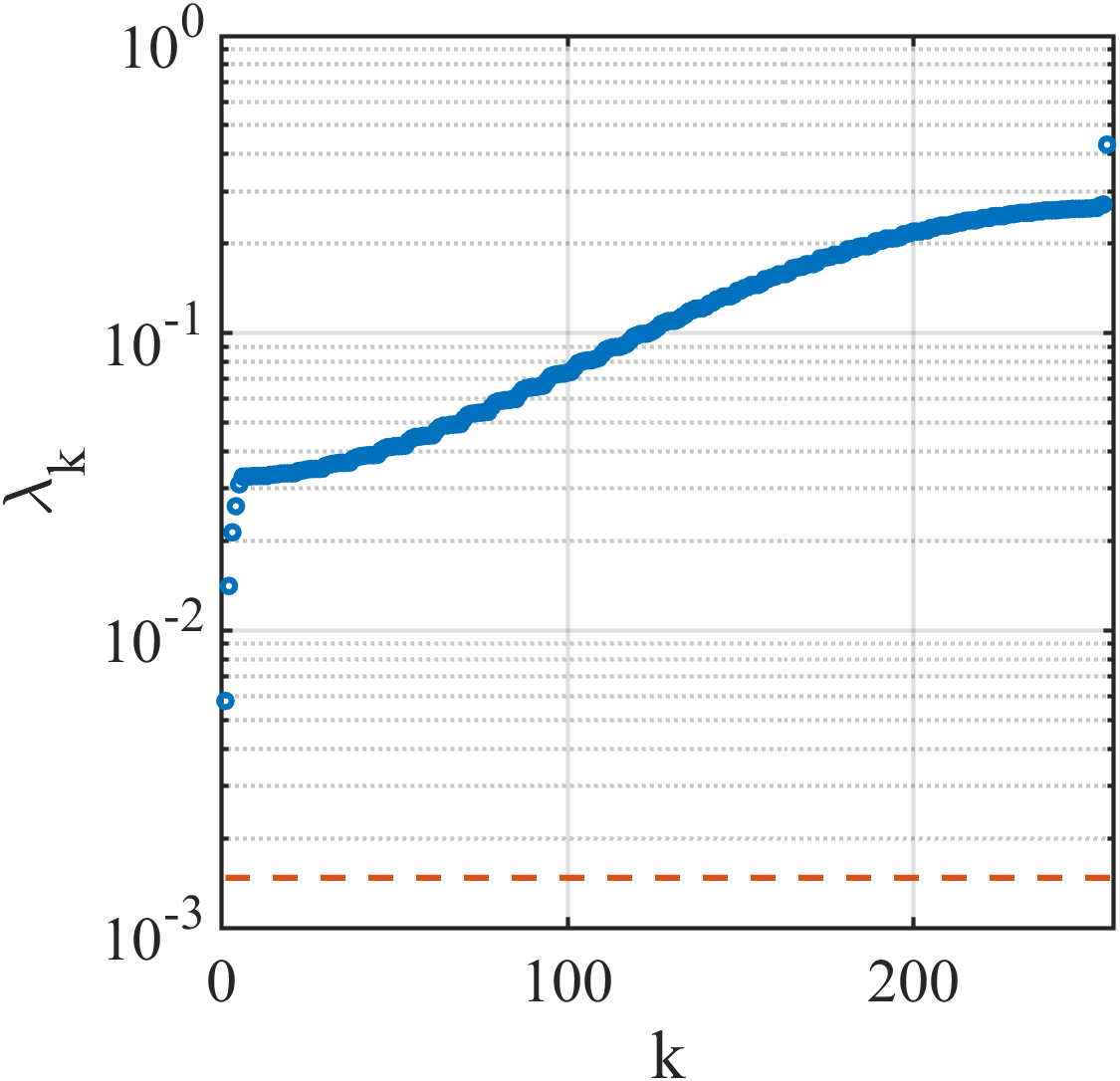}
    \caption{Plot of the generalized eigenvalues $\lambda_k$ versus $k$ for the 1d test problem \eqref{Eq:Testproblem1_var} with $\eta_0 = h_0 = 1$ and mesh resolutions $n=64$ (left), $n=256$ (right), as well as the lower bound $\kapub^{-1}$ (dashed line). The figures provide numerical evidence that the generalized eigenvalues of the operators $\mG_{\eta,h}$ and $\mA$ have a continuous spectrum confined to an interval $\sim [3 \times 10^{-2}, 3 \times 10^{-1}]$, with $\sim 6$ isolated eigenvalues outside this interval.}\label{Fig:1d_gen_eigvalues}
\end{figure}

When solving \eqref{Eq:DiscreteGEq}, we initialize PCG with $\mathbf{u}^{(0)} = \mathbf{0}$, so that the initial error is $\varepsilon(\mathbf{0})$.  In addition to an upper bound on the continuum operators $\kappa_{\mA}( \mG_{\eta,h}) \leq \kapub$, it can also be shown that the same upper bound holds for the discrete matrices $\mathbf{A}, \mathbf{G}$ defined by \S~\ref{Subsec:Spatialdiscretization}, i.e., $\kappa(\mathbf{A}^{-1} \mathbf{G}) \leq \kapub$ (see Remark~\ref{rmk:meshInd}).  As a result, if $\{ \mathbf{u}^{(j)}\}_{j \geq 0}$ denotes the iterates of PCG, then \eqref{eq:CGRate} yields the upper bound
\begin{align}\label{Eq:pcg_upperbound}
    \varepsilon\left(\mathbf{u}^{(j)}\right) \leq 2 \left( \frac{\sqrt{\kapub} - 1}{\sqrt{\kapub} + 1} \right)^j \varepsilon( \mathbf{0}) \, . 
\end{align}
Figure~\ref{Fig:1d_gen_eigvalues} plots the eigenvalues \eqref{Eq:DiscreteGenEigvalue} for the choice of parameters $\eta_0 = h_0 = 1$ and two mesh sizes---a coarse ($n=64$) and fine ($n=256$) one. In both cases, the eigenvalues lie in the interval $[\kapub^{-1}, 1]$, as expected by Theorem~\ref{Thm:CondNumBnd}.

To test the robustness and performance of the preconditioner $\mathbf{A}$, we use PCG to solve \eqref{Eq:DiscreteGEq} for increasingly challenging problems generated from large gradients in bathymetry (varying $\hvar$) and large contrast ratio (varying $\etavar$).  As discussed in \S~\ref{Sec:BoundonG} large values of $\etavar, \hvar$ yield two independent mechanisms for increasing the conditioning number $\kappa_{\mA}(\mG_{\eta,h})$. 

\begin{figure}[htp!]
    \includegraphics[width=.325\textwidth]{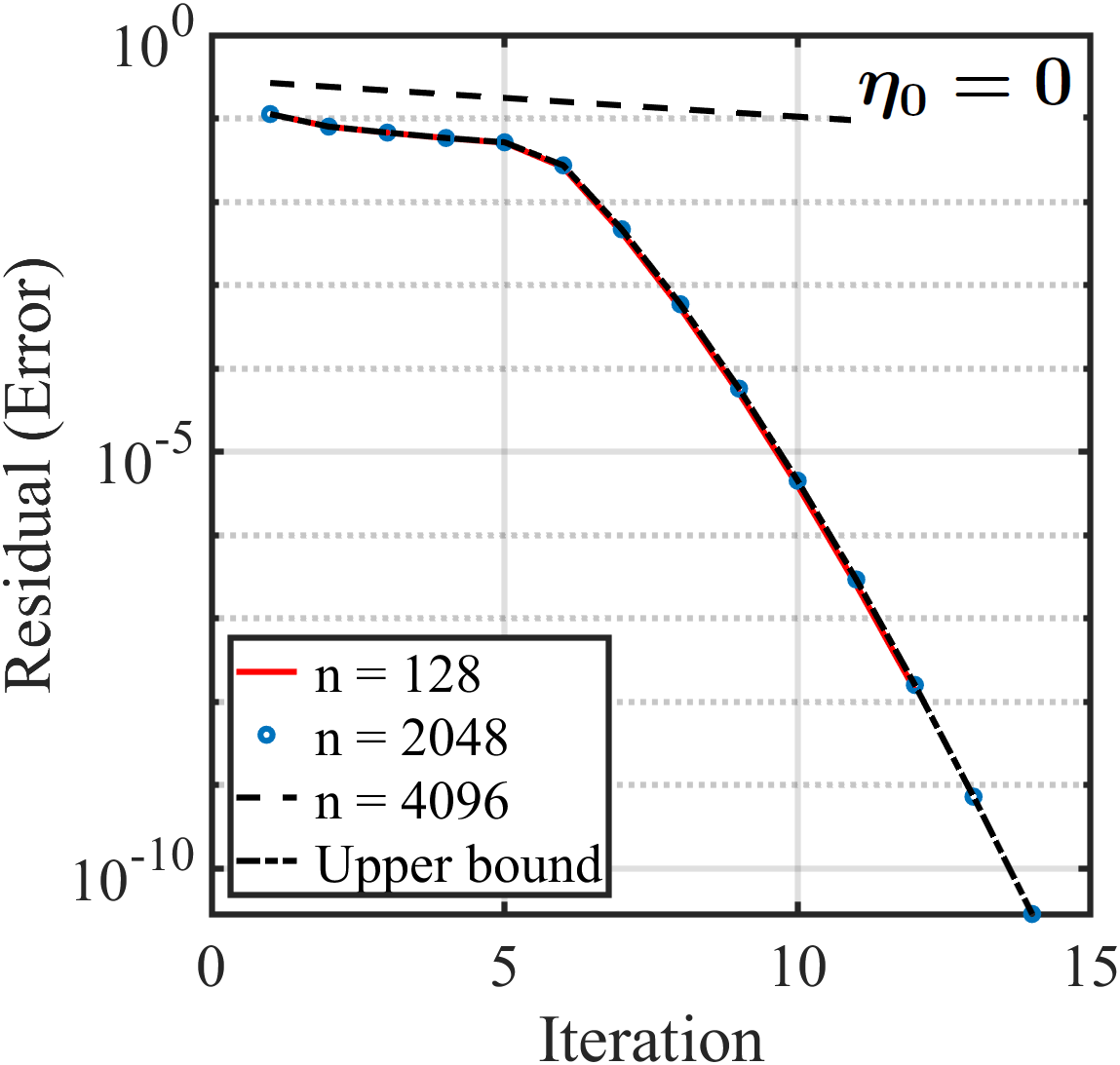}
    \includegraphics[width=.325\textwidth]{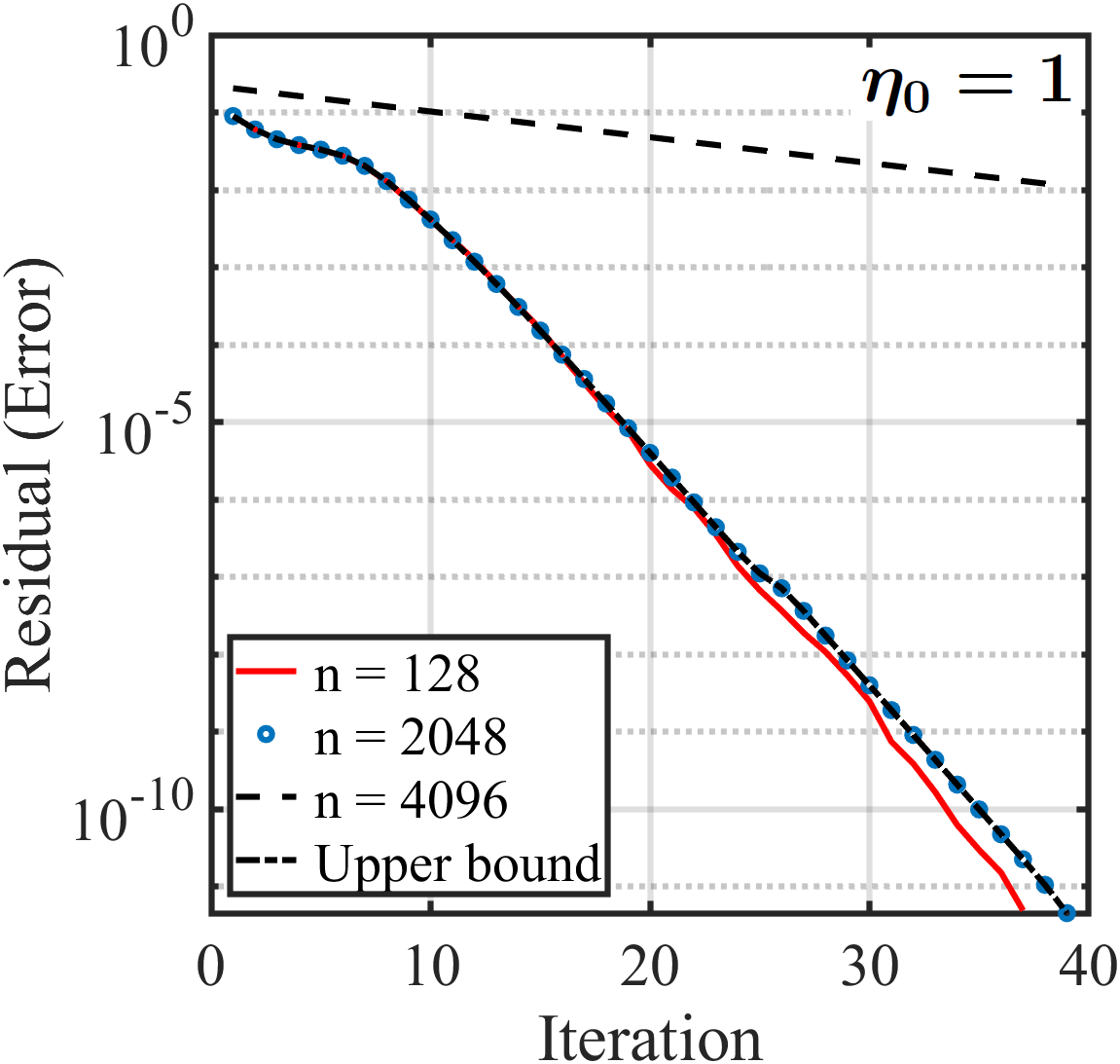}
    \includegraphics[width=.325\textwidth]{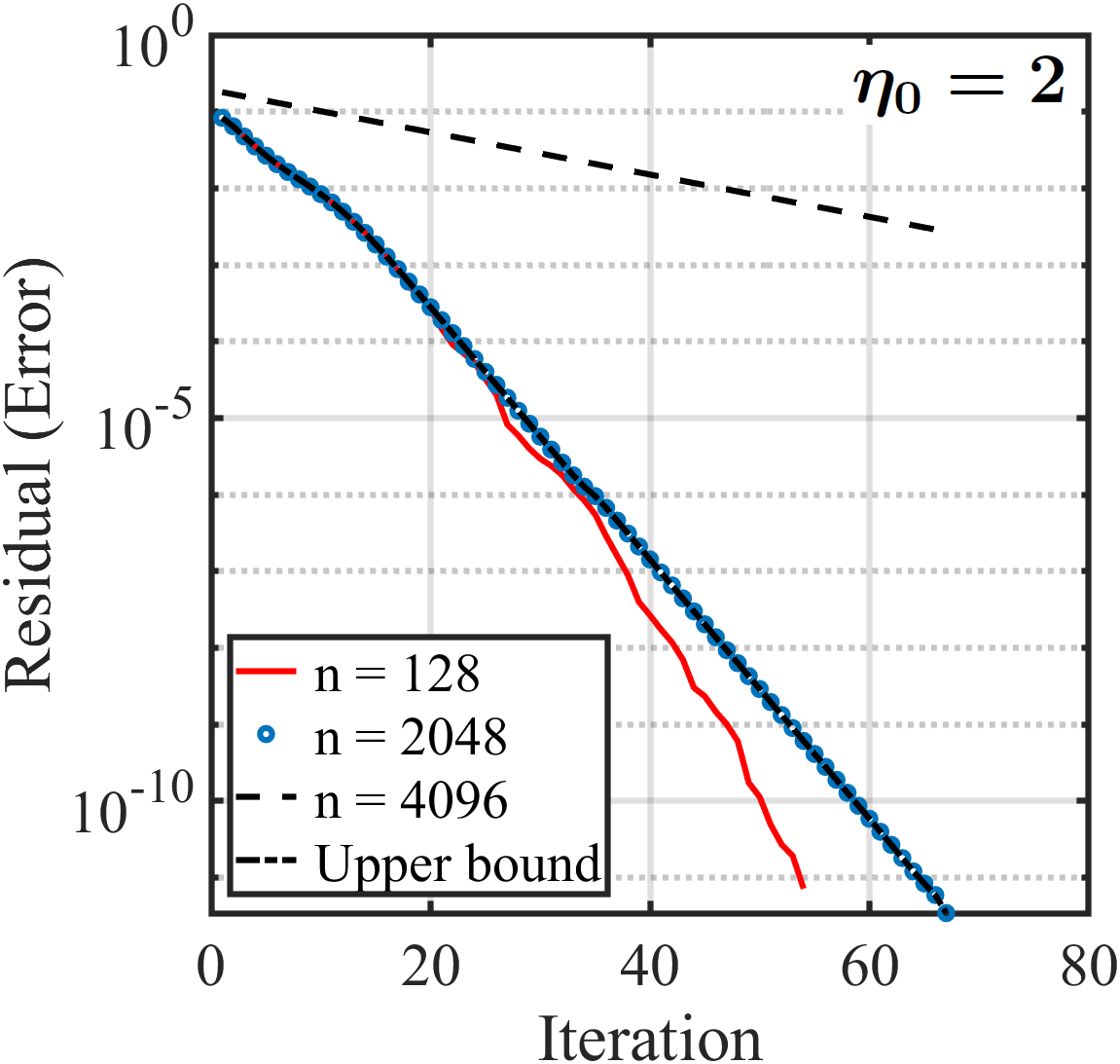}\\
    \includegraphics[width=.325\textwidth]{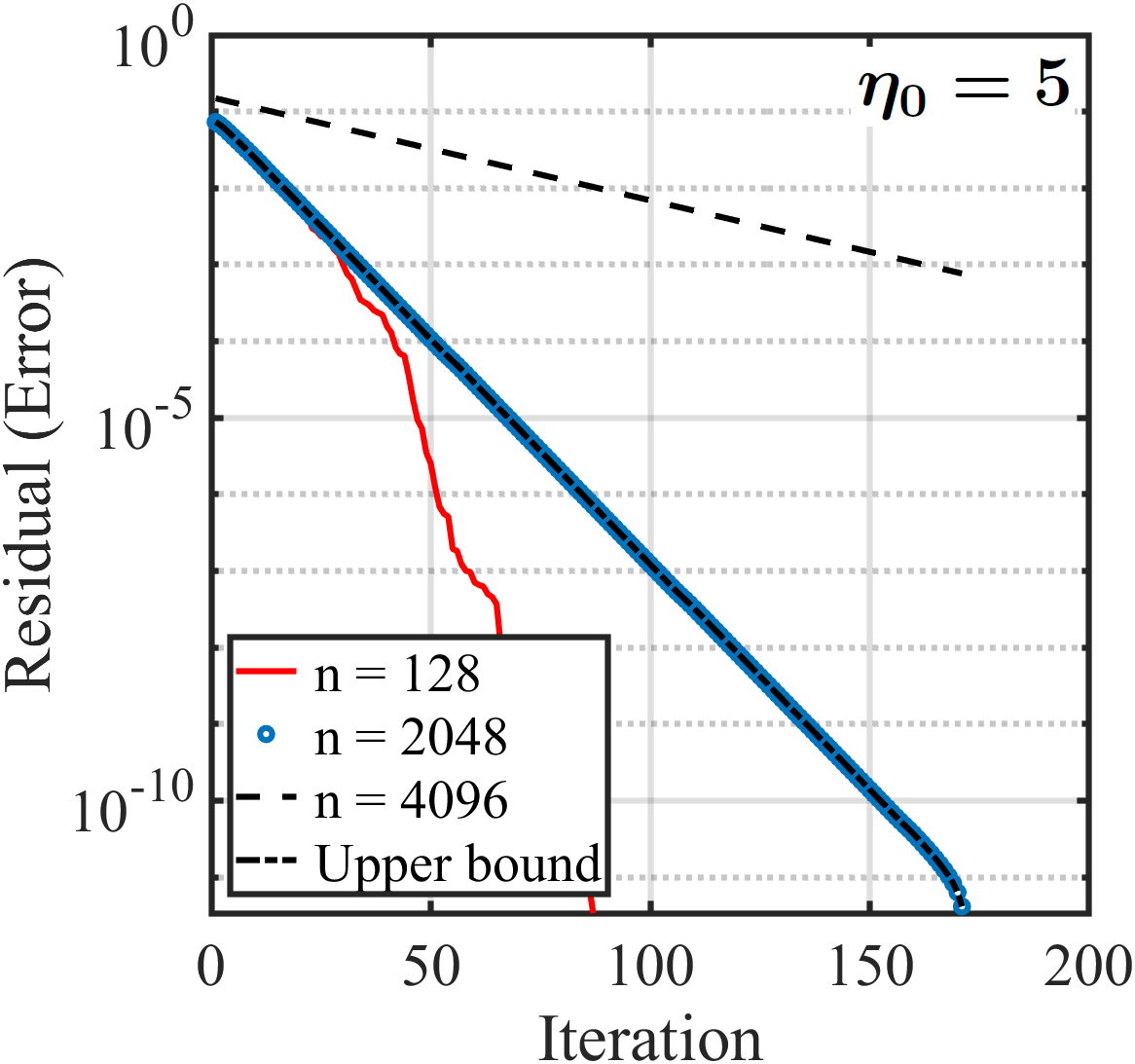}
    \includegraphics[width=.325\textwidth]{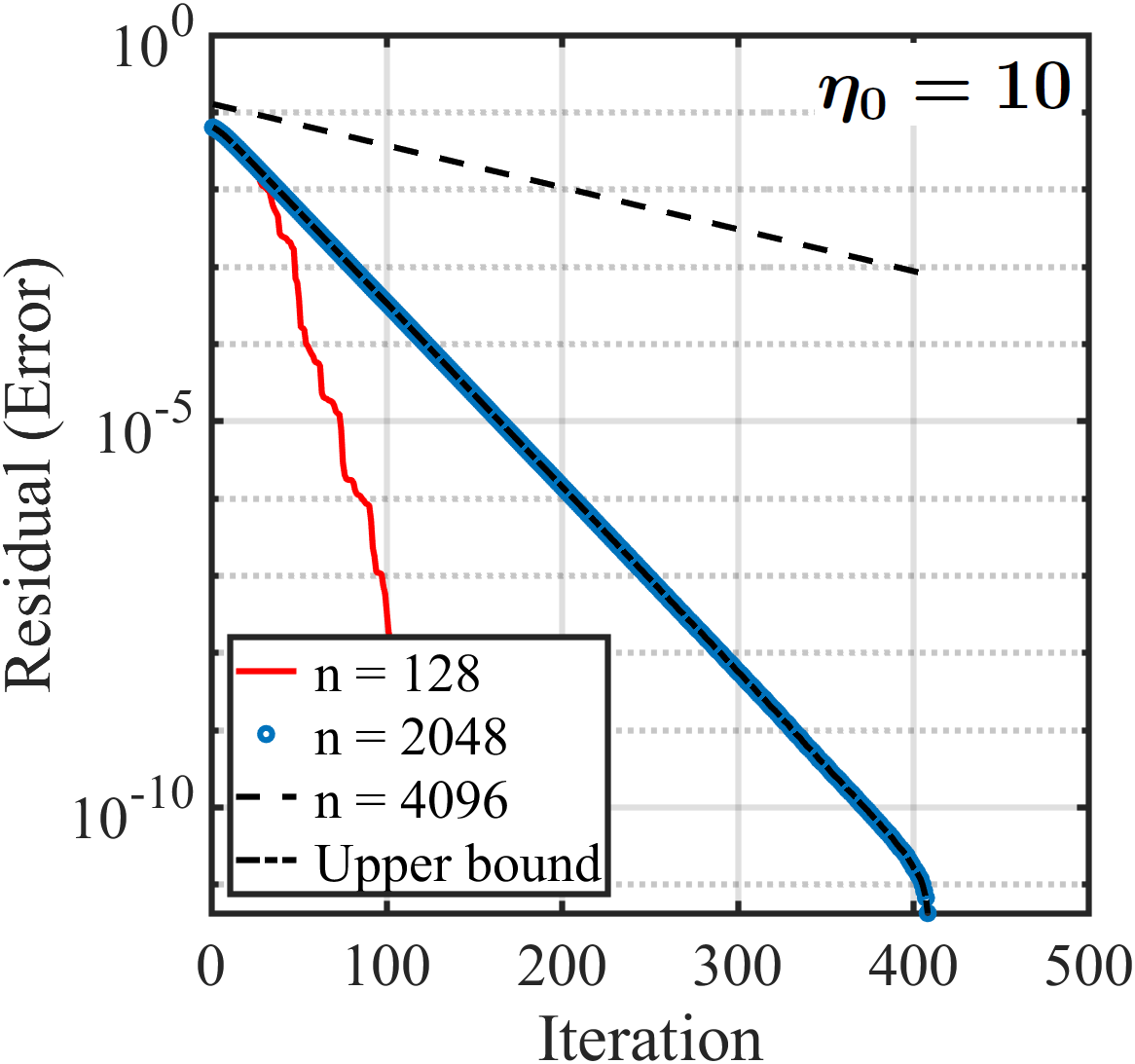}
    \includegraphics[width=.325\textwidth]{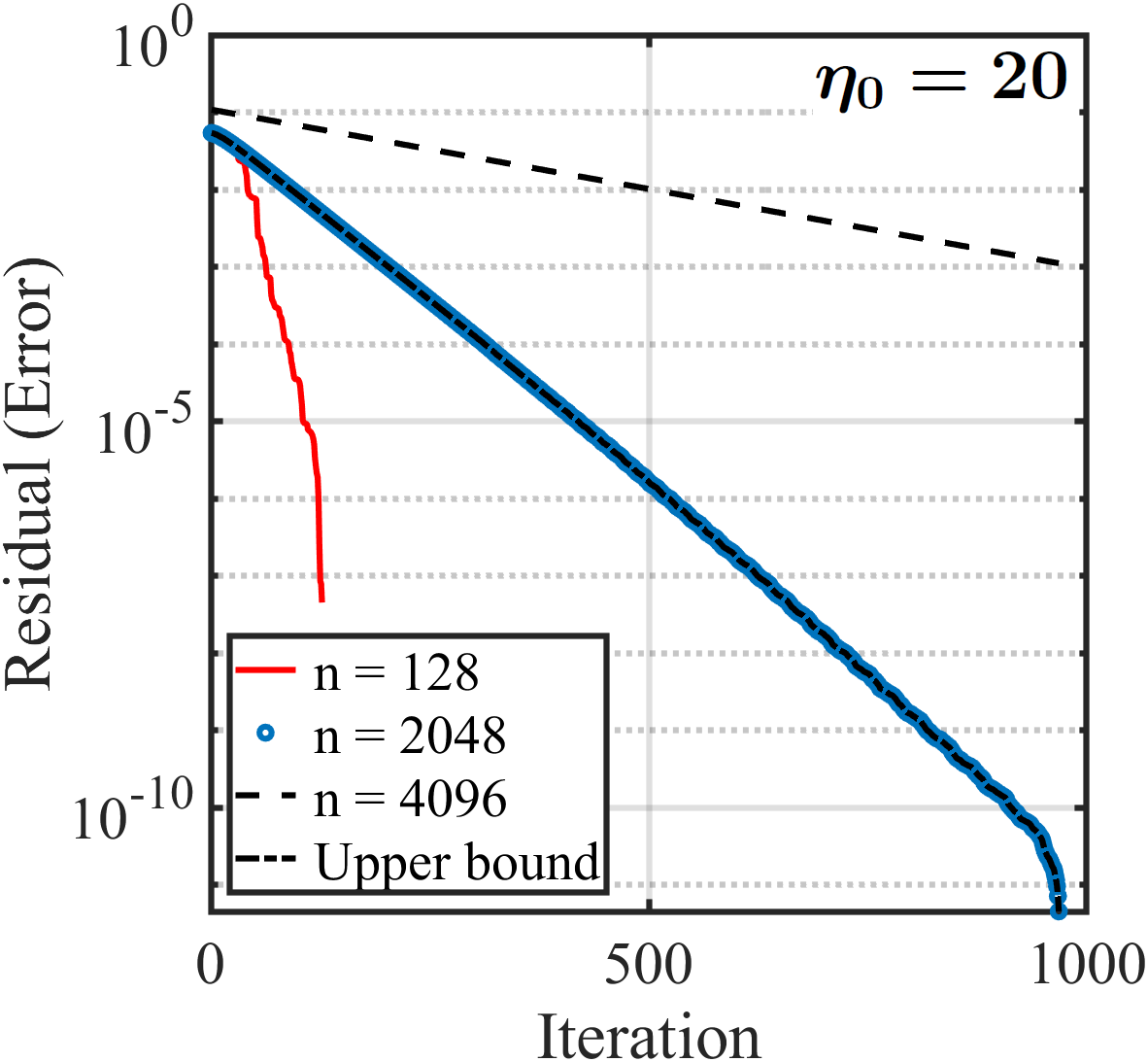}
    \caption{PCG residual error $\varepsilon$ (defined in \eqref{def:discrete_Gerr}) versus iteration for the 1d problem \eqref{Eq:Testproblem1_var}. Subfigures show the performance with increasing water depth contrast ratios $\eta_0$, holding $h_0 = 1$ fixed. The iterations increase with $\eta_0$ but remain bounded with increasing mesh resolution ($n$). The upper bound (dashed line) is \eqref{eq:CGRate} using $\kapub$.}\label{Fig:1d_pcg_cvg_eta0test}
\end{figure}

\begin{figure}[htp!]
    \includegraphics[width=.325\textwidth]{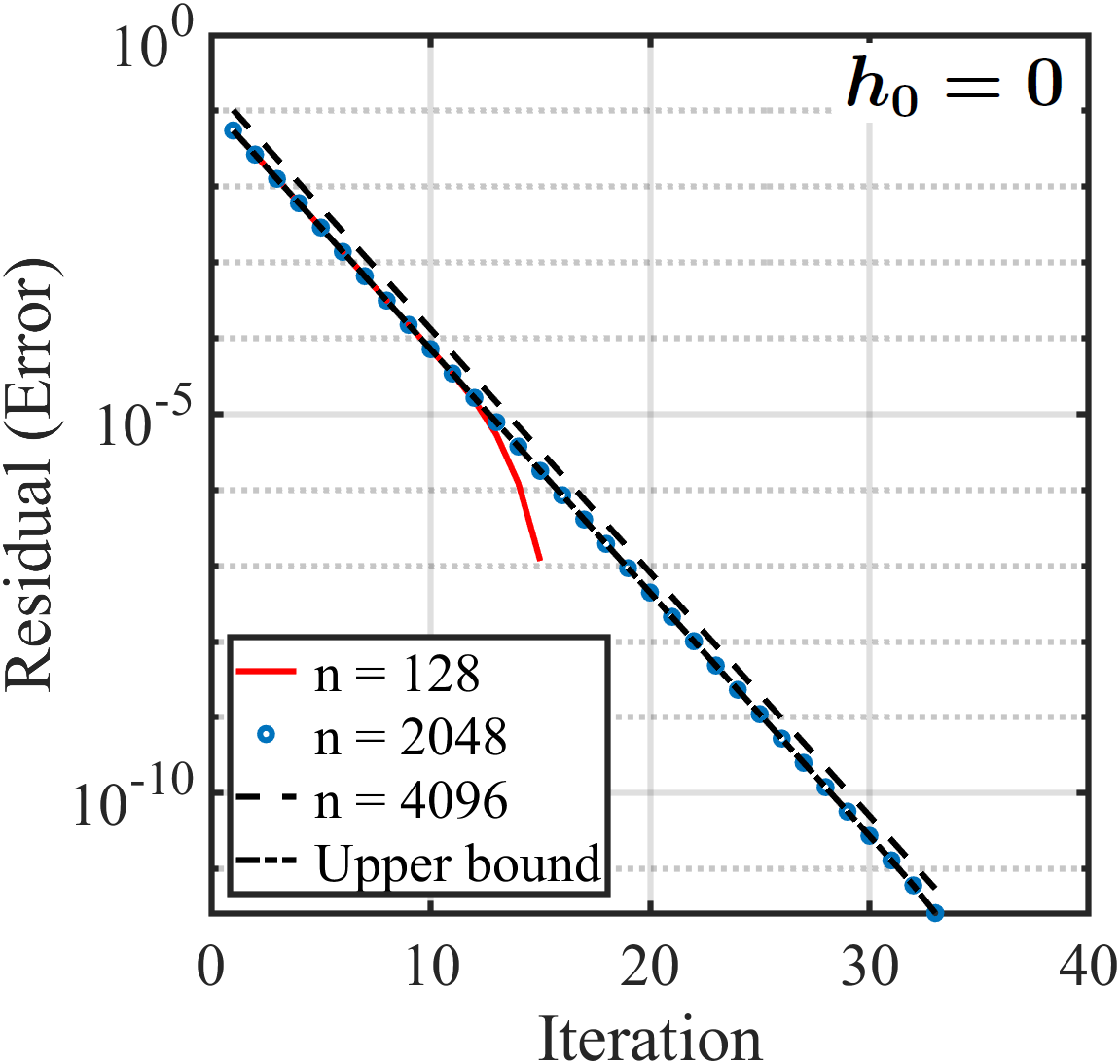}
    \includegraphics[width=.325\textwidth]{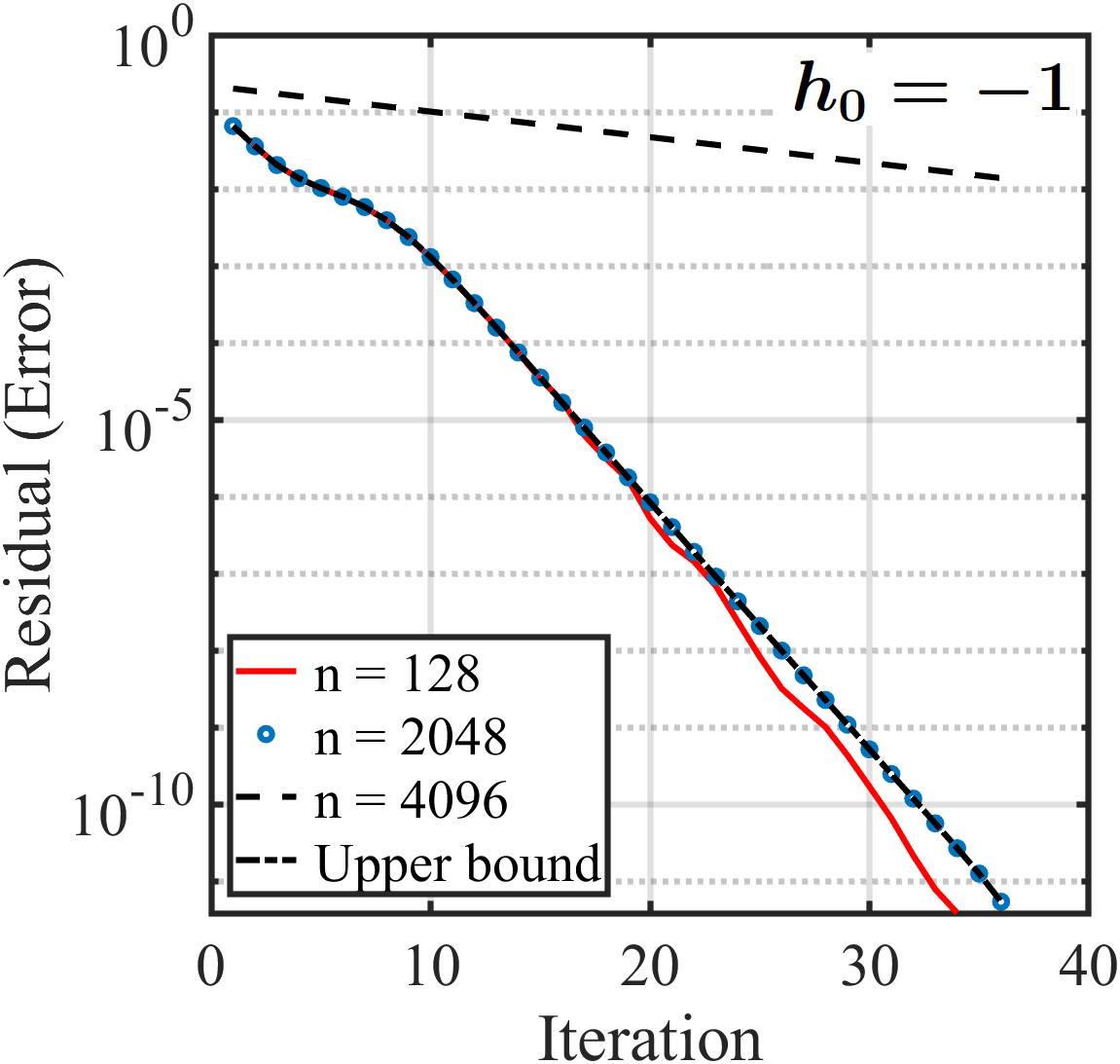}
    \includegraphics[width=.325\textwidth]{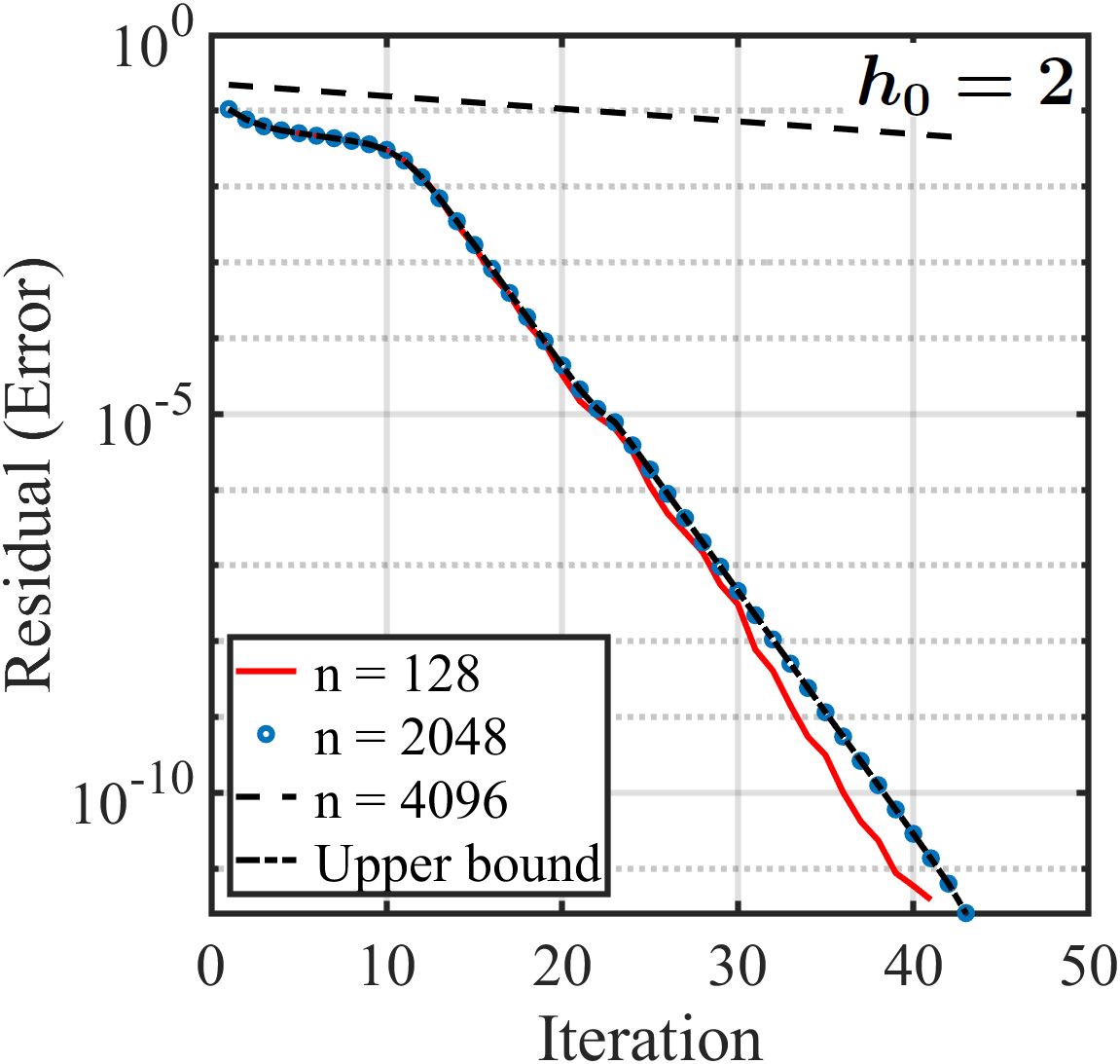}\\
    \includegraphics[width=.325\textwidth]{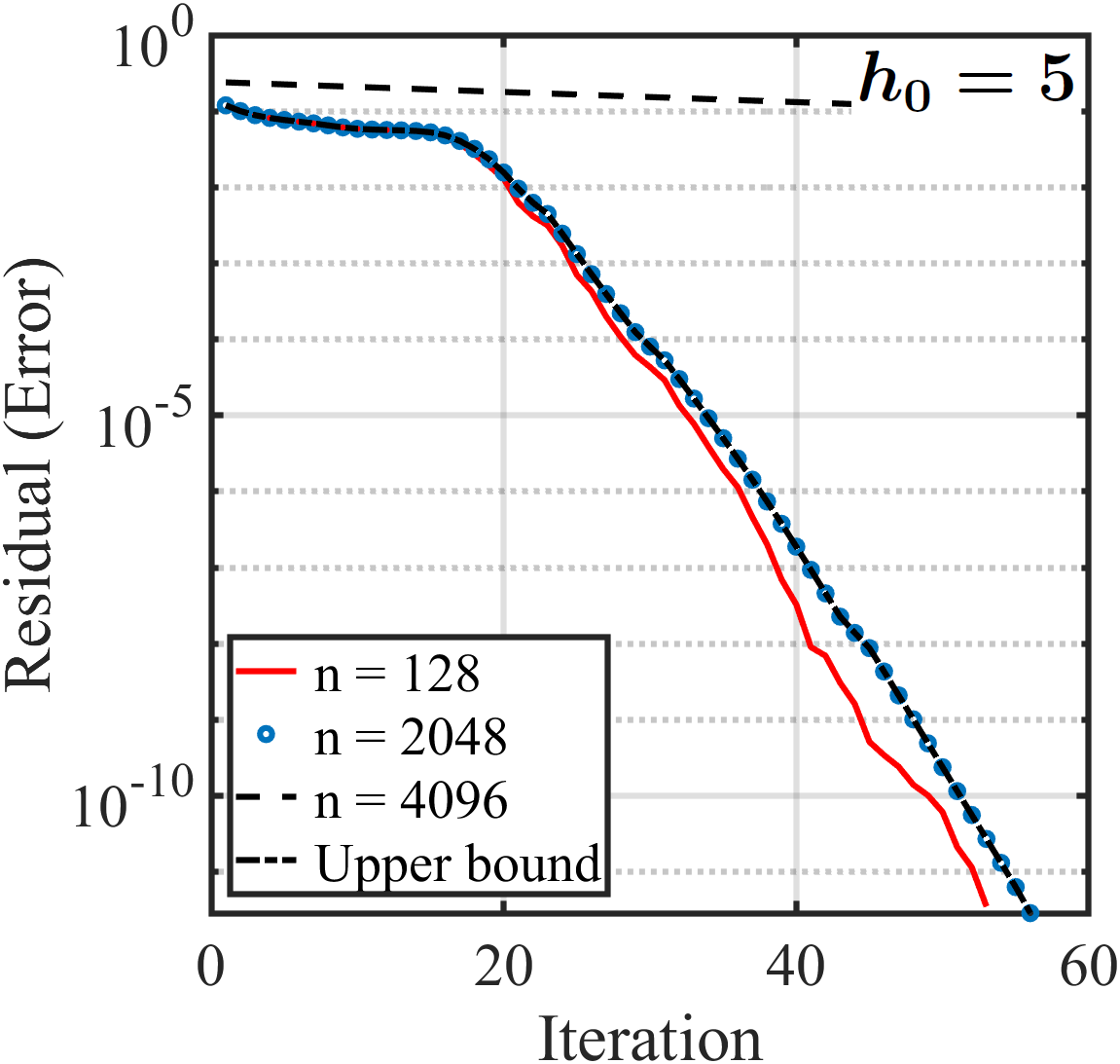}
    \includegraphics[width=.325\textwidth]{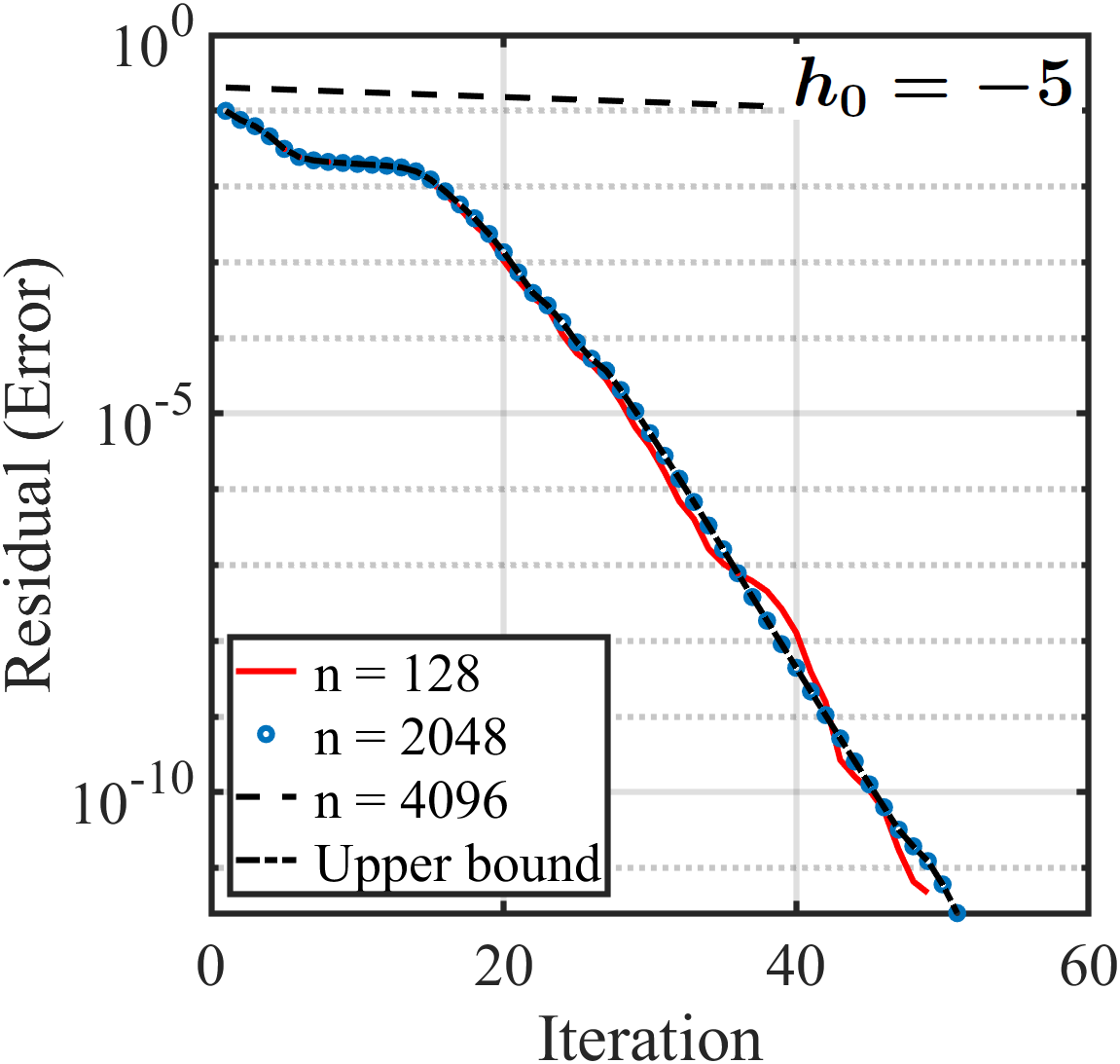}
    \includegraphics[width=.325\textwidth]{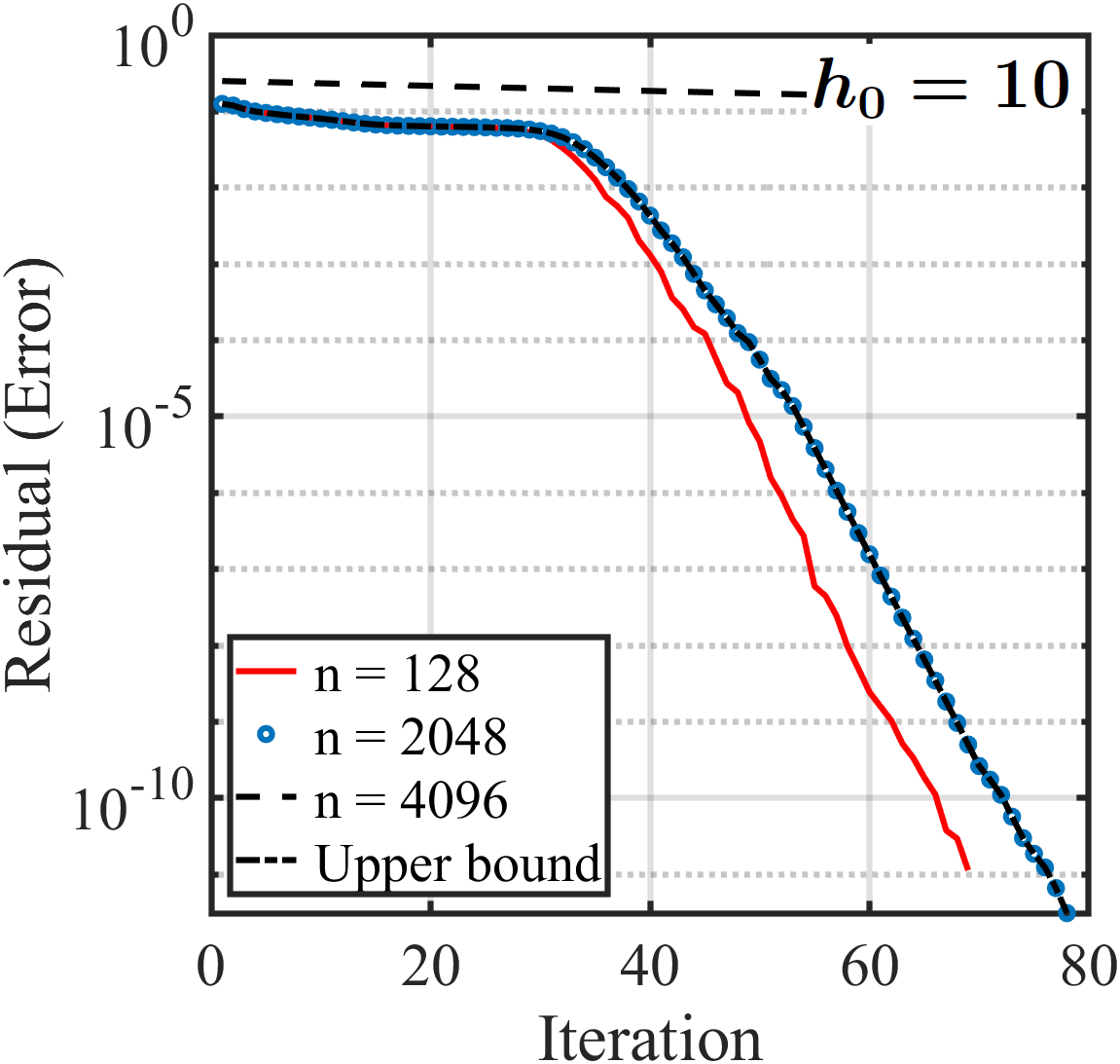}
    \caption{PCG residual error $\varepsilon$ (defined in \eqref{def:discrete_Gerr}) versus iteration for the 1d problem \eqref{Eq:Testproblem1_var}.  Subfigures show the performance with increasing bathymetry gradients $h_0$ (holding $\eta_0 = 1$ fixed).  Again, the iterations increase with $h_0$ but remain bounded with increasing mesh resolution ($n$). The upper bound (dashed line) is \eqref{eq:CGRate} using $\kapub$. }\label{Fig:1d_pcg_cvg_h0test}
\end{figure}

Figures~\ref{Fig:1d_pcg_cvg_eta0test} and \ref{Fig:1d_pcg_cvg_h0test} plot the error $\varepsilon(\mathbf{u}^{(i)})$ and upper bound \eqref{Eq:pcg_upperbound} versus iteration $i$.  When $\vec{b}$ is specified as problem data, we use the fully converged PCG output as a reference solution for $\vu^*$ to compute $\varepsilon$.  Figure~\ref{Fig:1d_pcg_cvg_eta0test} examines the effect of increasing $\etavar$ holding $\hvar$ fixed, while Figure~\ref{Fig:1d_pcg_cvg_h0test} varies $\hvar$ holding $\etavar$ fixed.  In both figures, the number of PCG iterations increase with $\etavar, \hvar$ (which is expected). However, the figures also highlight that $\mathbf{A}$ remains an effective preconditioner over a wide range of $\etavar, \hvar$, with little to no effect due to the mesh resolution $n$. Lastly, it is worth noting that the PCG outperforms the (apriori) upper bound \eqref{eq:CGRate}, and is partially explained in the subsequent remark.

\begin{remark}\label{Rmk:Eigenvalues}
The fact that the PCG convergence, e.g. Figure~\ref{Fig:1d_pcg_cvg_eta0test}, for test case \eqref{Eq:Testproblem} exhibits an ``elbow'' behavior with two characteristic rates may partially be understood via the (apostori) computation of the generalized eigenvalues, shown in Figure~\ref{Fig:1d_gen_eigvalues}. Figure~\ref{Fig:1d_gen_eigvalues} provides numerical evidence that the generalized eigenvalues of the operators $\mG_{\eta,h}$ and $\mA$ admit a continuous spectrum confined to an interval $\sim [3 \times 10^{-2}, 3 \times 10^{-1}]$, with $\sim 6$ isolated eigenvalues outside of this interval.  An improved upper bound for the PCG error can then be used. After roughly 6 iterations (one for each isolated eigenvalue), the PCG error decays with a new effect upper bound \eqref{Eq:pcg_upperbound} with $\kapub \sim (3\times 10^{-1}) / (3 \times 10^{-2}) = 10$. This tighter upper bound more closely explains the convergence behavior in Figures~\ref{Fig:1d_pcg_cvg_eta0test}--\ref{Fig:1d_pcg_cvg_h0test} that occurs after the initial phase. \myremarkend  
\end{remark}

\begin{figure}
\begin{minipage}[t]{0.49\textwidth} 
\centering 
\includegraphics[width=.8\textwidth]{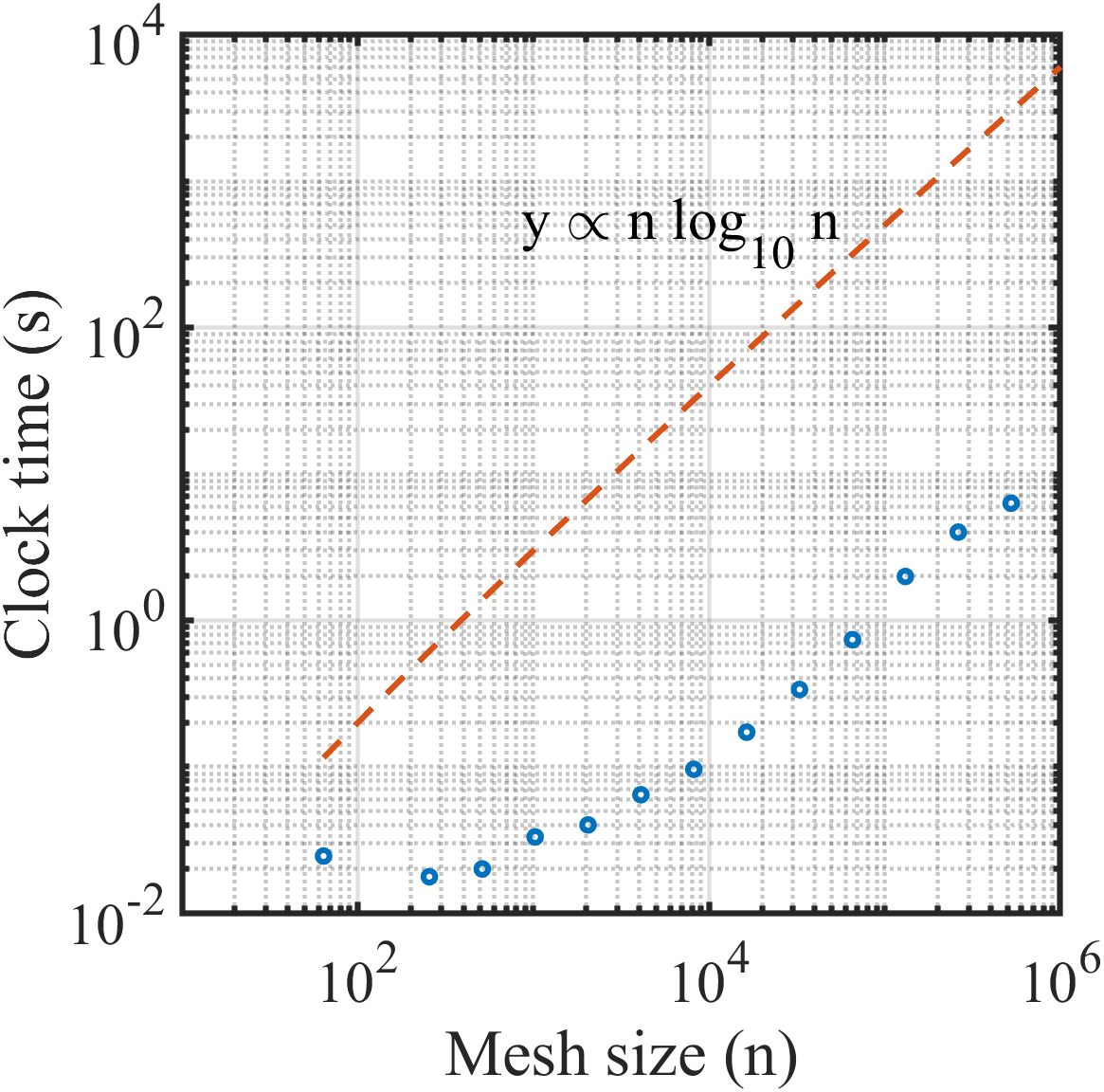}
\caption{PCG clock time ($\circ$; in seconds) required to solve for a residual of $\varepsilon <  10^{-12}$ versus mesh size $n$ for the 1d test problem \eqref{Eq:Testproblem1_var} with parameters $\eta_0 = h_0 = 1$. The reference line (dashed) presents the FLOP scaling for the FFT. The increase in computational time is due to the FFT scaling and not an increase in PCG iterations.} \label{Fig:1d_pcg_clocktime_test}
\end{minipage}
\hfill
\begin{minipage}[t]{0.49\textwidth} 
\centering 
    \includegraphics[width=.8\textwidth]{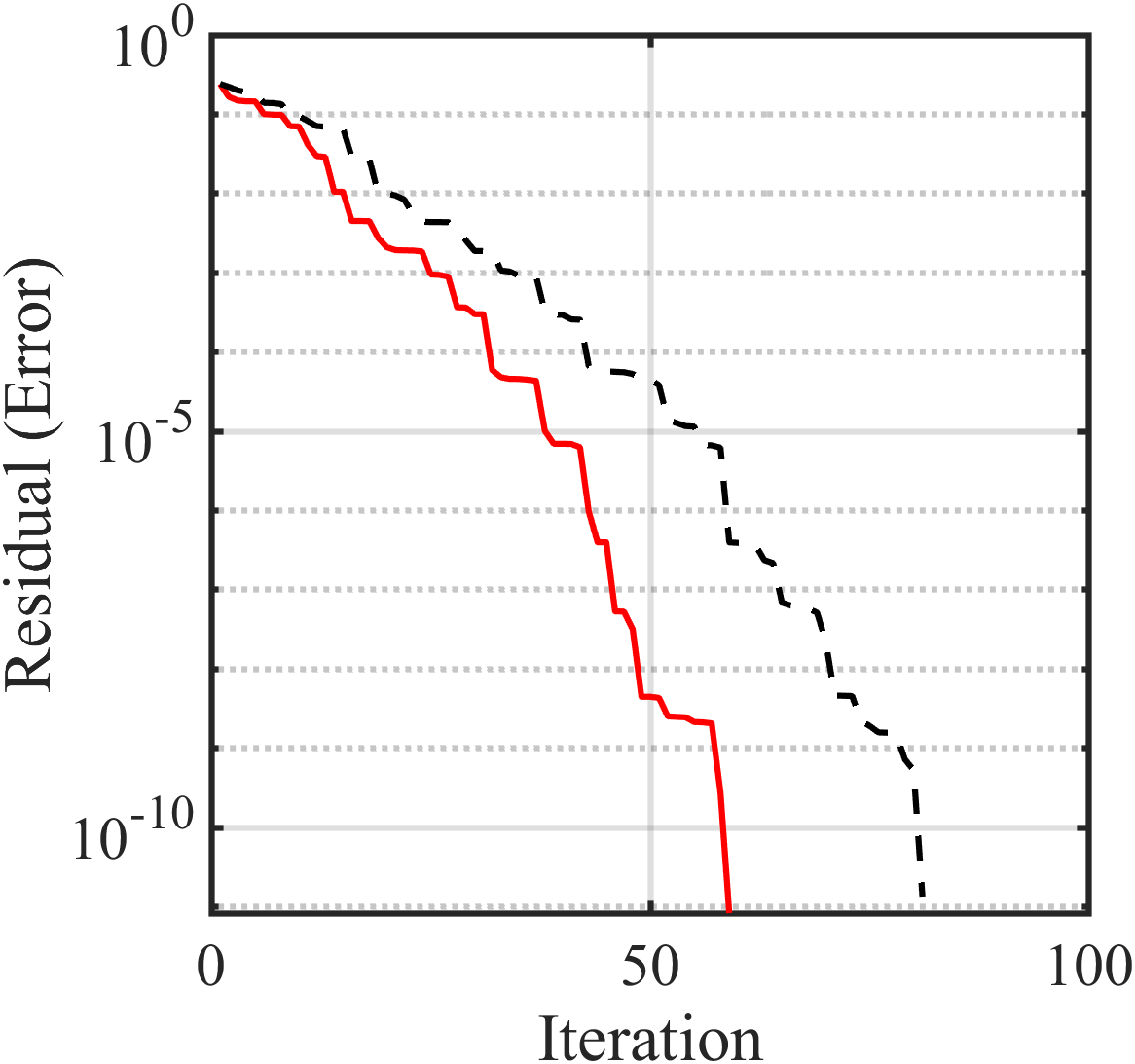}
    \caption{Comparison of the PCG residual for the coefficients \eqref{Eq:Optcoeff} (red; solid line) versus the simple (less optimal) coefficients in Remark~\ref{rmk:replacementcoeff} (black; dash line).  To achieve a tolerance $\varepsilon < 10^{-8}$, requires only 49 iterations versus 71. } \label{Fig:preconditioner_coeff_comp}
\end{minipage}
\end{figure}

Figure~\ref{Fig:1d_pcg_clocktime_test} examines the clock time required to solve \eqref{Eq:Testproblem1_var} with $\etavar = \hvar = 1$ at different mesh resolutions ($2^{5} \leq n \leq 2^{20}$).  The figure demonstrates that clock time scales well with the expected floating point operations of the FFT ($\propto n \log_{10} n$).  

Lastly, it is interesting to compare the optimal coefficient \eqref{Eq:FlatbathCoeff}--\eqref{Eq:Optcoeff} with the coefficients in Remark~\ref{rmk:replacementcoeff}, which are often easier to estimate for use in a large-scale computation.  As preconditioners for conjugate gradient, the two sets of coefficients perform very similarly unless 
\begin{align}\label{eq:replacementcoeff_cond}
    \max_{x \in \Omega} \Big(\eta(x) + \lambda_{+} \eta(x)|h_x(x)|^2 \Big) \ll \eta_{\rm max} (1 + \lambda_{+} |h_x|_{\rm max}^2)
\end{align}
i.e., the value of $\eta(x)$ is small in locations where $h_x$ is large.  In settings where \eqref{eq:replacementcoeff_cond} hold, the coefficient $\sigma^{\star\star}$ in Remark~\ref{rmk:replacementcoeff} will be larger than $\sigma^{\star}$ and may lead to a a drop in the preconditioner performance. The next test case specifically designs such an occurrence where \eqref{eq:replacementcoeff_cond} holds.

Set $\hvar = 1$ and $h(x)$ as in \eqref{Eq:Testproblem1_var}, and take a discontinuous water depth given by a square wave (the coefficients in the weak form of the operators $\mG_{\eta, h}$ and $\mA$, $\eta$ need not be continuous)
\begin{align*}
    \eta(x) = \left\{ \begin{array}{cc}
    0.1 & \textrm{if} \; x \in S \\
    0.1 \eta_0 & \textrm{if} \; x \notin S
    \end{array}\right.  \qquad \textrm{where} \qquad S = \{x \in \Omega \; | \; |h_x(x)|^2 \leq 0.2\} \,.
\end{align*}
Figure~\ref{Fig:preconditioner_coeff_comp} demonstrates the improved performance of the coefficients \eqref{Eq:Optcoeff} (red) over the coefficients in Remark~\ref{rmk:replacementcoeff} (dashed black).

\subsection{PCG Tests in 2d: Curvature of the Bathymetry}\label{subsec:2dpreconditioner}
We now examine two dimensional bathymetry effects on the PCG performance.  As a first test of our code, we verify every result in \S~\ref{subsec:1dpreconditioner} by setting $v = 0$ and extending the single variable $h$, $\eta$, $u$ to $d=2$ as constants in the $y-$direction. 

Next, to provide a (truly) two dimensional test case, we examine an elliptical bump bathymetry with varying aspect ratios 
\begin{align}\label{Eq:ellipticalbump}
	h(x,y) &= 
        \begin{cases}
            1-\tfrac{1}{2} \cos^2\left(\frac{\pi r}{2 r_0}\right),& |r|\leq r_0\\1,& |r|\geq r_0
        \end{cases}, \quad\quad 
        \textrm{where} \quad r=\sqrt{a^2 x^2+ b^2 y^2}, \quad r_0=\tfrac{1}{2} \,  \\
    \eta(x,y) &= h(x,y) + \textrm{exp}\big(\cos(2\pi x) \big) + \tfrac{1}{4}\sin(4\pi y) \, .
\end{align}
Note that both the magnitude and direction of $\nabla h$ vary throughout space to ensure the presence of anisotropic effects in the off-diagonal block entries of $\mG_{\eta,h}$.  The right hand side vector $\vec{b}$ is set to be the  discretization of the velocity $\vec{U} = \big(\cos(4\pi x), \cos( 4\pi y) \big)$. 

\begin{figure}[htp!]
    \includegraphics[width=.325\textwidth]{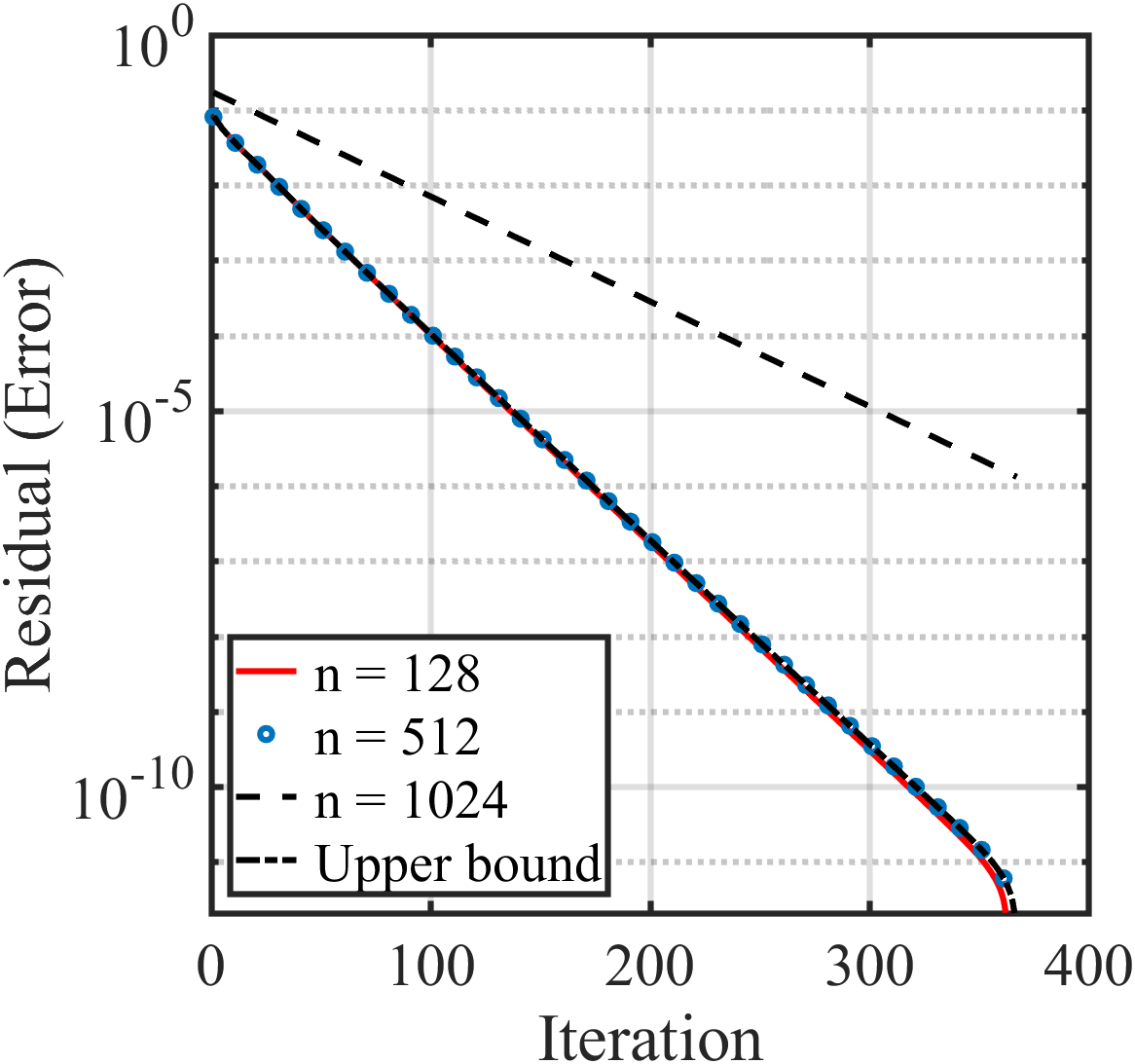}
    \includegraphics[width=.325\textwidth]{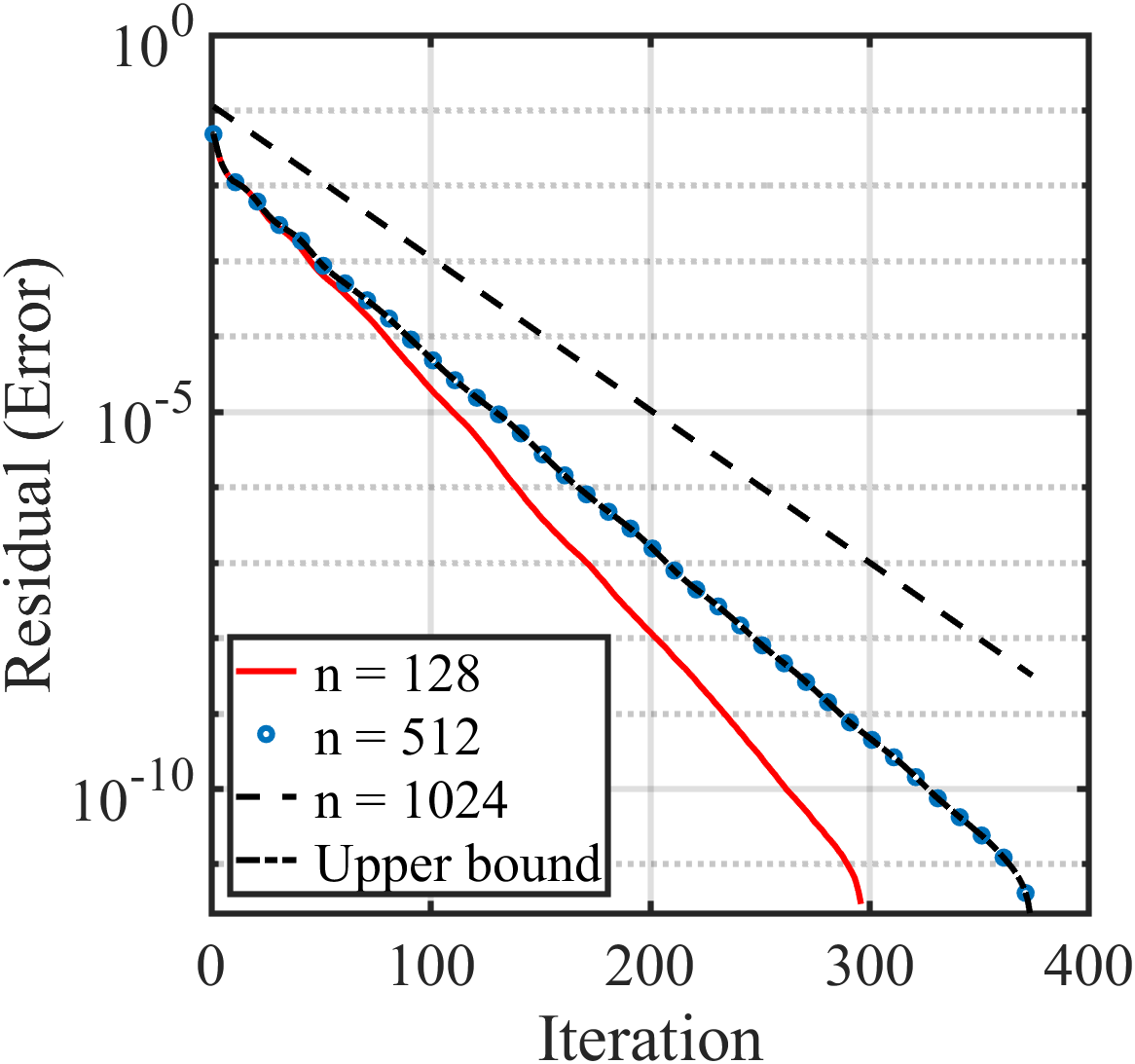}
    \includegraphics[width=.325\textwidth]{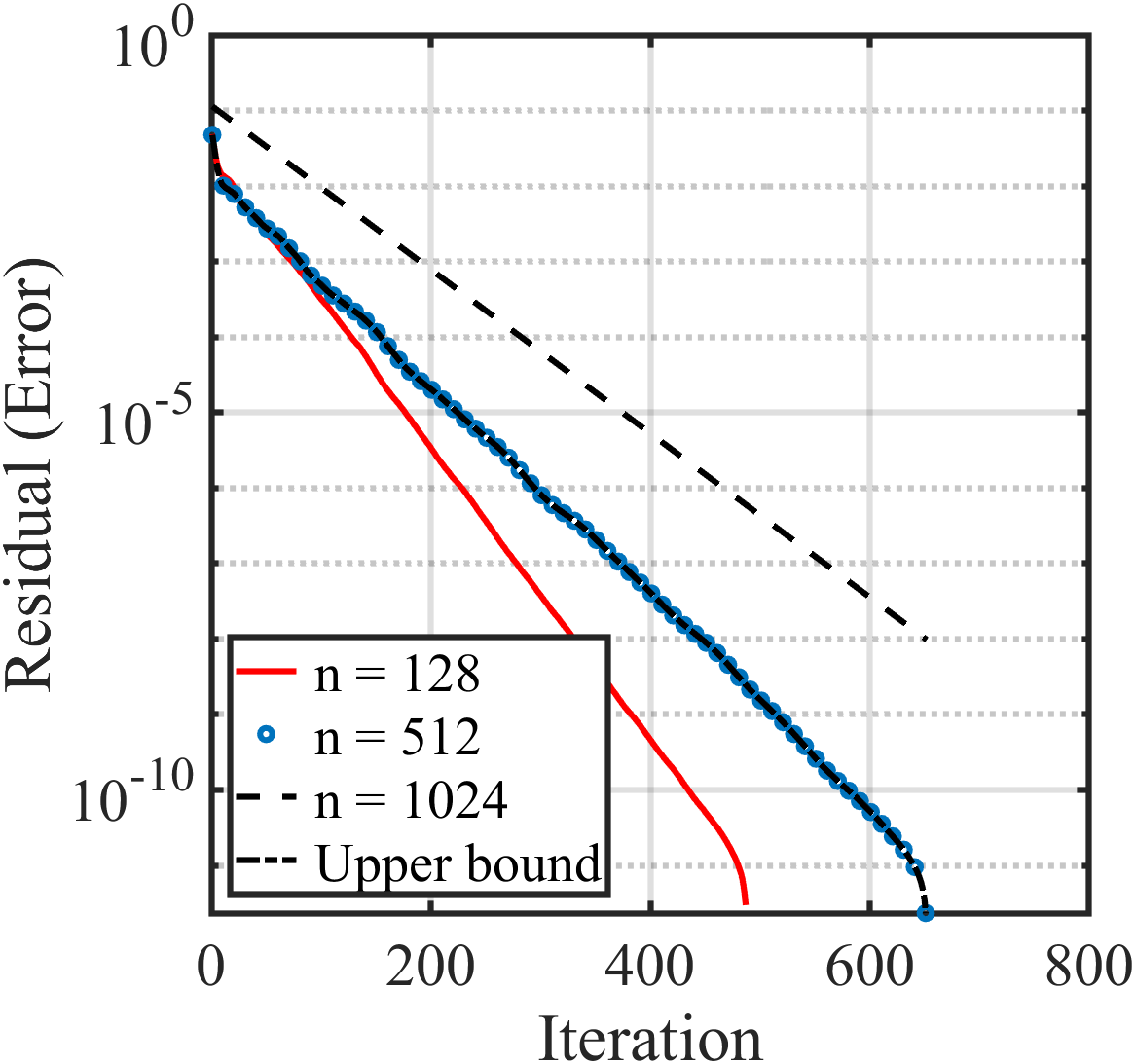}
    \caption{Test case for a 2d bathymetry \eqref{Eq:ellipticalbump}. Plots show the convergence of the PCG residual error $\varepsilon$ defined in \eqref{def:discrete_Gerr} with coefficients $(a,b) = (1,1)$ (left), $(1,20)$ (middle) and $(1,40)$ (right). The upper bound (dashed line) is \eqref{eq:CGRate} using the upper bound $\kapub$. Note the increase in iterations due to the effect of increasing aspect ratio and $\nabla h$.}\label{Fig:2d_pcg_cvg_bump}
\end{figure}

We test an increasing sequence of aspect ratios in the bathymetry that generate progressively larger gradients $\nabla h$.  Larger $\nabla h$ values increase the upper bound $\kapub$, which can be viewed as a proxy for the conditioning of $\mat{A}^{-1} \mat{G}$.  Figure~\ref{Fig:2d_pcg_cvg_bump} plots the PCG convergence versus iteration for three elliptical aspect ratios.  Again, the performance decreases with increasing $\nabla h$, however the preconditioner remains effective independent of the mesh size.

\subsection{Full Solution of the SGN Equations in 1d}\label{subsec:dswe_1d}
Building on the studies of the operator $\mA$ and $\mG_{\eta,h}$ from the previous section, we provide here, numerical test cases and simulations of the full time-dependent SGN equations.  We apply two different time integration strategies to solve the constraint form of the SGN equations as discussed in \S~\ref{Sec:TimeStepping}.

We first validate our codes in a flat bathymetry setting by solving the initial value problem with a solitary wave \eqref{ExSleta}--\eqref{ExSlu} (taking the domain large enough so that the exponential decay of the $\sech^2$ profile reaches machine precision).

\paragraph{Test 1: Manufactured solution}\label{paragraph:1d_manu_sol}
To validate the time-stepping approaches and perform comparative studies in the presence of a variable bathymetry, we test the numerical approaches using the method of \emph{manufactured solutions}. That is, the exact solution (on $L = 1$) is fixed to be
\begin{align}
    h(x) = 2 + \sin(2\pi x) \, , \qquad\quad 
    \eta^{\star}(x, t) &= 2 + \sin(2\pi x) \sin(10t) \, ,  \qquad\quad
    u^{\star}(x, t) = \cos(2\pi x) \cos(10 t) \, .
\end{align}
The SGN equations are then solved with an external forcing term chosen so that $\eta^{\star}$, $u^{\star}$ solve the PDE. 

The solution, initialized to $\eta^{\star}$, $u^{\star}$, is evolved to a final time $t_f = 1$. Since, $\eta^{\star}$, $u^{\star}$ are periodic-in-time, the external forcing is also periodic-in-time and induces a Floquet instability in the equations, which we confirmed by computing the spectral radius of the Floquet monodromy matrix. Evolving the SGN equations for a longer time with this external forcing allows the instability to grow in time.  The final time $t_f = 1$ in our tests is specifically chosen small enough to test numerical convergence and computational efficiency via clock time while avoiding unnecessary complications from instabilities that could manifest in the equations over longer times. 

\begin{figure}[htp!]
    \includegraphics[width=.325\textwidth]{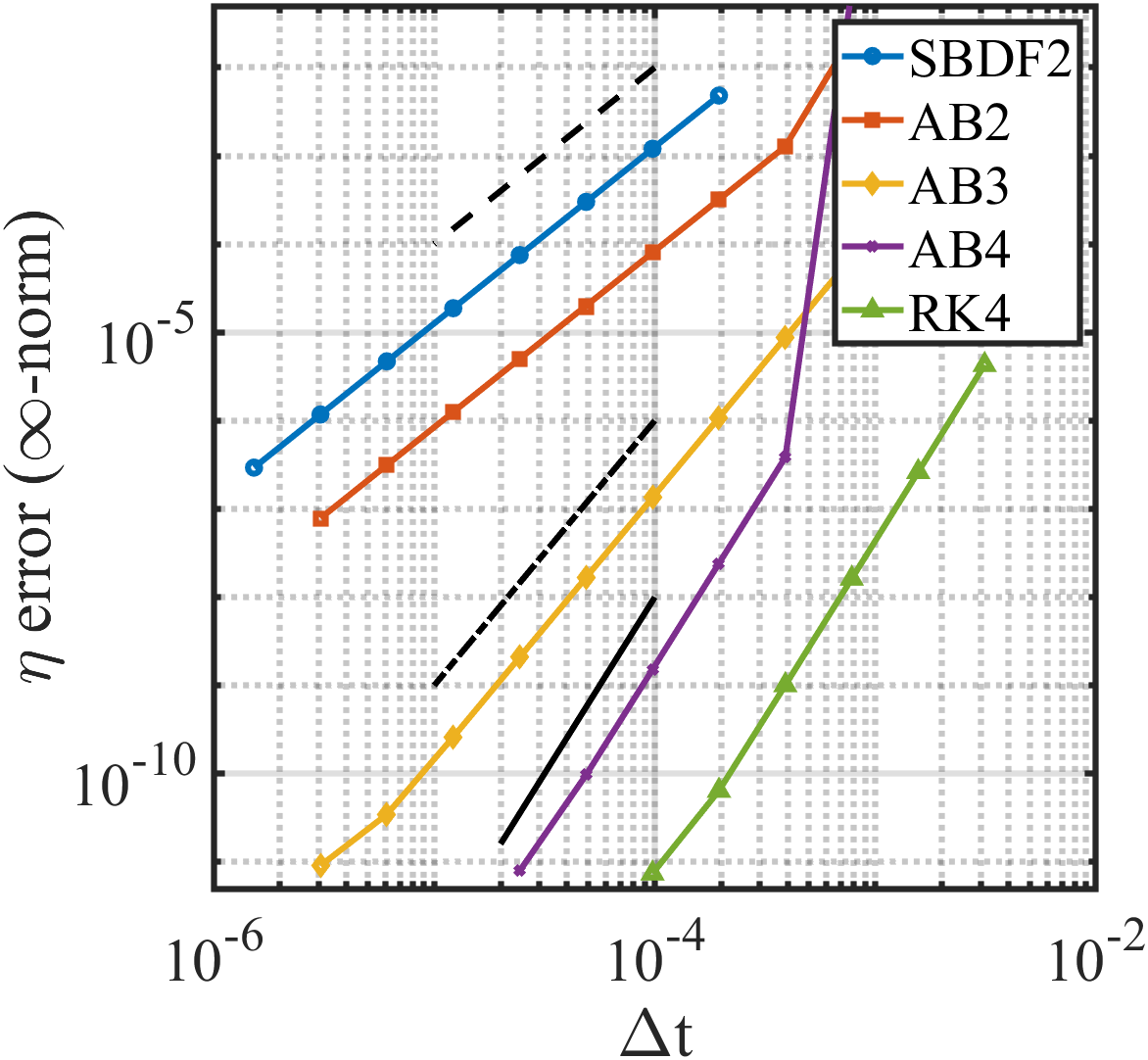}
    \includegraphics[width=.325\textwidth]{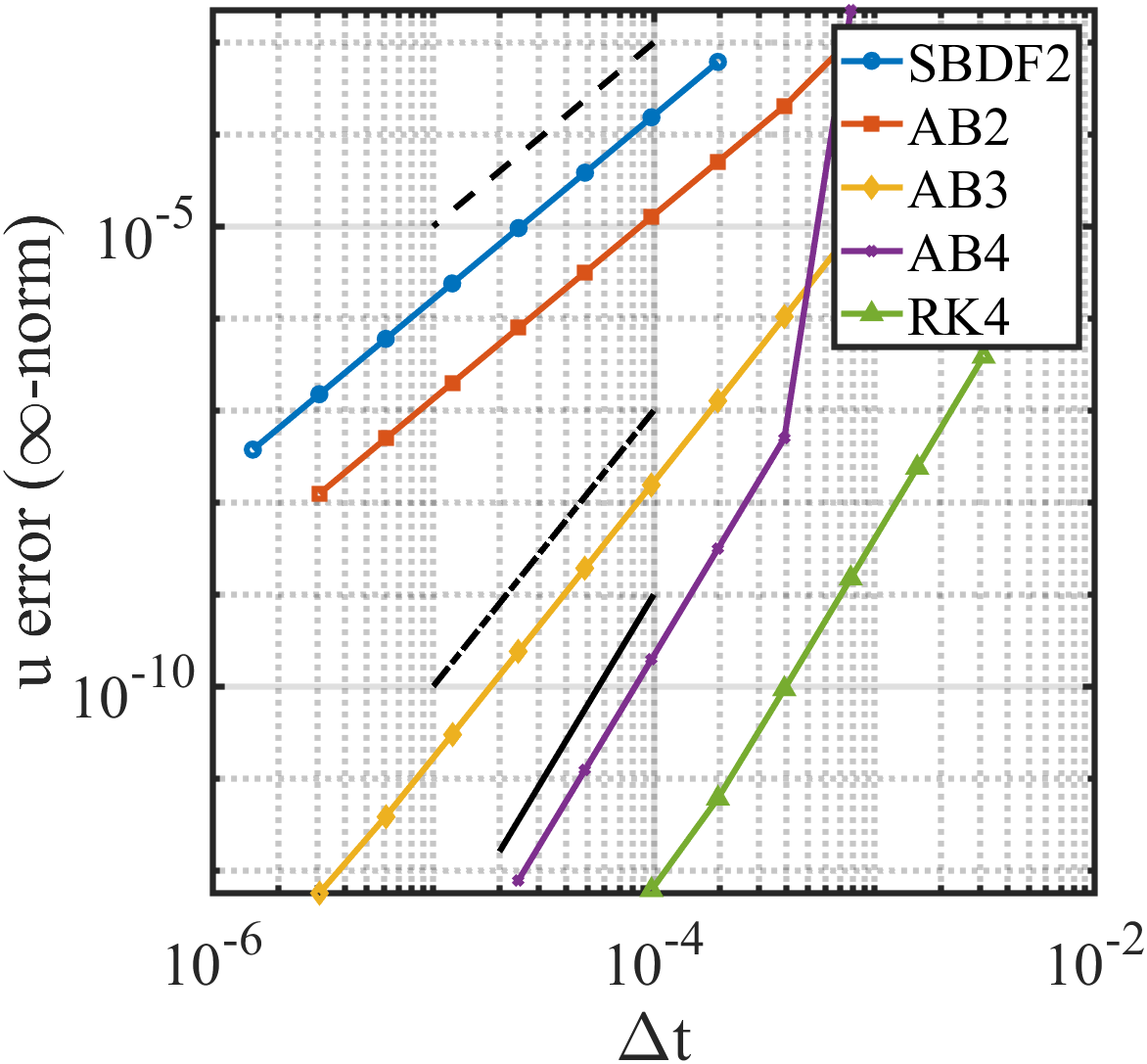}
    \includegraphics[width=.325\textwidth]{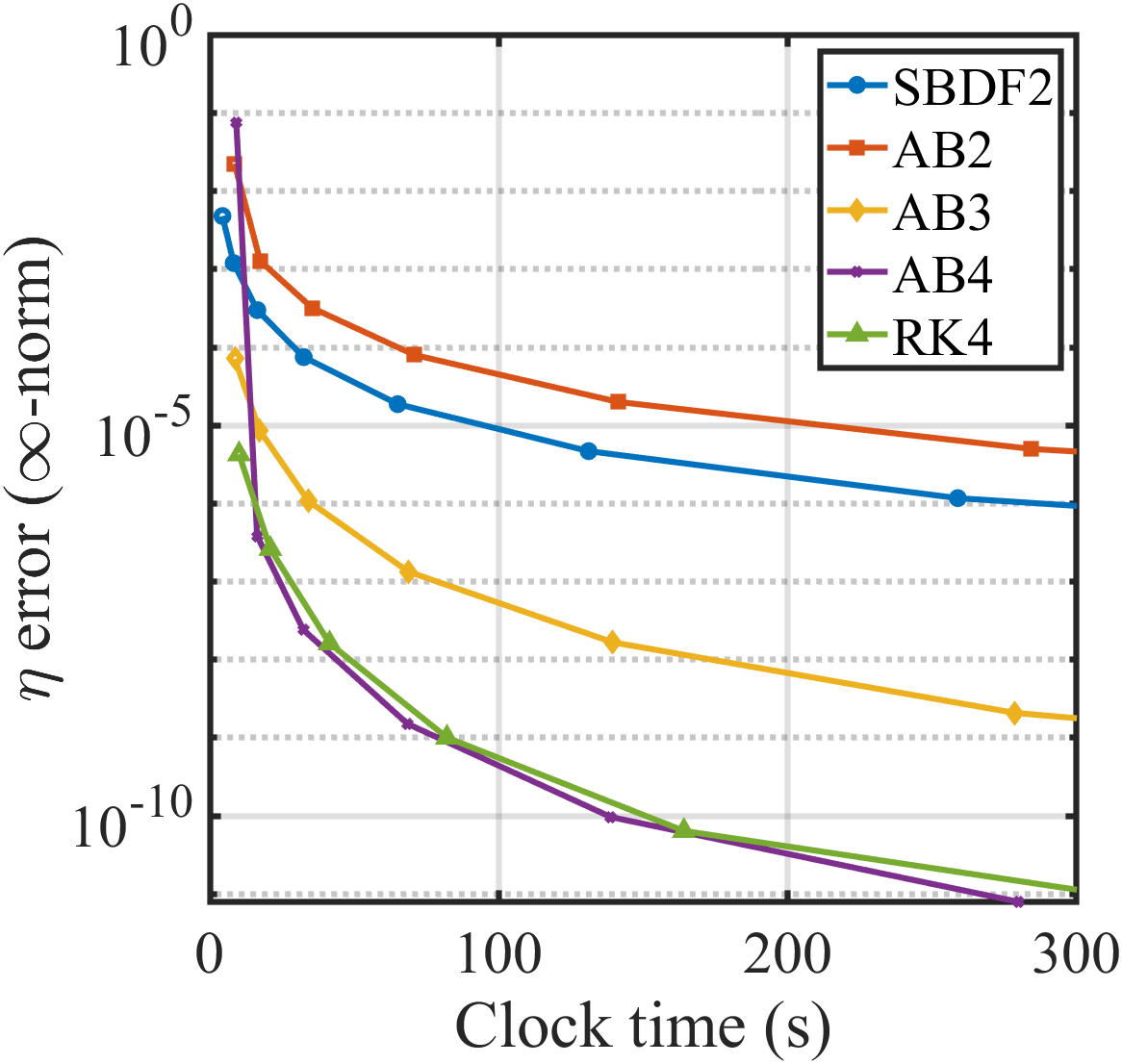}
    \caption{Figures provide an $L^{\infty}$ convergence study for $\eta$ (left), and $u$ (middle) as well as a clock time comparison (right) for the \S~\ref{paragraph:1d_manu_sol} Test Case 1. The figures compare the linearly implicit SBDF2 approach versus explicit schemes with PCG (AB2, AB3, AB4 and RK4). All methods converge at the expected rates, and higher order methods provide greater accuracy for a given clock time; the linearly implicit SBDF2 outperforms AB2 due to the single linear solve involving $\mA$.}\label{Fig:1d_cvg_man_sol}
\end{figure}

All numerical computations are run using a fixed mesh with $n = 256$ grid points, which is small enough to ensure spatial discretization errors remain well below temporal errors, i.e., are close to machine precision. For explicit time-stepping, the PCG tolerance is set to $10^{-14}$.  Figure~\ref{Fig:1d_cvg_man_sol} provides a convergence plot in $L^{\infty}$ for $\eta$ (left), $u$ (middle) demonstrating that the methods converge and do so at the expected rates. The CFL number which dictates the largest values of $\dt = C \Delta x$ ($\Delta x = 1/n$) are method dependent, and range from $C = .2$ for RK4 to $C = .05$ for SBDF2. Figure~\ref{Fig:1d_cvg_man_sol} (right) also provides a clock time comparison. As expected, the high order methods outperform lower order methods, while the linearly implicit SBDF2 outperforms AB2 with PCG.

\paragraph{Test 2: Simulation of a solitary wave approaching an underwater shelf}\label{paragraph:1d_shelf}

\begin{figure}
\begin{minipage}[t]{0.49\textwidth} 
\centering 
\includegraphics[width=.9\textwidth]{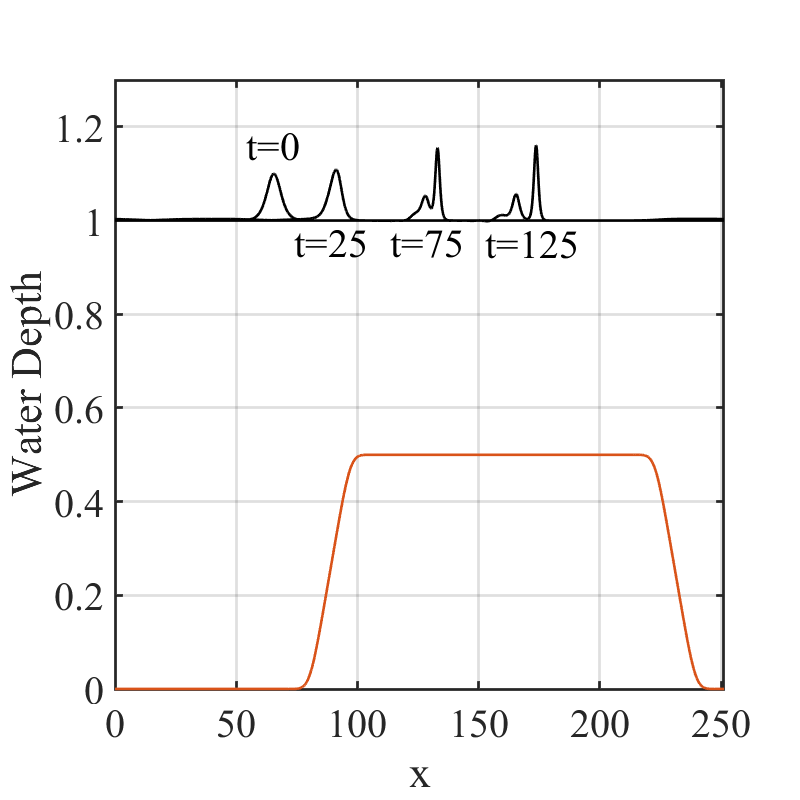}
\end{minipage}
\hfill
\begin{minipage}[t]{0.49\textwidth} 
\centering 
    \includegraphics[width=.9\textwidth]{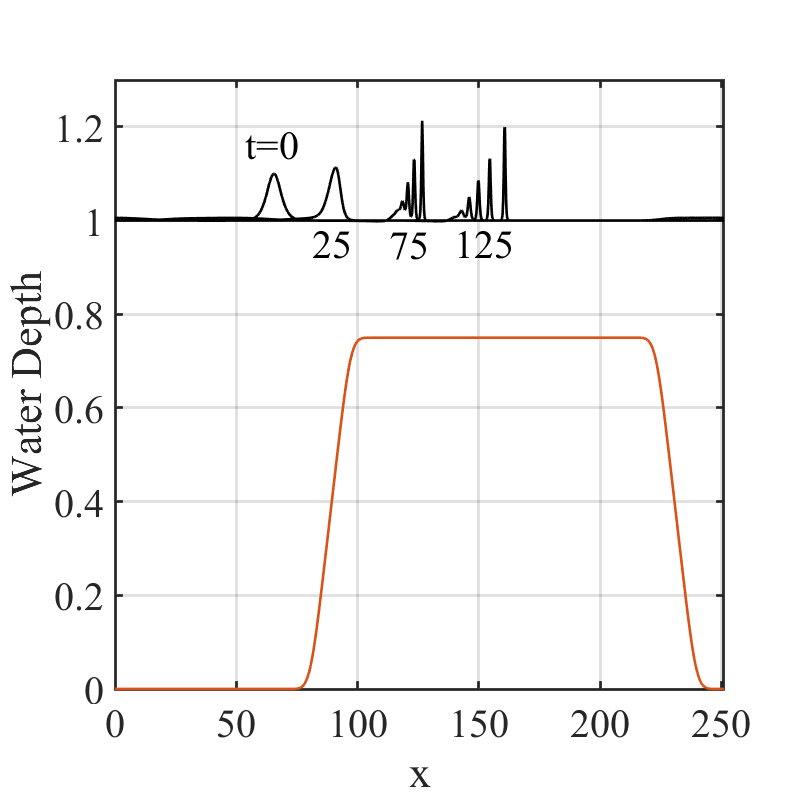}
\end{minipage}
\caption{Solitary wave traveling to the right over an underwater shelf (orange; solid line) with height $0.5$ (left) and $0.75$ (right). Figures~\ref{Fig:SolitonOnBeachZeta1} and \ref{Fig:SolitonOnBeachZeta2} show enlarged images of the waves. }\label{Fig:WavesOverShelf}
\end{figure}

\begin{figure}[htp!]
    \includegraphics[width=.5\textwidth]{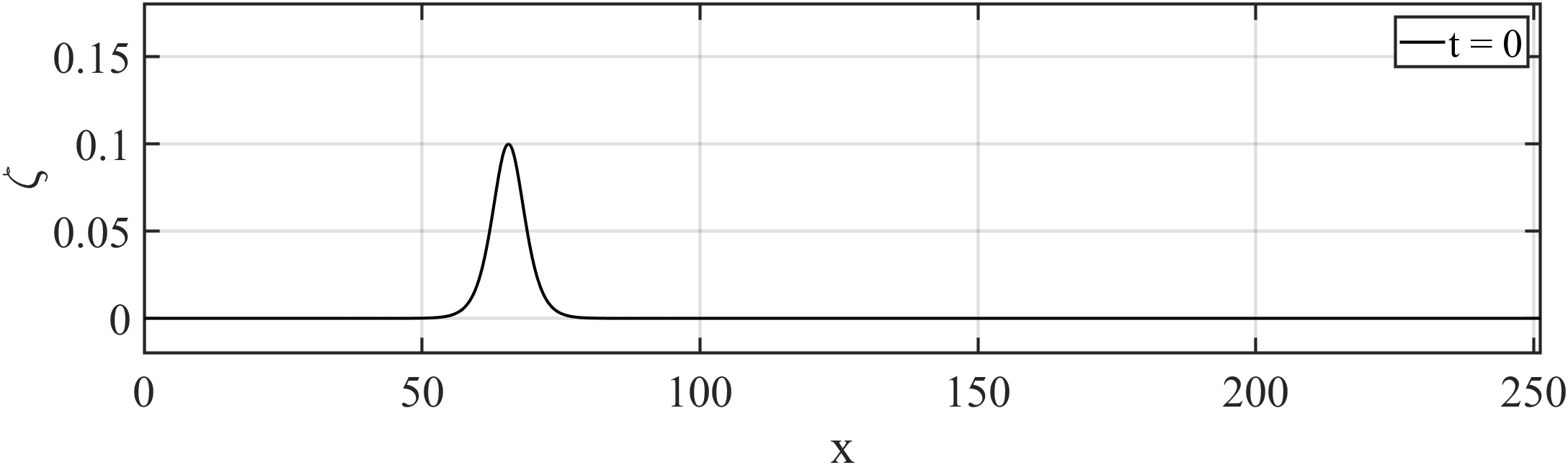}
    \includegraphics[width=.5\textwidth]{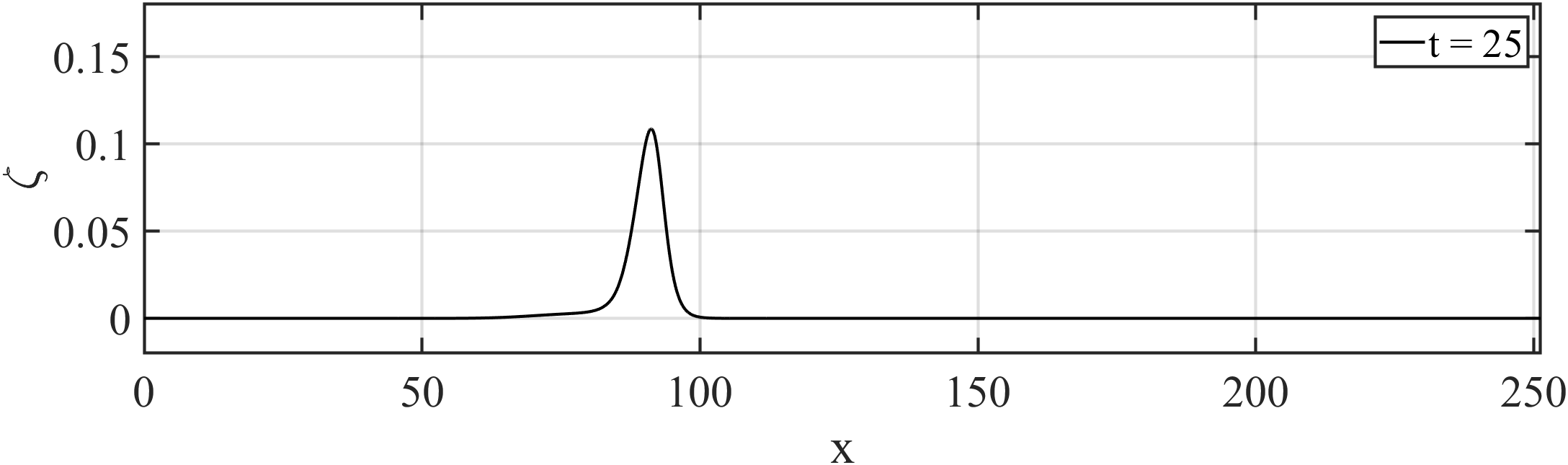}\\
    \includegraphics[width=.5\textwidth]{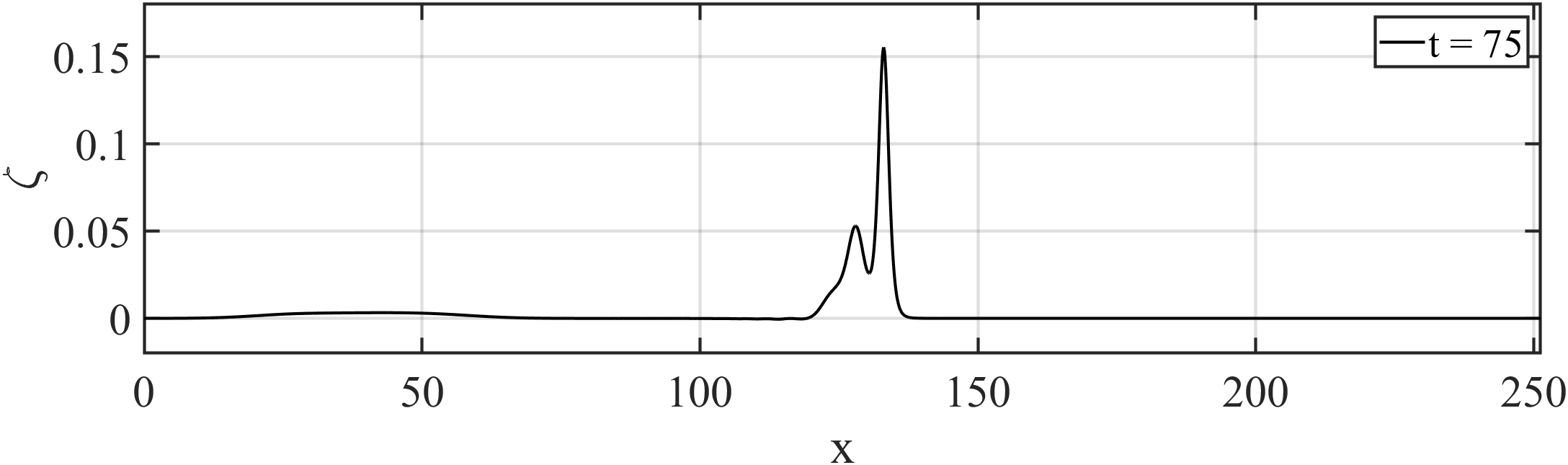}
    \includegraphics[width=.5\textwidth]{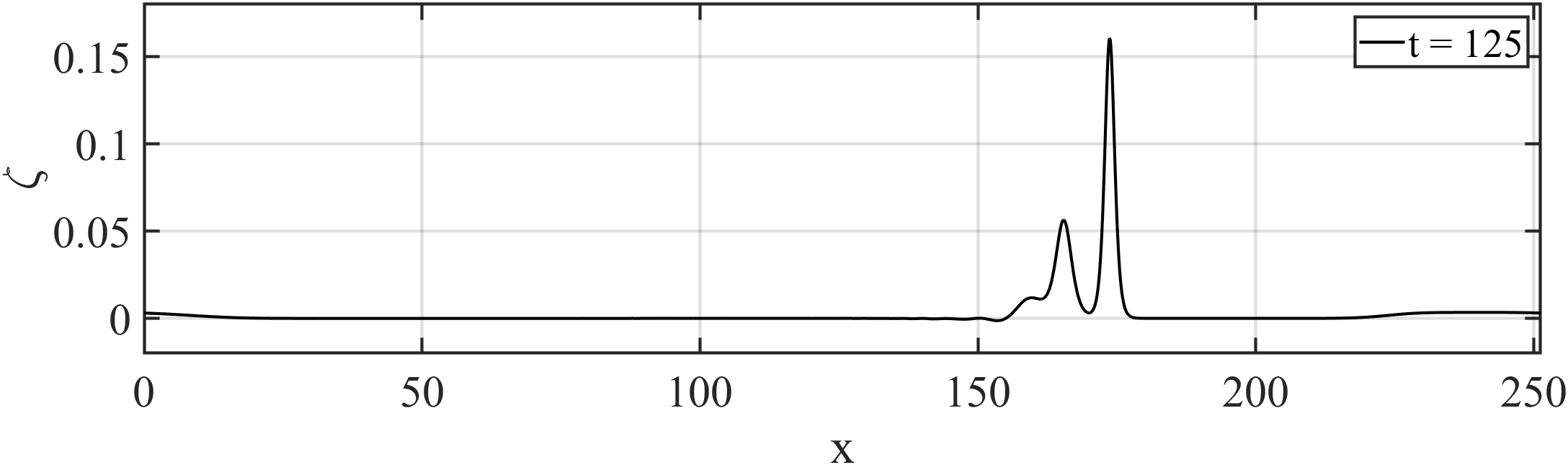}
    \caption{Enlargements of the solutions plotted in Figure~\ref{Fig:WavesOverShelf} (left). Note the reflected wave.}\label{Fig:SolitonOnBeachZeta1}
\end{figure}

\begin{figure}[htp!]
    \includegraphics[width=.5\textwidth]{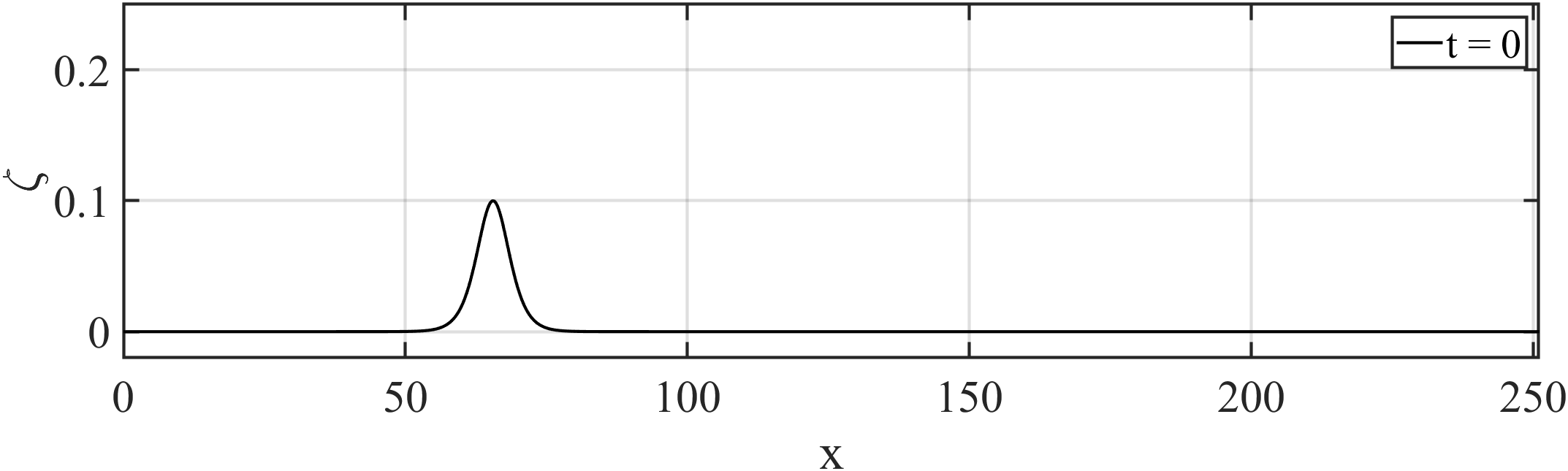}
    \includegraphics[width=.5\textwidth]{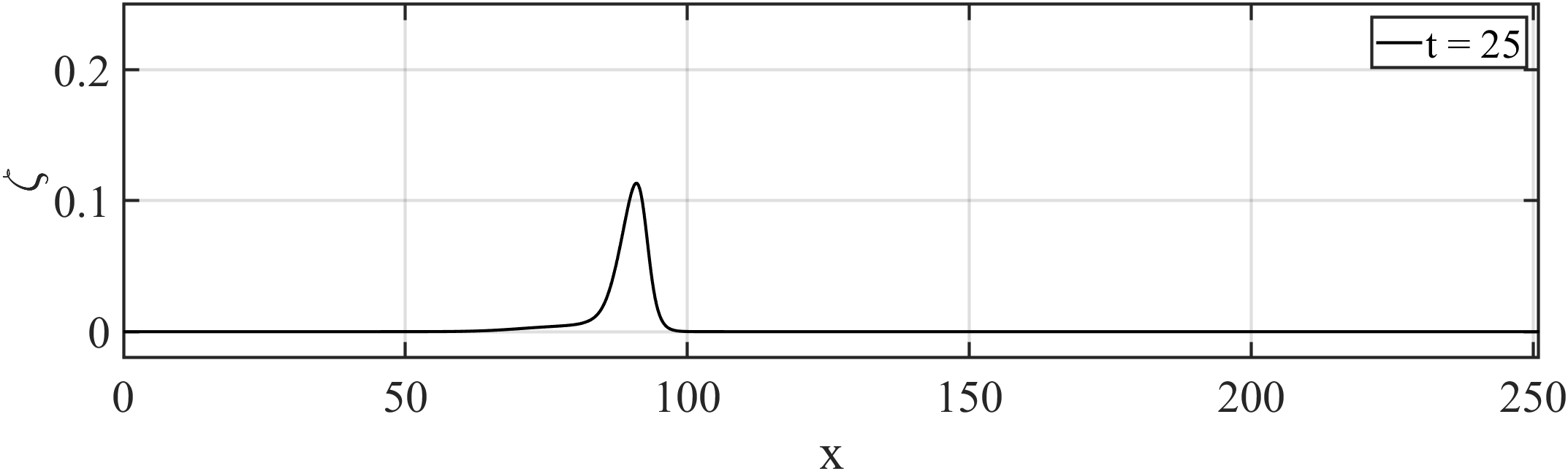}\\
    \includegraphics[width=.5\textwidth]{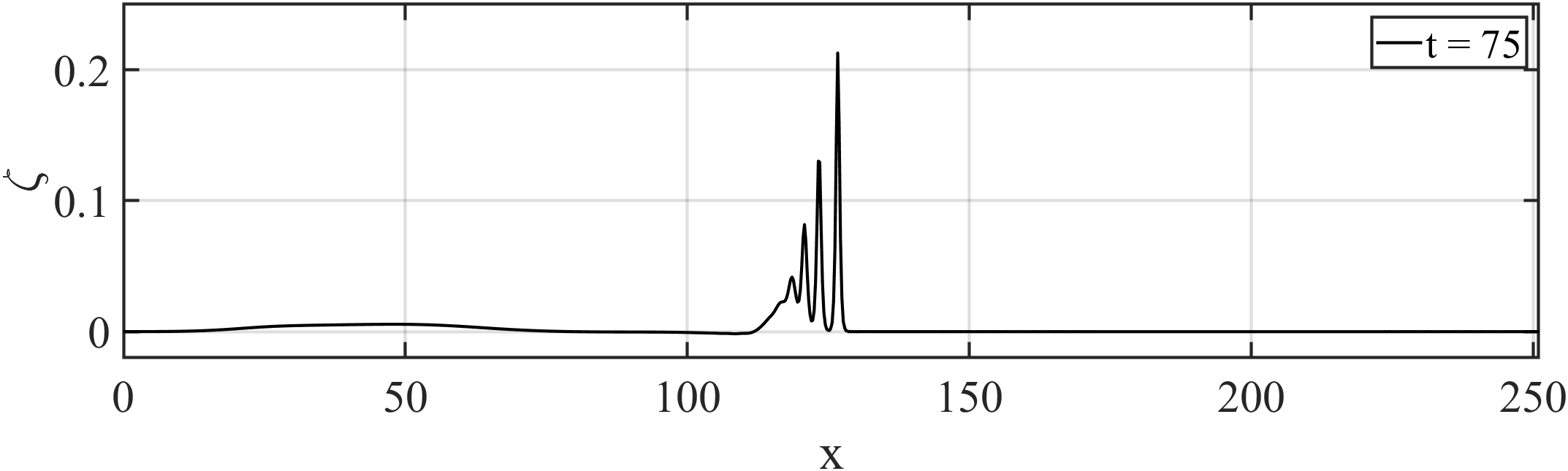}
    \includegraphics[width=.5\textwidth]{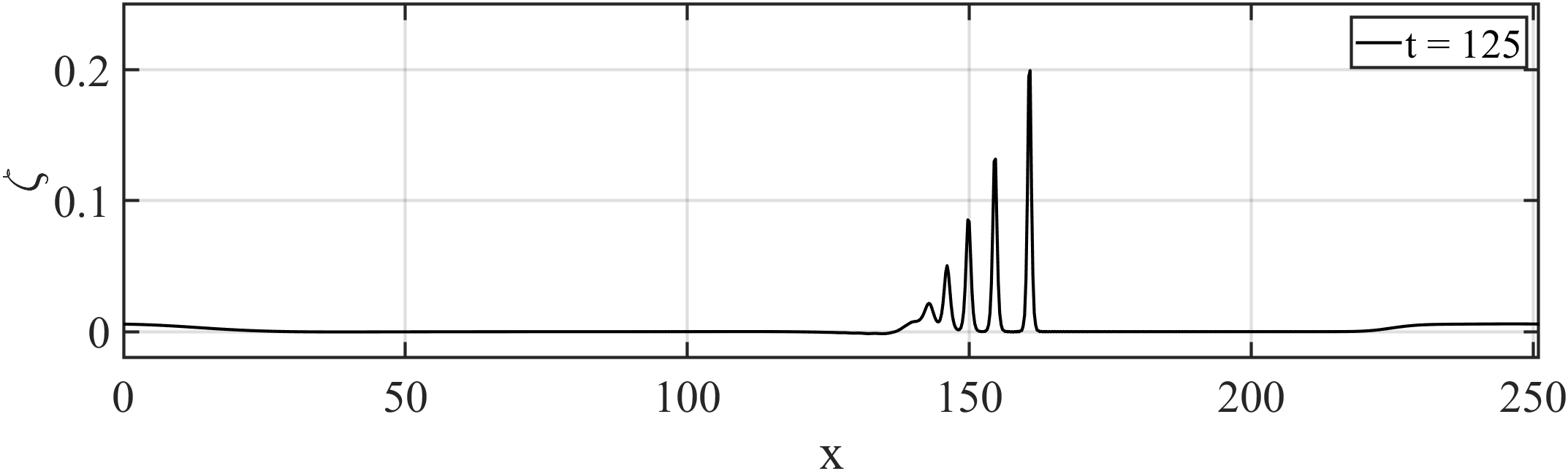}
    \caption{Enlargements of the solutions plotted in Figure~\ref{Fig:WavesOverShelf} (right). Note the reflected wave.}\label{Fig:SolitonOnBeachZeta2}
\end{figure}

In this test, we simulate a solitary wave traveling over an underwater shelf (here $L = 80\pi$). The initial data is set to the traveling solitary wave solution \eqref{ExSleta}--\eqref{ExSlu} with the water depth $h_0$ set to the height preceding the shelf.  The exact form of the shelf is a piecewise linear square wave convolved with a smooth bump function.  Figure~\ref{Fig:WavesOverShelf} visualizes the effect of dispersion for two shelf heights. Figures~\ref{Fig:SolitonOnBeachZeta1}--\ref{Fig:SolitonOnBeachZeta2} provide enlarged views of the waves for times $t = 0, 25, 75, 125$, which also shows the presence of a small reflected wave. 

Aliasing can present a source of error when spectral methods are used to solve nonlinear problems.  In all of our simulations, we have not made use of dealiasing. We found that regardless of the choice of time integration, solitary wave solutions with larger amplitudes can destabilize over long times scales.  We also found that dealiasing through the introduction of padding the high frequency wave numbers with additional zeros, can aid in avoiding the instability.

\subsection{Full Solution of the SGN equations in 2d}\label{subsec:dswe_2d}
This test problem solves the SGN equations on the doubly periodic square with length $L = 1$. To test the codes, we first performed a convergence study on the time integration schemes in two dimensions using a solitary wave solution. The results were analogous to the one dimensional setting. 

For the two dimensional test, we set a variable bathymetry to be $h(x,y) = h_0 - b(x,y)$ where $h_0 = 1.5a$ and $b(x,y)$ is the periodic extension of   
\begin{align}\label{Eq:2dBarrier}
    b(x,y) = 0.75 \,a \, \textrm{exp}( -r(x,y)^2 / (0.2)^2) 
    \qquad 
    r(x,y) = \sqrt{(x - 0.5)^2 + (y-0.5)^2} \, .
\end{align}
The surface $\eta$ and velocity $\vec{u}$ are initialized to be a traveling solitary wave front with an angle of $45$ degrees and $a = 1/100$.  If the bathymetry were constant, the solution would remain a solitary wave, however the profile loses its shape upon traveling over the variable $h(x,y)$.  

\begin{figure}[htp!]
    \includegraphics[width=.325\textwidth]{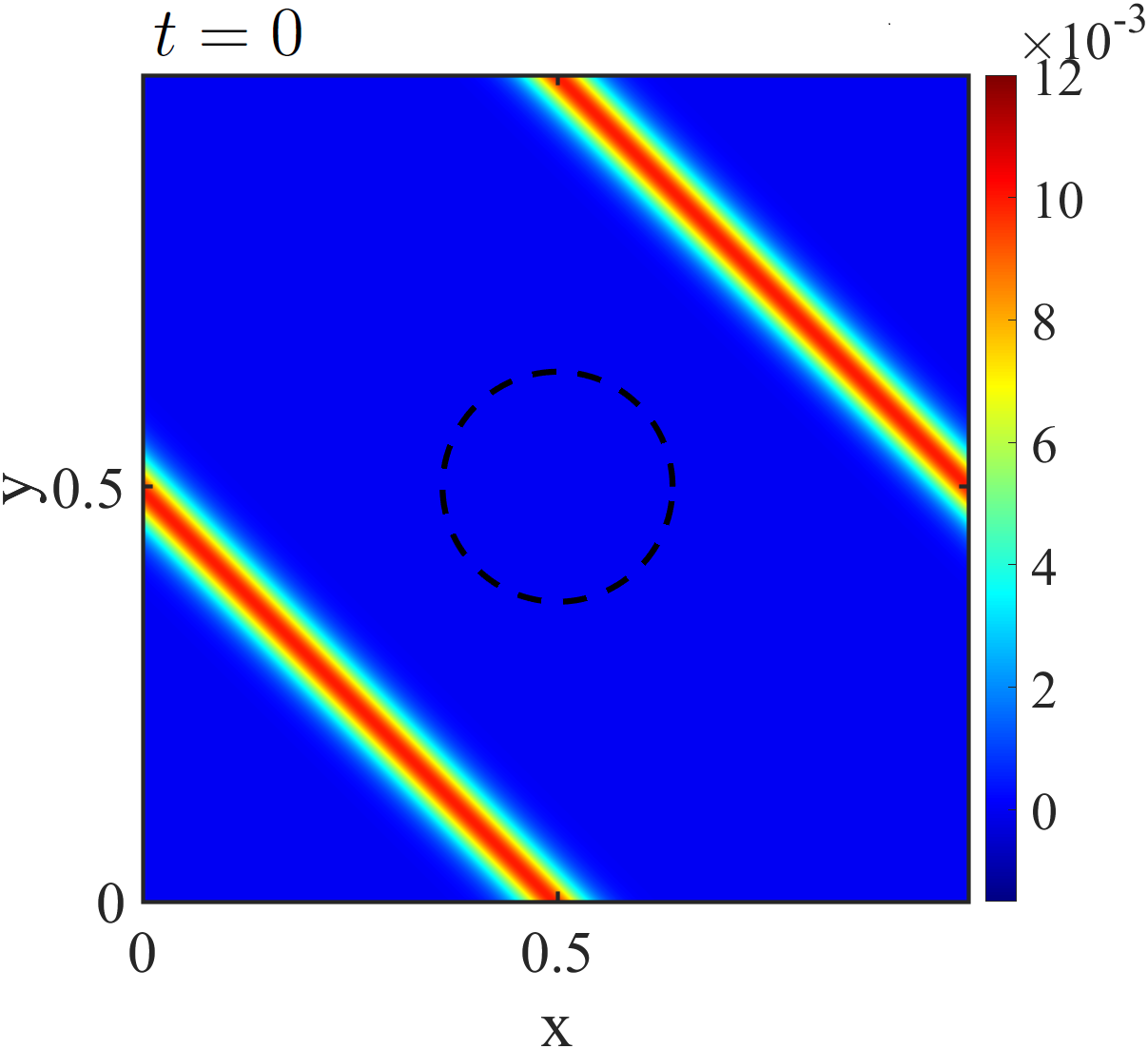}
    \includegraphics[width=.325\textwidth]{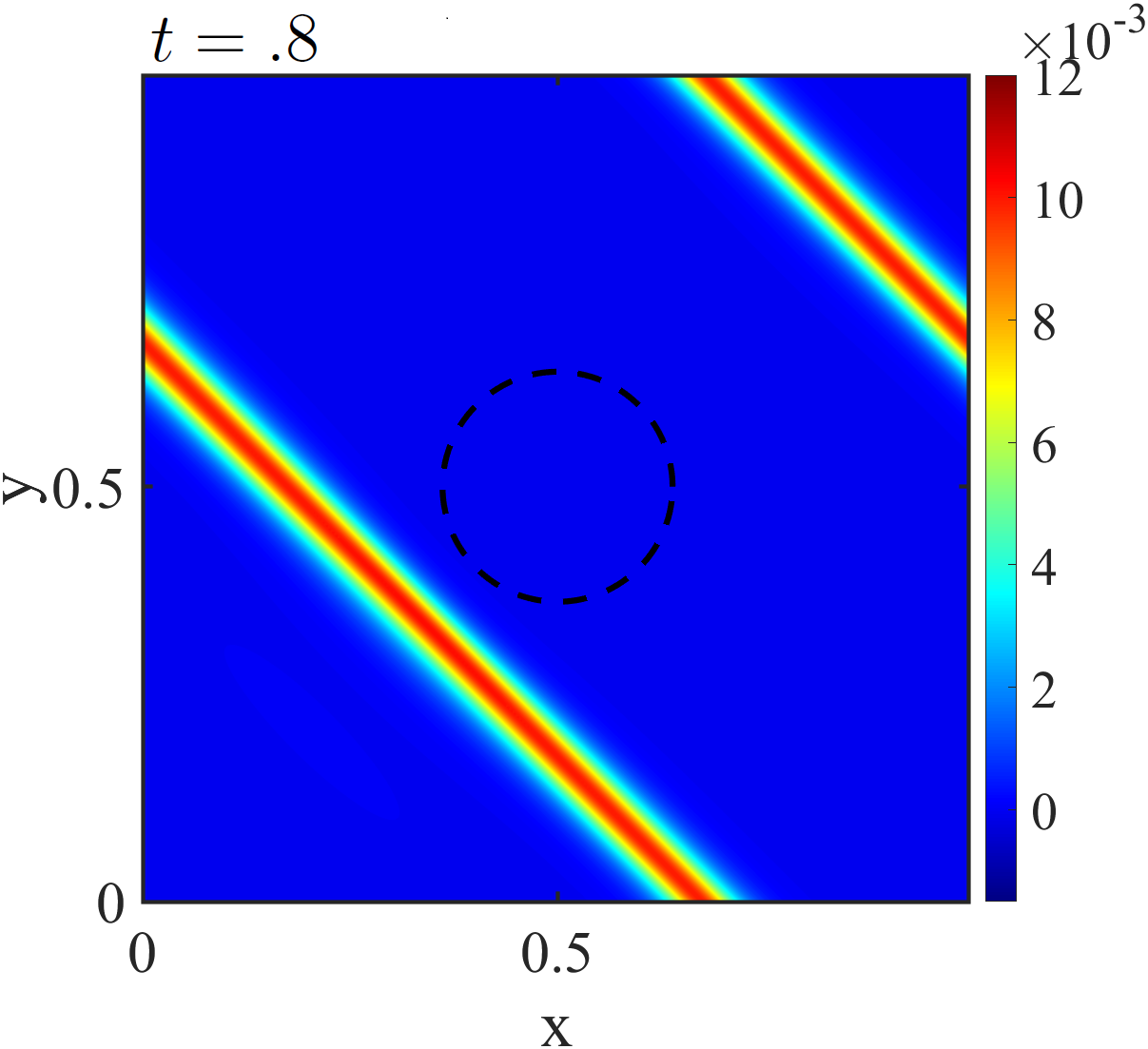}
    \includegraphics[width=.325\textwidth]{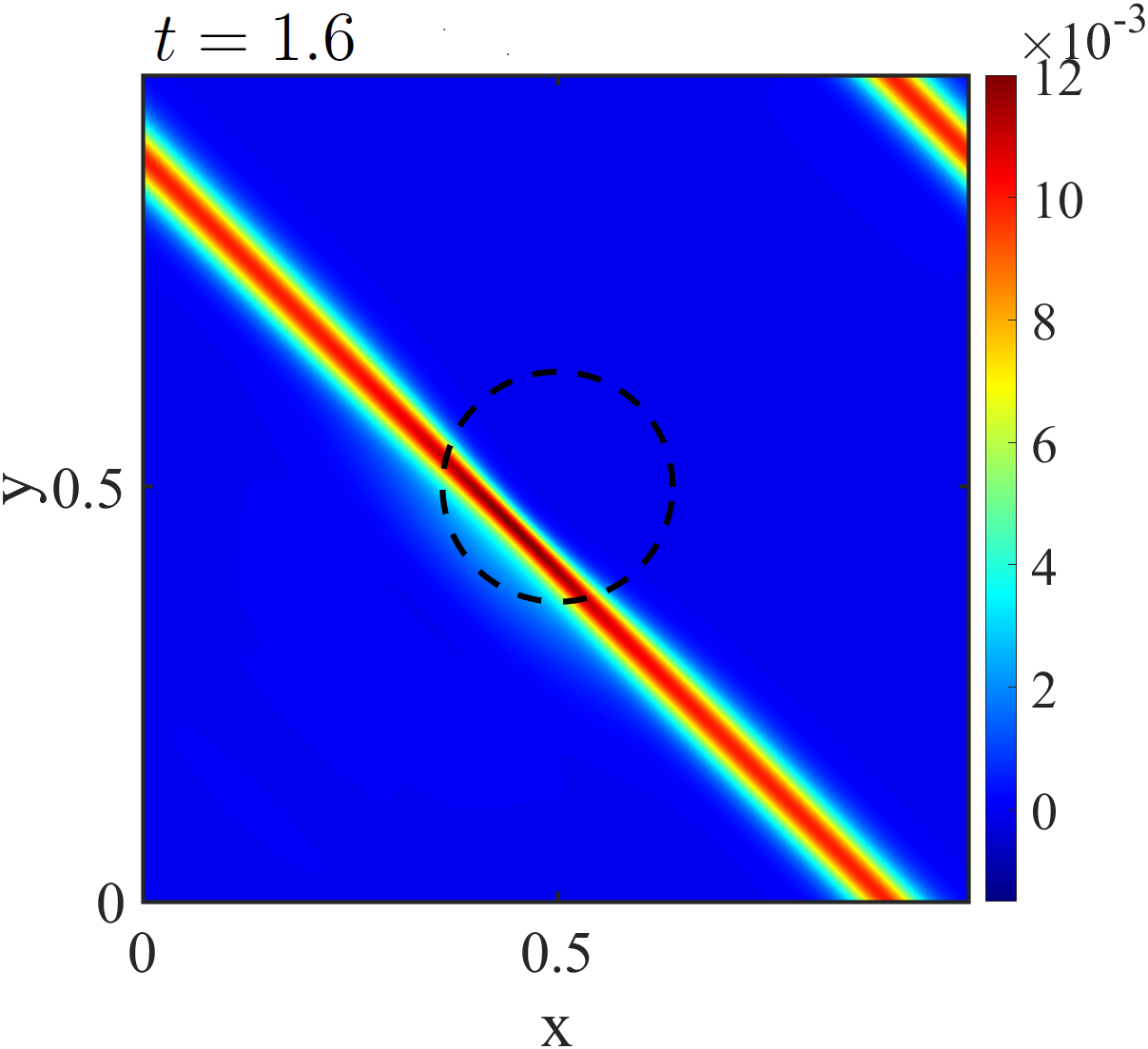}\\
    \includegraphics[width=.325\textwidth]{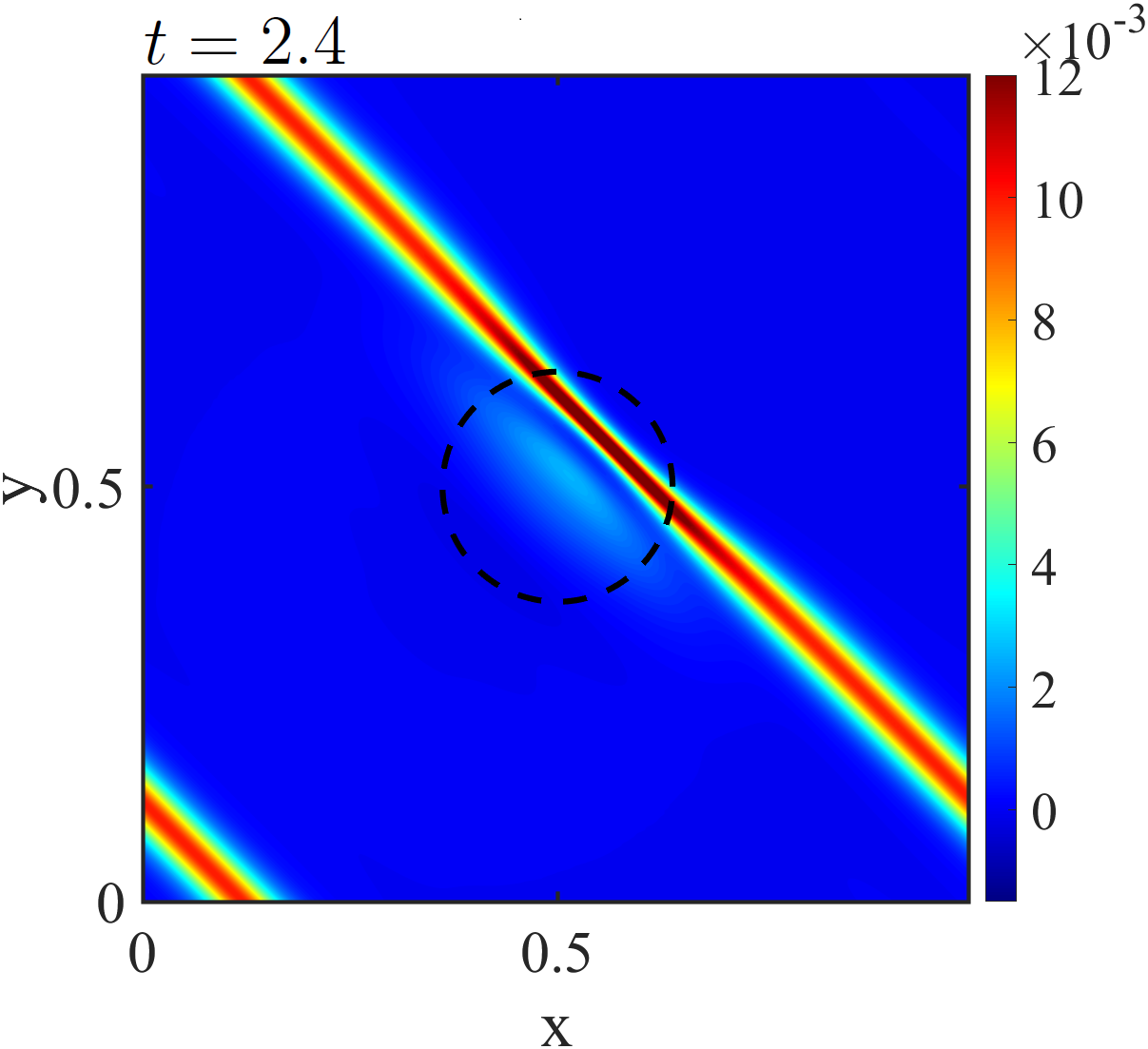}
    \includegraphics[width=.325\textwidth]{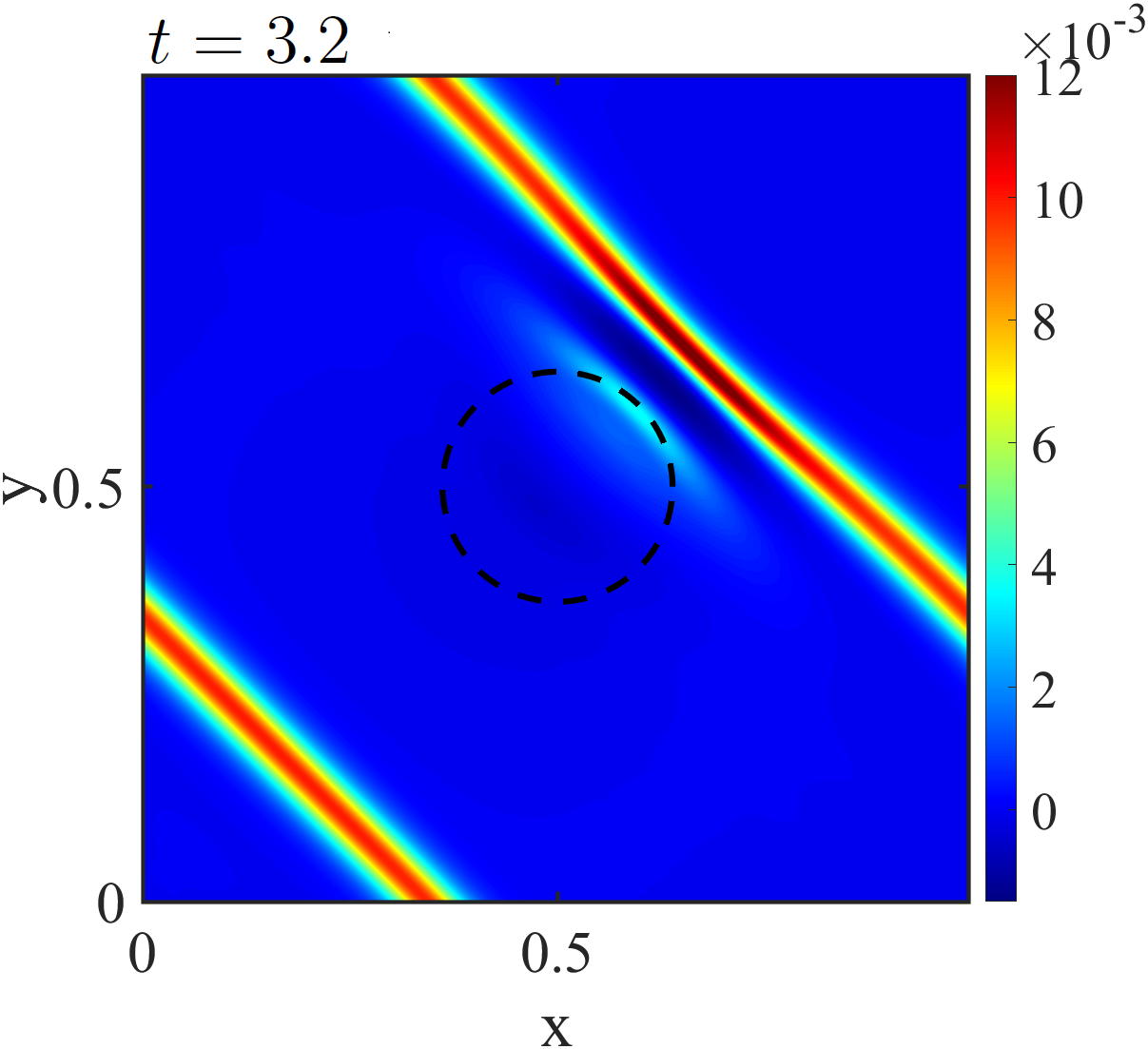}
    \includegraphics[width=.325\textwidth]{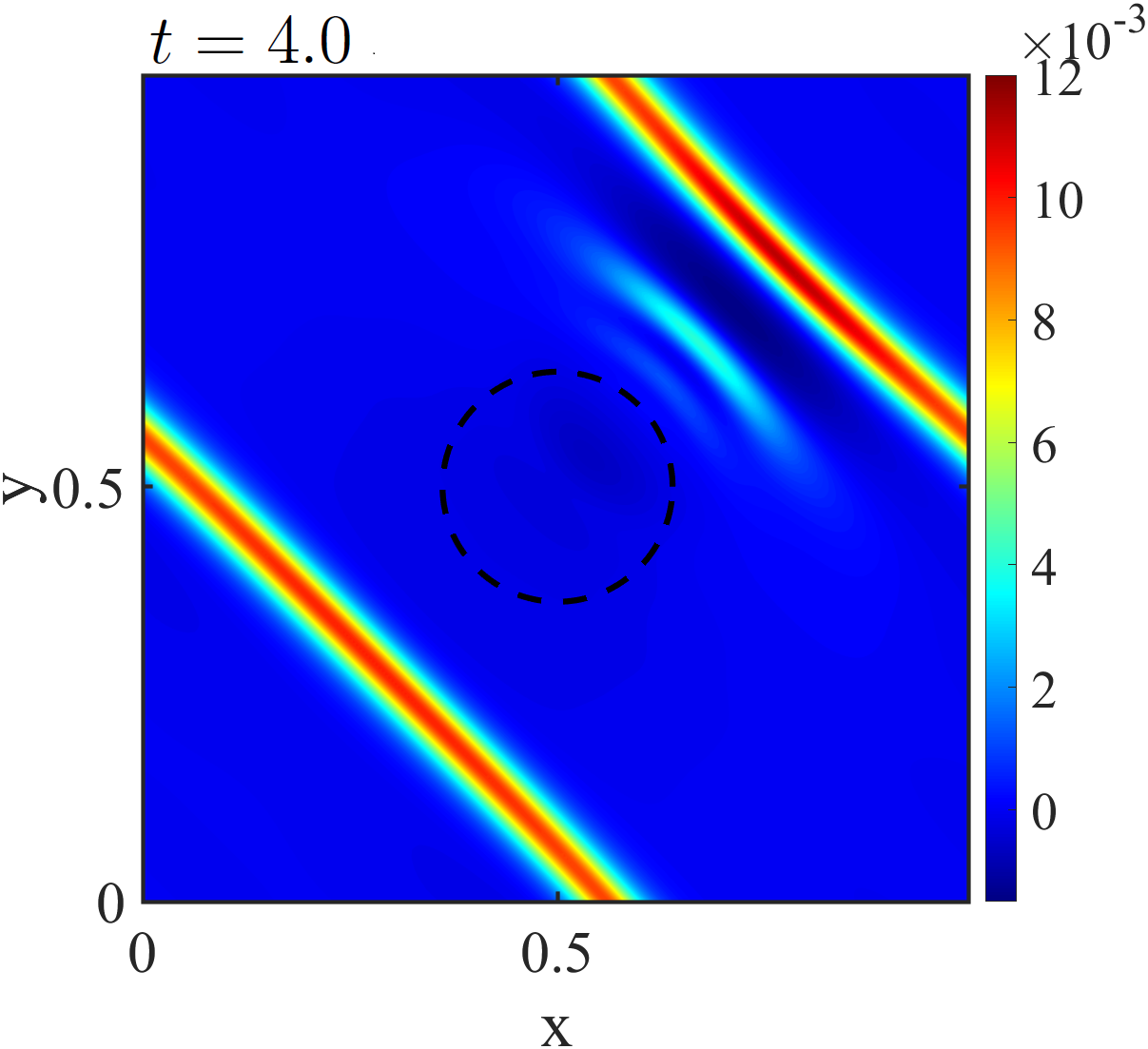}\\
    \caption{Contour plots of $\zeta$ showing a solitary wave passing over an underwater circular 2d barrier given by \eqref{Eq:2dBarrier}. The dashed line shows the level surface $b(x,y) = \frac{1}{2} 0.75 a$. Snapshots (top-bottom, left-right) correspond to $t=0, .8, 1.6, 2.4, 3.2, 4.0$.}\label{Fig:2d_simulation}
\end{figure}

Figure \ref{Fig:2d_simulation} shows the wave height solution $\zeta(x,y,t)$ at various times using AB2 with PCG. A $255 \times 256$ grid and time step $\Delta t = \frac{1}{3} \Delta x$ were used in the computation. Other time integration schemes produce visibly similar results on the scale presented, but yield a greater accuracy. 

The semi-implicit SBDF2 approach makes use of a single set of coefficients $\sigma^{\star}, \alpha^{\star}$ over the entire computation. To ensure that at each time the generalized eigenvalues of $\mA$ and $\mG_{\eta,h}$ lie in the stability region discussed in \S~\ref{subsec:zerostability}, we use the values of $\sigma^{\star}$ and $\alpha^{\star}$ given by the initial data. These values are close to the largest values provided by the formulas $\sigma^{\star}, \alpha^{\star}$ over the duration of the simulation. In general, if $\eta^{\star}(x,t)$ is the exact water depth solution of the SGN over the solution time $0 \leq t \leq T$, then the coefficients $\sigma^{\star}$ and $\alpha^{\star}$ in $\mA$ should be taken as the largest values, maximized over time, provided by \eqref{Eq:FlatbathCoeff}--\eqref{Eq:Optcoeff}.

\section{Discussion and Conclusions} 
In this work we have introduced a constant coefficient operator $\mA$ as a means to devise time integration approaches for the SGN equations that require only matrix-vector products involving $\mG_{\eta,h}$, i.e., are \emph{matrix-free}. The coefficients are rooted in theory, are shown to be quasi-optimal (within a class of constant coefficient operators) and enable mesh-independent performance of preconditioned conjugate gradient methods.  Furthermore, the operator $\mA$ enables the use of linearly implicit time integration that reduces the total number of linear solves involving $\mA$ to one per time step.  The computational theory is further supported by a suite of test cases.

Solving linear systems involving the constant coefficient operator $\mA$ is advantageous over direct solves involving $\mG_{\eta, h}$ since the latter is variable coefficient in both space and time. The formulas for $\mA$ and bounds provided here may also provide avenues for other efficient computational approaches. These potentially include successive convolution where the analytic and fast application of the Green's function for $\mA$ can be applied; exponential integrators since $\exp(\mA)$ may be readily evaluated using the Fourier transform, or boundary integral methods.  Such approaches could be particularly useful in the context of large scale computations (e.g., modeling the ocean) as the Green's function for $\mA$ can be analytically approximated and locally applied thereby avoiding the need to solve large linear systems. This would help to enable the use of physical models that incorporate dispersion for the simulation of waves in oceanography. 

\section*{Acknowledgments}
    The authors would like to thank Benjamin Seibold and Donna Calhoun for helpful conversations.  
    This material is based upon work support by the National Science Foundation under Grants No. DMS--1719693 (LF, DS); DMS--2012268, DMS--2309727 (DS); DMS-2108524 (WC). Any opinions, findings, and conclusions or recommendations expressed in this material are those of the authors and do not necessarily reflect the views of the National Science Foundation. 
    
\appendix
\section{Proof of Theorem~\ref{Thm:CondNumBnd}}\label{Appendix:Pf_MainResult}
This subsection provides a theoretical justification for the coefficients \eqref{Eq:FlatbathCoeff}--\eqref{Eq:Optcoeffb} by proving the convergence bound \eqref{Eq:CondNumBound}.  In the proof it will be useful to write both $A[\cdot,\cdot]$ and $B[\cdot,\cdot]$ into block matrix form as
\begin{align}\label{Eq:A_Blockform}
    A[\vu, \vv] &= \int_{\Omega}  \begin{pmatrix}
            \vu^T \; \; \nabla \cdot \vu
        \end{pmatrix} 
        \begin{pmatrix}
            \sigma^{\star} \mat{I} & 0 \\
            0 & \alpha^{\star}
        \end{pmatrix}
        \begin{pmatrix}
            \vv \\
            \nabla \cdot \vv
        \end{pmatrix} \, \du \vx  \, , \\  \label{Eq:B_Blockform}
    B[\vu, \vv] &= \int_{\Omega} 
        \begin{pmatrix}
            \vu^T \; \; \nabla \cdot \vu
        \end{pmatrix} \mat{M}(\vx) 
        \begin{pmatrix}
            \vv \\
            \nabla \cdot \vv
        \end{pmatrix} \, \du \vx  \, ,
\end{align}
where $\mat{M}(\vx)$ is the variable coefficient matrix
\begin{align}\label{Def:MatrixM}
    \mat{M}(\vx) = 
        \eta(\vx) \begin{pmatrix}
            \mat{I} + \nabla h(\vx) \, \nabla h^T(\vx) & \tfrac{1}{2} \eta(\vx) \, \nabla h(\vx) \\
            \tfrac{1}{2} \eta(\vx) \, \nabla h(\vx)^T & \tfrac{1}{3} \eta^2(\vx)
        \end{pmatrix}\, .
\end{align}

The following is a proof of the main result.

\begin{proof}[Of Theorem~\ref{Thm:CondNumBnd}] 
    Step 1: We first show that $\mat{M}$ and the coefficients \eqref{Eq:FlatbathCoeff}--\eqref{Eq:Optcoeff} satisfy the following matrix inequality
    \begin{align*}     
        \mathrm{(MI)} \quad \quad         
        (\gb)^{-1} \begin{pmatrix}
            \sigma^{\star} \mat{I} & 0 \\
            0 & \alpha^{\star}
        \end{pmatrix} \preceq     
    \mat{M}(\vx) \preceq
    \begin{pmatrix}
            \sigma^{\star} \mat{I} & 0 \\
            0 & \alpha^{\star} 
        \end{pmatrix}      
        \, \quad \textrm{holds} \; \textrm{for}\, \vx \in \Omega \, . 
\end{align*}
Here $\mat{I}$ is the $2\times 2$ identity matrix, $\mat{M}(\vx)$ is defined in \eqref{Def:MatrixM}, and for square symmetric matrices $\mat{X} \preceq \mat{Y}$ means $\mat{Y} - \mat{X}$ is positive semidefinite.  

When $|\nabla h|_{\rm max} = 0$, corresponding to a flat bathymetry, the matrix $\mat{M}$ is diagonal so that
\begin{align}\label{Eq:FlatBathMatrix}
    \begin{pmatrix}
        \eta_{\rm min} \mat{I} & 0 \\
        0 & \tfrac{1}{3}\eta_{\rm min}^3 
    \end{pmatrix} \preceq \begin{pmatrix}
        \eta(\vec{x}) \mat{I} & 0 \\
        0 & \tfrac{1}{3} \eta^3(\vec{x})  
    \end{pmatrix} \preceq 
    \begin{pmatrix}
        \eta_{\rm max} \mat{I} & 0 \\
        0 & \tfrac{1}{3} \eta_{\rm max}^3 
    \end{pmatrix} = \begin{pmatrix}
        \sigma^{\star} \mat{I} & 0 \\
        0 & \alpha^{\star}
    \end{pmatrix} \, . 
\end{align}
Combining \eqref{Eq:FlatBathMatrix} with \eqref{Eq:FlatbathCoeff} then yields (MI); the lower bound on $\mat{M}$ follows since $\eta_{\rm max} \geq \eta_{\rm min} > 0$. 

Assume now that  $|\nabla h|_{\rm max} \neq 0$. Let $\mat{R}(\vx)$ be the  Householder matrix that reflects $\nabla h(\vx)$ to align with the second coordinate axis  (if $\nabla h(\vec{x}) = 0$, take $\mat{R} = \mat{I}$), i.e.,  
\begin{align*}
    \mat{W}(\vx) = \begin{pmatrix}
        \mat{R}(\vx)^T & 0 \\
        0 & 1    
    \end{pmatrix}  \qquad 
    \mat{R}^T(\vx) \, \nabla h(\vx) = \begin{pmatrix}
        0 \\
         |\nabla h(\vx)|
    \end{pmatrix}  \qquad 
    \mat{R}(\vx)^T \mat{R}(\vx) = \mat{I}.
\end{align*}
With this $\mat{W}(\vx)$ 
\begin{align}\label{Eq:RotateM}
    \mat{W}(\vx)^T \mat{M}(\vx) \mat{W}(\vx) 
    = \begin{pmatrix}
        \eta(\vx) & 0 \\
        0 & \widetilde{\mat{M}}(\vx)
    \end{pmatrix} \, , 
\end{align}
where
\begin{align}\label{Def:Mtilde}    
    \widetilde{\mat{M}}(\vx) = 
    \begin{pmatrix}
            \eta(\vx) & 0 \\
            0 & 0 
        \end{pmatrix} +
        \eta(\vx)  
        \begin{pmatrix}
            |\nabla h (\vx)| & 0 \\
            0 & \eta(\vx) 
        \end{pmatrix}
        \begin{pmatrix}
            1 & \tfrac{1}{2} \\
            \tfrac{1}{2} & \tfrac{1}{3}
        \end{pmatrix}
        \begin{pmatrix}
            |\nabla h (\vx)| & 0  \\
            0 & \eta(\vx) 
        \end{pmatrix} \, . 
\end{align}
Moreover, from basic properties of semidefinite matrices, (MI) is equivalent to 
\begin{align*}
    \textrm{(MI-R)} \quad\quad &
    (\gb)^{-1} \begin{pmatrix}
        \sigma^\star \mat{I} & 0 \\
        0 & \alpha^{\star}
    \end{pmatrix}
    \preceq \mat{W}^T(\vx) \mat{M}(\vx) \mat{W}(\vx) 
    \preceq 
    \begin{pmatrix}
        \sigma^\star \mat{I} & 0 \\
        0 & \alpha^{\star}
    \end{pmatrix}
     \quad    \textrm{for}\,  \vx \in \Omega \, . 
\end{align*}
The equivalence of (MI) and (MI-R) holds since (i) $\mat{A} \preceq \mat{B} \Leftrightarrow \mat{Q}^T \mat{A} \mat{Q} \preceq \mat{Q}^T \mat{B} \mat{Q}$ for any invertible $\mat{Q}$ and pair of symmetric matrices $\mat{A}, \mat{B}$; and (ii) the diagonal matrix $\textrm{diag}(\sigma^{\star}, \sigma^{\star}, \alpha^{\star})$ is invariant when multiplied by $\mat{W}^T(\vx)$ and $\mat{W}(\vx)$. Thus, it is sufficient to show that (MI-R) holds. 

The matrix of numbers in $\widetilde{\mat{M}}$ has eigenvalues $\lambda_{\pm}$, i.e., 
\begin{align}\label{Eq:NumberMatrix}
    \lambda_- \mat{I} \preceq \begin{pmatrix}
        1 & \tfrac{1}{2} \\
        \tfrac{1}{2} & \tfrac{1}{3} 
    \end{pmatrix} \preceq 
    \lambda_+ \mat{I}  \qquad \textrm{where} \qquad \lambda_{\pm} = \tfrac{1}{6}(4 \pm \sqrt{13}) > 0 \, . 
\end{align}
Combining \eqref{Eq:NumberMatrix} with the expression for $\widetilde{\mat{M}}$ yields the bounds
\begin{align}\label{Eq:BoundsMtilde}
    \begin{pmatrix}
        \eta(\vec{x}) + \lambda_- \eta(\vec{x})|\nabla h(\vec{x})|^{2} & 0 \\
        0 & \lambda_- \eta^3(\vec{x}) 
    \end{pmatrix}
    \preceq \widetilde{\mat{M}}(\vx) \preceq 
    \begin{pmatrix}
        \eta(\vec{x}) + \lambda_+ \eta(\vec{x})|\nabla h(\vec{x})|^{2} & 0 \\
        0 & \lambda_+ \eta^3(\vec{x}) 
    \end{pmatrix} \,.
\end{align}
Inserting \eqref{Eq:BoundsMtilde} into the expression for $\mat{M}(\vx)$ in \eqref{Eq:RotateM} and maximizing over $\vx$ yields an upper bound
\begin{align*}
    \mat{W}(\vec{x})^T \mat{M}(\vx) \mat{W}(\vx) &\preceq 
    \begin{pmatrix}
        \eta_{\rm max} &  &  \\
         & \max_{\vx \in \Omega}\left(\eta(\vec{x}) + \lambda_+ \eta(\vec{x})|\nabla h(\vec{x})|^{2}\right) &  \\
         &  & \lambda_+ \eta_{\rm max}^3
    \end{pmatrix} \, , \\
    &\preceq \textrm{diag}\big( \sigma^{\star}, \sigma^{\star}, \alpha^{\star}\big) \,.
\end{align*}
The second line follows since $\lambda_{+} \eta(\vx)|\nabla h(\vx)|^2 \geq 0$ implies $\eta_{\rm max} \leq  \sigma^{\star}$.  Therefore $\mat{M}(\vx)$ satisfies the right hand side of (MI-R) and thus (MI).

In a similar spirit, inserting \eqref{Eq:BoundsMtilde} into the expression for $\mat{M}(\vx)$ and minimizing over $\vx$ yields a lower bound

\begin{align*}
     \mat{W}(\vec{x})^T \mat{M}(\vx) \mat{W}(\vx) &\succeq 
     \begin{pmatrix}
        \eta_{\rm min} &  &  \\
         & \min_{\vx \in \Omega}(\eta(\vec{x}) + \lambda_- \eta(\vec{x})|\nabla h(\vec{x})|^{2}) &  \\
         &  & \lambda_- \eta_{\rm min}^3 
    \end{pmatrix} \, ,    \\
    &\succeq \textrm{diag}\big( \eta_{\rm min}, \, \eta_{\rm min}, \, \lambda_{-}\eta_{\rm min}^3 \big) \,, \\
    &\succeq (\gb)^{-1} \textrm{diag}\big( \sigma^{\star}, \, \sigma^{\star}, \, \alpha^{\star}\big) \,.
\end{align*}
In the last line we have used the fact that 
\begin{align*}
    (\gb)^{-1} = \min\left\{ \frac{\eta_{\rm min}}{\sigma^{\star}}, \frac{\lambda_{-}}{\lambda_{+}}\left( \frac{\eta_{\rm min}}{\eta_{\rm max}}\right)^3 \right\} \, ,
\end{align*}
which simultaneously implies $(\gb)^{-1} \sigma^{\star} \leq \eta_{\rm min}$ and  $(\gb)^{-1} \alpha^{\star} \leq \lambda_{-}\eta_{\rm min}^3$. Therefore, $\mat{M}(\vx)$ satisfies the lower bound in (MI) as well. 

Step 2: Next, from (MI) the bilinear forms $A[\cdot, \cdot], B[\cdot, \cdot]$ satisfy the inequality
\begin{align}\label{Eq:VarIneq}
        (\gb)^{-1} \, A[\vu, \vu] \leq B[\vu, \vu] \leq A[\vu,\vu] \, \quad \forall \vu \in H_{\rm div}(\Omega) \, .
\end{align}
For if $\vu \in H_{\rm div}(\Omega)$ then $(\vu^T, \nabla \cdot \vu)^T \in (L^2(\Omega))^3$. We can then contract both sides of the left inequality in (MI) with the vector $(\vu^T \;\; \nabla \cdot \vu)^T$ to yield
    \begin{align} \label{Eq:Ineq1}   
        \begin{pmatrix}
            \vu^T \; \; \nabla \cdot \vu
        \end{pmatrix}
        \mat{M}(\vx)
        \begin{pmatrix}
            \vu \\
            \nabla \cdot \vu
        \end{pmatrix} 
        \leq 
        \begin{pmatrix}
            \vu^T \; \; \nabla \cdot \vu
        \end{pmatrix}
        \begin{pmatrix}
            \sigma^{\star} \mat{I} & 0 \\
            0 & \alpha^{\star}
    \end{pmatrix} 
        \begin{pmatrix}
            \vu \\
            \nabla \cdot \vu
        \end{pmatrix} \, , 
    \end{align}
which holds for all $\vu \in H_{\rm div}(\Omega)$ and (almost every) $\vx \in \Omega$. Integrating \eqref{Eq:Ineq1} over $\Omega$ and combining with \eqref{Eq:A_Blockform}--\eqref{Eq:B_Blockform} yields $B[\vu,\vu] \leq A[\vu,\vu]$. The lower bound in \eqref{Eq:VarIneq} follows from an identical argument.

Finally, \eqref{Eq:VarIneq} proves the desired result since
\begin{itemize}
    \item The conditioning number satisfies $\kappa_{\mA}( \mG_{\eta,h}) \leq \gb$. This follows immediately from dividing \eqref{Eq:VarIneq} by $A[\vu,\vu]$ and applying the definitions in \eqref{Exp:RayleighQuotients}. 
    \item Every eigenvalue $\lambda$ of $\mA^{-1} \mG_{\eta,h}$ (equivalently generalized eigenvalue of \eqref{Eq:Weakform}) is real and satisfies $(\gb)^{-1} \leq \lambda \leq 1$. Set $\vv = \vu$ in \eqref{Eq:WeakGenEV}, where $(\vu, \lambda)$ is a generalized eigenvector and eigenvalue, and substitute the resulting $B[\vu, \vu] = \lambda A[\vu,\vu]$ into \eqref{Eq:VarIneq} to obtain:
    \begin{align}
        (\gb)^{-1} A[\vu, \vu] \leq \lambda A[\vu, \vu] \leq A[\vu, \vu] \, .
    \end{align}
    Dividing by $A[\vu, \vu] \neq 0$ yields the result.         
\end{itemize}
\end{proof}

\section{Quasi-optimality of the Conditioning Number Bound}\label{Appendix:QuasiOpt}
The proof of the main result relies on $(\sigma^{\star}, \alpha^{\star}, \gb)$ satisfying (MI). One might then wonder if there is a set of coefficients satisfying (MI) that yield a (fundamentally) tighter conditioning number bound on $\mA^{-1} \mG_{\eta, h}$? To answer this question, consider the optimal conditioning number bound obtained by coefficients satisfying (MI)
\begin{alignat*}{2}
    (P) \qquad &&\textrm{Minimize:} & \quad \gamma \\
                &&\textrm{Subject to:} & \quad 
                \mat{M}(\vx) \preceq
                \begin{pmatrix}
                    \sigma \mat{I} & 0 \\
                    0 & \alpha
                \end{pmatrix} 
            \preceq 
            \gamma \, 
            \mat{M}(\vx)
                \, \quad \textrm{holds for} \; \vx \in \Omega \, , 
\end{alignat*}
Note that (P) is a convex problem since the objective function $\gamma$ is linear and the constraints on $\sigma$ and $\alpha$ are an intersection of a family of linear matrix inequalities. Denote the solution to (P) as $\gamma^{\star}_{\rm P}$. 

The answer turns out to be that no other set of coefficients satisfying (MI) can avoid the scalings of $\nabla h$ and $\eta$ in $\gb$. 
\begin{theorem}\label{Thm:Quasiopt} For any functions $\eta, \nabla h \in L^{\infty}(\Omega)$ with $\eta_{\rm min} > 0$, the bound $\gb$ satisfies
\begin{align}\label{Eq:Quasiopt_Inequality}
    \gamma^{\star}_{\rm P} \leq \gb \leq \left(\tfrac{4 + \sqrt{13}}{4 - \sqrt{13}}\right) \gamma^{\star}_{\rm P} \, .
\end{align}
In the case of a flat bathymetry ($\nabla h = 0$) then $\gb = \gamma^{\star}_{\rm P}$. 
\end{theorem}
The inequality \eqref{Eq:Quasiopt_Inequality} shows that $\gb$ is always within a constant factor of $\gamma^{\star}_{\rm P}$. Thus the fundamental scaling of $\eta$ and $\nabla h$ in $\gb$ cannot be avoided via other coefficients satisfying (MI).

To establish the result, consider a fixed pair $\eta$ and $\nabla h$ and introduce the feasible constraint set in (P) as
\begin{align}\label{Def:S}
    \mathcal{S} &\equiv \left\{ (\sigma, \alpha, \gamma) \in \mathbb{R}^3  \; \Big| \; 
     \mat{M}(\vx) \preceq
    \begin{pmatrix}
            \sigma \mat{I} & 0 \\
            0 & \alpha
        \end{pmatrix} 
    \preceq 
    \gamma \, 
    \mat{M}(\vx)    
    \right\} \,, \\
    &= \left\{ (\sigma,\alpha, \gamma) \in \mathbb{R}^3 \; \Big| \; 
    \eta(\vx) \leq \sigma \leq \gamma \, \eta(\vx), 
    \quad 
    \widetilde{\mat{M}}(\vx) \preceq
    \begin{pmatrix}
        \sigma & 0 \\
        0 & \alpha 
    \end{pmatrix} \preceq
    \gamma \, \widetilde{\mat{M}}(\vx) 
    \quad \forall \; \vx \in \Omega    
    \right\}\, .
\end{align}
Here $\widetilde{M}(\vx)$ is defined in \eqref{Def:Mtilde} and the second line follows from the same Householder reflection in Appendix~\ref{Appendix:Pf_MainResult}. The set $\mathcal{S}$ is also non-empty since $\eta_{\rm max}, |\nabla h|_{\rm max}$ are assumed to be finite. 

\medskip
\noindent\fbox{Case 1: A flat bathymetry $(\nabla h = 0)$}  When $\nabla h = 0$ the coefficients $(\sigma^\star, \alpha^\star, \gb)$ solve (P).

In this case the matrix inequality in $\mathcal{S}$ reduces to two linear inequalities defined by the diagonals of $\mat{M}(\vx)$,
\begin{alignat}{3}\nonumber 
    \mathcal{S} &=& \Big\{ (\sigma, \alpha, \gamma) \in \mathbb{R}^3 \; \big| \; &\eta(\vx) \leq \sigma \leq \gamma \; \eta(\vx) \quad &&\textrm{and} \quad
    \tfrac{1}{3}\eta^3(\vx) \leq \alpha \leq \gamma \; \tfrac{1}{3}\eta^3(\vx) \, \quad \forall x \in \Omega \Big\} \, ,  \\ \label{Eq:CompWiseFlatBottom}
    &=& \Big\{ (\sigma, \alpha, \gamma) \in \mathbb{R}^3 \; \big| \; &\eta_{\rm max} \leq \sigma \leq \gamma \, \eta_{\rm \min} \quad &&\textrm{and} \quad \tfrac{1}{3} \eta^3_{\rm max} \leq \alpha \leq \gamma \, \tfrac{1}{3} \eta^3_{\rm min} \Big\} \, . 
\end{alignat}
The second line follows from taking the (essential) supremum and infimum over $\vec{x}$. Thus, $\mathcal{S}$ is polyhedral---it consists of two pairs of linear inequalities.

The coefficients $(\sigma^\star, \alpha^\star, \gb)$ for $\mA$ then exactly solve (P), since 
\begin{align*}
    \sigma^\star = \eta_{\rm max} \, , \quad \alpha^\star = \tfrac{1}{3} \eta^3_{\rm max} \, , \quad \textrm{and} \quad 
    \gamma^{\star}_{\rm P} = \max\left\{ 
        \frac{\eta_{\rm max}}{\eta_{\rm min}}, 
        \left(\frac{\eta_{\rm max}}{\eta_{\rm min} } \right)^{3}
    \right\}     
    =
    \left(\frac{\eta_{\rm max}}{\eta_{\rm min}} \right)^{3} = \gb\, .
\end{align*}
While the choice of $\alpha^\star$ is unique, any value $\sigma \in [\eta_{\rm max}, \gb \eta_{\rm min}]$ yields the optimal $\gb$ value.

\medskip
\noindent\fbox{Case 2: A variable bathymetry $(\nabla h \neq 0)$} 
Here the set $\mathcal{S}$ is no longer polyhedral, and the coefficients $(\sigma^\star, \alpha^\star, \gb)$ for $\mA$ will be quasi-optimal solutions of (P). 

We first bound $\mathcal{S}$ from the inside and outside with polyhedral sets. 
\begin{proposition}[Polyhedral bounds on $\mathcal{S}$]\label{Prop:PolyhedralApprx}  
    For a fixed $\eta(\vx)$ and $\nabla h(\vx)$, we have:
    \begin{enumerate}
        \item[(a)] (Outer polyhedral bound of $\mathcal{S}$) Any $(\sigma, \alpha, \gamma) \in \mathcal{S}$ necessarily satisfies
        \begin{equation}\label{Eq:OuterApprx}
            \left. 
            \begin{aligned}
                \max_{\vx \in \Omega} \, \left( \eta(\vx) + \eta(\vx)|\nabla h(\vx) |^2 \right) &\leq 
                \sigma \leq \gamma \, \eta_{\rm min}   \\
                \tfrac{1}{3}\eta^3_{\rm max} &\leq \alpha \leq \gamma \, \tfrac{1}{3} \eta^3_{\rm min} 
            \end{aligned} \;  \right\} \, . 
        \end{equation}
        \item[(b)] (Inner polyhedral bound of $\mathcal{S}$) If $(\sigma, \alpha, \gamma)$ satisfies
        \begin{equation}\label{Eq:InnerAprx}
            \left. 
            \begin{aligned}
                \max_{\vx \in \Omega} \left( \eta(\vx) + \lambda_+ \eta(\vx) |\nabla h(\vx)|^2 \right) &\leq \sigma \leq \gamma \, \eta_{\rm min}  \\
                \lambda_+ \eta_{\rm max}^3 & \leq \alpha \leq \gamma \, \lambda_{-} \, \eta_{\rm min}^3 
            \end{aligned} \;  \right\} \, ,
        \end{equation}
        then $(\sigma, \alpha, \gamma) \in \mathcal{S}$.
        Here $\lambda_{\pm}$ are eigenvalues defined in \eqref{Eq:NumberMatrix}.                       
\end{enumerate}
\end{proposition}

The main idea of the proof is to bound $\widetilde{\mat{M}}(\vx)$ in \eqref{Def:Mtilde} with diagonal matrices.
\begin{proof}(Of Proposition~\ref{Prop:PolyhedralApprx}) 
\medskip
\noindent
(a) Let $(\sigma, \alpha, \gamma) \in \mathcal{S}$. Using the fact that semidefinite matrices have non-negative diagonal entries, comparing the diagonals in the matrix inequality for $\widetilde{\mat{M}}(\vx)$ in \eqref{Def:Mtilde}, shows that $(\sigma, \alpha, \gamma)$ satisfy: 
    \begin{equation}\label{Eq:Out_Inequalty}
        \left. \begin{aligned}
        \eta(\vx) \leq \sigma & \leq \gamma \, \eta(\vx)  \\ 
        \eta(\vx) + \eta(\vx) |\nabla h(\vx)|^2 \leq \sigma &\leq  \gamma \, \left( \eta(\vx) + \eta(\vx) |\nabla h(\vx)|^2 \right) \\ 
        \tfrac{1}{3} \eta^3(\vx) \leq \alpha &\leq  \gamma \, \tfrac{1}{3} \eta^3(\vx)  
        \end{aligned} \right\} \quad \forall \; \vx \in \Omega \, , 
    \end{equation} 
    Note that the first inequality in \eqref{Eq:Out_Inequalty} imposes a tighter upper bound on $\sigma$ while the second inequality in \eqref{Eq:Out_Inequalty} imposes a tighter lower bound on $\sigma$.  Maximizing and minimizing \eqref{Eq:Out_Inequalty} over $\vx$, and combining the two $\sigma$ inequalities together yields \eqref{Eq:OuterApprx}. 

\medskip 
\noindent
(b) Substituting the matrix inequality \eqref{Eq:NumberMatrix} into the definition of $\widetilde{\mat{M}}(\vx)$ shows that $\widetilde{\mat{M}}(\vx)$ (for almost every $\vx$) satisfies the following bound
\begin{align*}
    \begin{pmatrix}
        \eta(\vx) + \lambda_{-} \eta(\vx) |\nabla h(\vx)|^2 & 0  \\
        0 & \lambda_{-} \eta^3(\vx)
    \end{pmatrix} \, 
    &\preceq \widetilde{\mat{M}}(\vx) \preceq  
    \begin{pmatrix}
        \eta(\vx) + \lambda_{+} \eta(\vx) |\nabla h(\vx)|^2 & 0  \\
        0 & \lambda_{+} \eta^3(\vx)
    \end{pmatrix} \, .
\end{align*}
Thus, the matrix linear inequality in $\mathcal{S}$ for $\widetilde{\mat{M}}(\vx)$ can be ensured if $(\sigma, \alpha, \gamma)$ satisfy
\begin{align*}
    \begin{pmatrix}
        \sigma & 0\\
        0 & \alpha 
    \end{pmatrix} \preceq 
    \gamma \begin{pmatrix}
        \eta(\vx) + \lambda_{-} \eta(\vx) |\nabla h(\vx)|^2 & 0  \\
        0 & \lambda_{-} \eta^3(\vx)
    \end{pmatrix} \qquad \forall \, \vx \in \Omega \,,
\end{align*}
and
\begin{align*}    
    \begin{pmatrix}
        \eta(\vx) + \lambda_{+} \eta(\vx) |\nabla h(\vx)|^2 & 0  \\
        0 & \lambda_{+} \eta^3(\vx)
    \end{pmatrix} \preceq 
    \begin{pmatrix}
        \sigma & 0\\
        0 & \alpha 
    \end{pmatrix} \qquad \forall \, \vx \in \Omega \,.
\end{align*}
We thus observe that if $(\sigma, \alpha, \gamma)$ satisfy
\begin{equation}\label{Eq:SuffCondition}
    \left. \begin{aligned}
    \eta(\vx) &\leq \sigma \leq \gamma \, \eta(\vx)  \\
    \eta(\vx) + \lambda_{+} \eta(\vx) |\nabla h(\vx)|^2 & \leq \sigma \leq \gamma \, \left(\eta(\vx) + \lambda_{-} \eta(\vx) |\nabla h(\vx)|^2\right)  \\
    \lambda_{+} \eta^3(\vx) &\leq \alpha \leq \gamma \, \lambda_{-} \, \eta^3(\vx) 
    \end{aligned} \; \right\} \quad \forall \; \vx \in \Omega \, ,  
\end{equation}
then  $(\sigma, \alpha, \gamma) \in \mathcal{S}$.  Just as with, \eqref{Eq:Out_Inequalty}, the first inequality in \eqref{Eq:SuffCondition} imposes a tighter upper bound on $\sigma$ while the second imposes a tighter lower bound.  Maximizing and minimizing \eqref{Eq:SuffCondition} over $\vx$, and combining the two $\sigma$ inequalities together yields \eqref{Eq:InnerAprx}. 
\end{proof}

The coefficients for $\mA$ are in fact not only feasible, but yield the tightest conditioning bounds provided by the inner approximation
    \begin{align}\label{Eq:Optcoeff1}        
        \sigma^{\star}\equiv \, \max_{\vx \in \Omega} \, \left( \eta(\vx) + \lambda_{+} \eta(\vx)|\nabla h(\vx) |^2\right) \, , \quad  \quad 
        \alpha^{\star} \equiv \lambda_{+} \, \eta_{\rm max}^3 \, , \quad 
    \end{align}    
    along with 
    \begin{align}\label{Eq:DefGammaBnd}
        \gb \equiv \, \max\left\{ \frac{\sigma^{\star}}{\eta_{\rm min}},\frac{\lambda_{+}}{\lambda_{-}} \left(\frac{\eta_{\rm max}}{\eta_{\rm min}} \right)^{3} \right\} \, .
    \end{align}    

   \begin{proof}[Of Theorem~\ref{Thm:Quasiopt}]
       The lower bound in \eqref{Eq:Quasiopt_Inequality} follows immediately since $(\sigma^{\star}, \alpha^{\star}, \gb)$ is a feasible point of (P) and $\gamma^{\star}_{\rm P}$ is the global optimum of (P). 

       The upper bound in \eqref{Eq:Quasiopt_Inequality} is a consequence of the polyhedral bound \eqref{Eq:OuterApprx} in Proposition~\ref{Prop:PolyhedralApprx}(a).  Namely, every set of points in $\mathcal{S}$, including the optimal value $\gamma^{\star}_{\rm P}$ satisfies 
       \begin{align}\label{Def:BigGamma}
                \gamma^{\star}_{\rm P} \geq \Gamma \equiv \max\left\{ \frac{1}{{\eta_{\rm min}}} \max_{\vx \in \Omega} \, \left( \eta(\vx) + \eta(\vx)|\nabla h(\vx) |^2 \right), \left(\frac{\eta_{\rm max}}{\eta_{\rm min}} \right)^{3} \right\} \, ,
        \end{align}
        Here $\Gamma$ is the larger of the two lower bounds on $\gamma$ provided by \eqref{Eq:OuterApprx}. 

        Next, note that we can bound $\gb$ from above by a constant factor of $\Gamma$. Since $\lambda_+ > 1 > \lambda_- > 0$ and $\eta$ is non-negative, we have
        \begin{align*}
            \sigma^{\star} \leq \lambda_+ \max_{\vx \in \Omega} \, \left( \eta(\vx) + \eta(\vx)|\nabla h(\vx) |^2 \right) 
             \leq \frac{\lambda_+}{\lambda_-} \max_{\vx \in \Omega} \, \left( \eta(\vx) + \eta(\vx)|\nabla h(\vx) |^2 \right) \, .
        \end{align*}
        Hence, 
        \begin{align*}
            \gb 
            &\leq  \max\left\{ \frac{\lambda_+}{\lambda_{-} \eta_{\rm min}}\max_{\vx \in \Omega} \, \left( \eta(\vx) + \eta(\vx)|\nabla h(\vx) |^2 \right) , \frac{\lambda_{+}}{\lambda_{-}} \frac{\eta_{\rm max}^3}{\eta_{\rm min}^3} \right\} \\
            &= \frac{\lambda_+}{\lambda_-}\Gamma \,.
        \end{align*}
        Combining this with \eqref{Def:BigGamma} yields the upper bound in \eqref{Eq:Quasiopt_Inequality}. 
   \end{proof}
     
    \begin{remark}[Alternatives to $\lambda_{\pm}$]\label{Rmk:OptCoeff}
        The choice of eigenvalues $\lambda_{\pm}$ for upper and lower bounds in \eqref{Eq:NumberMatrix} is not unique. In fact, the following diagonal matrix bound
        \begin{align}\label{Def:CoefMatIneq}
            \frac{1}{2(2 + \sqrt{3})} \begin{pmatrix}
                1 & 0 \\
                0 & \tfrac{1}{3} 
            \end{pmatrix} \preceq
            \begin{pmatrix}
                1 & \tfrac{1}{2} \\
                \tfrac{1}{2} & \tfrac{1}{3}
            \end{pmatrix}
            \preceq     
            \tfrac{1}{2}(2 + \sqrt{3}) \begin{pmatrix}
                1 & 0 \\
                0 & \tfrac{1}{3} 
            \end{pmatrix}             
        \end{align}         
        can be used to obtain (slightly) tighter inner and outer polyhedral bounds on $\mathcal{S}$. The bound \eqref{Def:CoefMatIneq} can also be used to obtain a set of $(\sigma, \alpha, \gamma) \in \mathcal{S}$ with a (slightly) improved guarantee $\gamma \leq (2+\sqrt{3})^2 \gamma^{\star}_{\rm P}$. However in practice, we found that there was no practical improvement using these coefficients.  \myremarkend
    \end{remark}

\renewcommand\refname{Reference}
\bibliographystyle{siam}
\bibliography{references_DSWE}

\begin{thebibliography}{10}

\bibitem{SamaniegoLannes2008}
{\sc B.~Alvarez-Samaniego and D.~Lannes}, {\em Large time existence for 3d
  water-waves and asymptotics}, Inventiones mathematicae, 171 (2008),
  pp.~485--541.

\bibitem{AscherRuuthWetton1995}
{\sc U.~Ascher, S.~J. Ruuth, and B.~Wetton}, {\em Implicit-explicit methods for
  time dependent partial differential equations}, SIAM J. Numer. Anal., 32
  (1995), pp.~797--823.

\bibitem{BonnetonChazelLannesMarcheTissier2011}
{\sc P.~Bonneton, F.~Chazel, D.~Lannes, F.~Marche, and M.~Tissier}, {\em A
  splitting approach for the fully nonlinear and weakly dispersive
  {G}reen-{N}aghdi models}, J. Comput. Phys., 230 (2011), pp.~1479--1498.

\bibitem{Boussinesq1872}
{\sc J.~Boussinesq}, {\em Th{\'e}orie des ondes et des remous qui se propagent
  le long d'un canal rectangulaire horizontal, en communiquant au liquide
  contenu dans ce canal des vitesses sensiblement pareilles de la surface au
  fond}, Journal de Math{\'e}matiques Pures et Appliqu{\'e}es. Deuxi{\`e}me
  S{\'e}rie, 17 (1872), pp.~55--108.

\bibitem{BustoDumbserEscalanteFavrieGavrilyuk2021}
{\sc S.~Busto, M.~Dumbser, C.~Escalante, N.~Favrie, and S.~Gavrilyuk}, {\em On
  high order {ADER} discontinuous galerkin schemes for first order hyperbolic
  reformulations of nonlinear dispersive systems}, J. Sci. Comput., 87 (2021),
  p.~48.

\bibitem{c22}
{\sc W.~Choi}, {\em High-order strongly nonlinear long wave approximation and
  solitary wave solution}, J. Fluid Mech., 945 (2022), p.~A15.

\bibitem{DuranMarche2017}
{\sc A.~Duran and F.~Marche}, {\em A discontinuous galerkin method for a new
  class of {Green-Naghdi} equations on simplicial unstructured meshes}, Applied
  Mathematical Modelling, 45 (2017), pp.~840--864.

\bibitem{DutykhClamondMilewskiMitsotakis2013}
{\sc D.~Dutykh, D.~Clamond, P.~Milewski, and D.~Mitsotakis}, {\em Finite volume
  and pseudo-spectral schemes for the fully nonlinear 1d {S}erre equations},
  European Journal of Applied Mathematics, 24 (2013), pp.~761--787.

\bibitem{FavrieGavrilyuk2017}
{\sc N.~Favrie and S.~Gavrilyuk}, {\em A rapid numerical method for solving
  {Serre-Green-Naghdi} equations describing long free surface gravity waves},
  Nonlinearity, 30 (2017), pp.~2718--2736.

\bibitem{GavrilyukShyue2023}
{\sc S.~Gavrilyuk and K.-M. Shyue}, {\em {2D Serre-Green-Naghdi} equations over
  topography: {E}lliptic operator inversion method}, J. Hydralic Engineering,
  150 (2023).

\bibitem{GillesTownsend2019}
{\sc M.~Gilles and A.~Townsend}, {\em Continuous analogues of {K}rylov subspace
  methods for differential operators}, SIAM J. Numer. Anal., 57 (2019),
  pp.~899--924.

\bibitem{gn76}
{\sc A.~E. Green and P.~M. Naghdi}, {\em Derivation of equations for wave
  propagation in water of variable depth}, J. of Fluid Mech., 78 (1976),
  pp.~237--246.

\bibitem{GuermondKeesPopovTovar2019}
{\sc J.-L. Guermond, C.~Kees, B.~Popov, and E.~Tovar}, {\em Robust explicit
  relaxation technique for solving the {G}reen-{N}aghdi equations}, J. Comput.
  Phys., 399 (2019), p.~108917.

\bibitem{GuermondKeesPopovTovar2022b}
\leavevmode\vrule height 2pt depth -1.6pt width 23pt, {\em Hyperbolic
  relaxation technique for solving the dispersive the
  {S}erre–{G}reen-{N}aghdi equations with topography}, J. Comput. Phys., 450
  (2022), p.~110809.

\bibitem{GuermondKeesPopovTovar2022a}
\leavevmode\vrule height 2pt depth -1.6pt width 23pt, {\em Well-balanced
  second-order convex limiting technique for solving the
  {S}erre–{G}reen-{N}aghdi equations}, Water Waves, 4 (2022), pp.~409--445.

\bibitem{wanner1996solving}
{\sc E.~Hairer and G.~Wanner}, {\em Solving ordinary differential equations
  II}, Springer Berlin Heidelberg, 1996.

\bibitem{JangSungPark2024}
{\sc T.~Jang, H.~Sung, and J.~Park}, {\em A non-local formulation for
  simulating the fully nonlinear {Serre-Green-Naghdi} equations for a solitary
  wave interaction with a variable slope}, Applied Ocean Research, 153 (2024),
  p.~104220.

\bibitem{Kirby2010}
{\sc R.~C. Kirby}, {\em From functional analysis to iterative methods}, SIAM
  Review, 52 (2010), pp.~269--293.

\bibitem{kdv95}
{\sc D.~J. Korteweg and G.~de~Vries}, {\em On the change of form of long waves
  advancing in a rectangular canal, and on a new type of long stationary
  waves}, Phil. Mag, 39 (1895), pp.~422--443.

\bibitem{LannesBonneton2009}
{\sc D.~Lannes and P.~Bonneton}, {\em Derivation of asymptotic two-dimensional
  time-dependent equations for surface water wave propagation}, Phys. Fluids,
  21 (2009), p.~016601.

\bibitem{LannesMarche2015}
{\sc D.~Lannes and F.~Marche}, {\em A new class of fully nonlinear and weakly
  dispersive {G}reen-{N}aghdi models for efficient 2{D} simulations}, J.
  Comput. Phys., 282 (2015), pp.~238--268.

\bibitem{MetayerGavrilyukHank2010}
{\sc O.~{Le Métayer}, S.~Gavrilyuk, and S.~Hank}, {\em A numerical scheme for
  the {G}reen–{N}aghdi model}, Journal of Computational Physics, 229 (2010),
  pp.~2034--2045.

\bibitem{Leveque2007}
{\sc R.~LeVeque}, {\em Finite Difference Methods for Ordinary and Partial
  Differential Equations: Steady-State and Time-Dependence Problems}, SIAM,
  first~ed., 2007.

\bibitem{LiHymanChoi2004}
{\sc Y.~A. Li, J.~M. Hyman, and W.~Choi}, {\em A numerical study of the exact
  evolution equations for surface waves in water of finite depth}, Stud. Appl.
  Math., 113 (2004), pp.~303--324.

\bibitem{MadsenMei1969}
{\sc O.~S. Madsen and C.~C. Mei}, {\em The transformation of a solitary wave
  over an uneven bottom}, J. Fluid Mech., 39 (1969), pp.~781--791.

\bibitem{mbl02}
{\sc P.~A. Madsen, H.~B. Bingham, and H.~Liu}, {\em A new approach to
  high-order {B}oussinesq models.}, J. Fluid Mech., 462 (2002), pp.~1--30.

\bibitem{ms98}
{\sc P.~A. Madsen and H.~A. Schaffer}, {\em Higher-order {B}oussinesq-type
  equations for surface gravity waves: derivation and analysis}, Proc. R. Soc.
  Lond. A, 356 (1998), pp.~3123--3184.

\bibitem{MalekStrakos2015}
{\sc J.~M\'{a}lek and Z.~Strak\u{o}s}, {\em Preconditioning and the Conjugate
  Gradient Method in the Context of Solving PDEs}, SIAM, 2015.

\bibitem{mat15}
{\sc Y.~Matsuno}, {\em Hamiltonian formulation of the extended {G}reen-{N}aghdi
  equations}, Physica D, 301-302 (2015), pp.~1--7.

\bibitem{mat16}
\leavevmode\vrule height 2pt depth -1.6pt width 23pt, {\em Hamiltonian
  structure for two-dimensional extended {G}reen-{N}aghdi equations}, Proc. R.
  Soc. A, 472 (2016), p.~20160127.

\bibitem{Monk2003}
{\sc P.~Monk}, {\em Finite element methods for {M}axwell's {E}quations}, Oxford
  University Press, 2003.

\bibitem{NoelleParisotTscherpel2022}
{\sc S.~Noelle, M.~Parisot, and T.~Tscherpel}, {\em A class of boundary
  conditions for the time-discrete {G}reen-{N}aghdi equations with bathymetry},
  SIAM J. Numer. Anal., 60 (2022), pp.~2681--2712.

\bibitem{nwo93}
{\sc O.~Nwogu}, {\em An alternative form of the {B}oussinesq equations for
  nearshore wave propagation.}, J. Water., Port, Coastal, Ocean Eng., 119
  (1993), pp.~618--638.

\bibitem{Parisot2019}
{\sc M.~Parisot}, {\em Entropy-satisfying scheme for a hierarchy of dispersive
  reduced models of free surface flow}, Internat. J. Numer. Methods Fluids, 91
  (2019), pp.~509--531.

\bibitem{Popinet2015}
{\sc S.~Popinet}, {\em A quadtree-adaptive multigrid solver for the
  {S}erre-{G}reen-{N}aghdi equations}, J. Comput. Phys., 302 (2015),
  pp.~336--358.

\bibitem{RanochaRicchiuto2024}
{\sc H.~Ranocha and M.~Ricchiuto}, {\em Structure-preserving approximations of
  the {Serre-Green-Naghdi} equations in standard and hyperbolic form}, 2024.
\newblock arxiv.org/pdf/2408.02665.

\bibitem{ray76}
{\sc J.~W.~S. Rayleigh}, {\em On waves}, Phil. Mag., Ser. 5, 1 (1876),
  pp.~257--279.

\bibitem{Theory}
{\sc R.~Rosales, B.~Seibold, D.~Shirokoff, and D.~Zhou}, {\em Unconditional
  stability for multistep {I}m{E}x schemes-theory}, SIAM Journal on Numerical
  Analysis, 55 (2017), pp.~2336--2360.

\bibitem{SamiiDawson2018}
{\sc A.~Samii and C.~Dawson}, {\em An explicit hybridized discontinuous
  galerkin method for {Serre-Green-Naghdi} wave model}, Computer Methods in
  Applied Mechanics and Engineering, 330 (2018), pp.~447--470.

\bibitem{Schechter2001}
{\sc M.~Schechter}, {\em Principles of Functional Analysis}, American
  Mathematical Society, second~ed., 2001.

\bibitem{Prac}
{\sc B.~Seibold, D.~Shirokoff, and D.~Zhou}, {\em Unconditional stability for
  multistep {I}m{E}x schemes-practice}, Journal of Computational Physics, 376
  (2019), pp.~295--321.

\bibitem{Serre1953}
{\sc F.~Serre}, {\em Contribution \`a l'\'etude des \'ecoulements permanents et
  variables dans les canaux}, La Houille Blanche, 39 (1953), pp.~830--872.

\bibitem{sg69}
{\sc C.~H. Su and C.~S. Gardner}, {\em Korteweg-de vries equation and
  generalizations. {III}: {D}erivation of the {K}orteweg-de {V}ries equation
  and {B}urgers equation}, J. Math. Phys., 10 (1969), pp.~536--539.

\bibitem{TrefethenBau1997}
{\sc L.~Trefethen and D.~Bau}, {\em Numerical Linear Algebra}, SIAM, first~ed.,
  1997.

\bibitem{wei95}
{\sc G.~Wei, J.~T. Kirby, S.~T. Grilli, and R.~Subramanya}, {\em A fully
  nonlinear {B}oussinesq model for surface waves. {P}art {I}. {H}ighly
  nonlinear unsteady waves}, J. Fluid Mech., 294 (1995), pp.~71--92.

\bibitem{wu99}
{\sc T.~Y. Wu}, {\em Modeling nonlinear dispersive water waves}, J. Eng. Mech.,
  11 (1999), pp.~747--755.

\end{thebibliography}


\end{document}